\documentclass[11pt,leqno]{article}

\usepackage{amsmath,amsfonts,amscd,amssymb,theorem}

\long\def\comment#1\endcomment{}

\comment
\endcomment

\makeatletter
\begingroup
\gdef\th@dotted{\normalfont\itshape
  \def\@begintheorem##1##2{%
        \item[\hskip\labelsep \theorem@headerfont ##1\ ##2.]}%
\def\@opargbegintheorem##1##2##3{%
   \item[\hskip\labelsep \theorem@headerfont ##1\ ##2\ (##3).]}}
\endgroup
\makeatother

\theoremstyle{dotted}

\newtheorem{theorem}{Theorem}[section]
\newtheorem{lemma}[theorem]{Lemma}

\newtheorem{prop}[theorem]{Proposition}
\newtheorem{corr}[theorem]{Corollary}

\makeatletter
\begingroup
\gdef\th@upshape{\normalfont
  \def\@begintheorem##1##2{%
        \item[\hskip\labelsep \theorem@headerfont ##1\ ##2.]}%
\def\@opargbegintheorem##1##2##3{%
   \item[\hskip\labelsep \theorem@headerfont ##1\ ##2\ (##3).]}}
\endgroup
\makeatother

\theoremstyle{upshape}

\newtheorem{defn}[theorem]{Definition}
\newtheorem{remark}[theorem]{Remark}
\newtheorem{exa}[theorem]{Example}

\makeatletter
\renewcommand{\subsection}{\@startsection{subsection}{2}{0pt}{-3ex
plus -1ex minus -0.2ex}{-2mm plus -0pt minus
-2pt}{\normalfont\bfseries}} 
\renewcommand{\subsubsection}{\@startsection{subsubsection}{3}{0pt}{-3ex
plus -1ex minus -0.2ex}{-2mm plus -0pt minus
-2pt}{\normalfont\bfseries}} 
\makeatother

\makeatletter
\@addtoreset{equation}{section}
\makeatother

\newcommand{\cntrct}                % contraction with a vector field
{\hspace{2pt}\raisebox{1pt}{\text{$\lrcorner$}}\hspace{2pt}}

\newcommand{\proof}[1][Proof.]{\smallskip\noindent{\em #1}}
\def\endproof{\hfill\ensuremath{\square}\par\medskip}

\def\eqref#1{\thetag{\ref{#1}}}

\let\latexref=\ref
\def\ref#1{{\normalfont{\latexref{#1}}}}

\newcommand{\wt}{\widetilde}
\newcommand{\wh}{\widehat}

\newcommand{\br}[1]{\langle #1 \rangle}

\newcommand{\dg}{\dagger}

\newcommand{\idot}{{\:\raisebox{1pt}{\text{\circle*{1.5}}}}}
\newcommand{\hdot}{{\:\raisebox{3pt}{\text{\circle*{1.5}}}}}
\newcommand{\eps}{\varepsilon}
\renewcommand{\phi}{\varphi}

\newcommand{\Spec}{\operatorname{Spec}}
\newcommand{\Fun}{\operatorname{Fun}}

\newcommand{\Hom}{\operatorname{Hom}}
\newcommand{\Ext}{\operatorname{Ext}}

\newcommand{\Coker}{\operatorname{Coker}}
\newcommand{\Ker}{\operatorname{Ker}}
\renewcommand{\Im}{\operatorname{Im}}

\newcommand{\gr}{\operatorname{\sf gr}}

\newcommand{\ppt}{{\sf pt}}
\newcommand{\id}{\operatorname{\sf id}}

\newcommand{\D}{\mathcal{D}}
\newcommand{\DF}{\mathcal{D}F}

\newcommand{\bCF}{\overline{CF}}

\newcommand{\Z}{\mathbb{Z}}
\newcommand{\Q}{\mathbb{Q}}

\newcommand{\F}{\mathbb{F}}

\newcommand{\K}{\mathbb{K}}

\newcommand{\lotimes}{\overset{\sf\scriptscriptstyle L}{\otimes}}

\newcommand{\C}{\mathcal{C}}

\newcommand{\I}{{\sf I}}

\newcommand{\vH}{\check{H}}
\newcommand{\vC}{\check{C}}
\newcommand{\vbC}{\check{\bar{C}}}
\newcommand{\wC}{\wt{C}}

\newcommand{\Tot}{\operatorname{\sf Tot}}
\newcommand{\tot}{\operatorname{\sf tot}}

\newcommand{\Aut}{\operatorname{\sf Aut}}

\newcommand{\Exp}{\operatorname{Exp}}

\newcommand{\Per}{\operatorname{Per}}
\newcommand{\cPer}{\operatorname{\overline{Per}}}
\newcommand{\per}{\operatorname{per}}

\newcommand{\bCP}{\overline{CP}}
\newcommand{\bHP}{\overline{HP}}

\newcommand{\wCC}{\wh{CC}}

\newcommand{\bCPH}{\overline{CPH}}

\newcommand{\Sets}{\operatorname{Sets}}

\newcommand{\bP}{\overline{P}}
\newcommand{\bQ}{\overline{Q}}
\newcommand{\bi}{\bar{\imath}}
\newcommand{\bLambda}{\overline{\Lambda}}

\newcommand{\bV}{W}%{{\overline{V}}}

\newcommand{\aug}{\operatorname{\sf aug}}

\newcommand{\vpi}{\check{\pi}}
\newcommand{\vbpi}{\check{\bar{\pi}}}

\newcommand{\bdelta}{\bar{\delta}}

\newcommand{\wj}{\wt{j}}
\newcommand{\bJ}{\bar{J}}

\newcommand{\hs}{\natural}

\def\dlim_#1{{\displaystyle\lim_{#1}}^\hdot}

\title{Co-periodic cyclic homology}

\author{D. Kaledin\thanks{Partially supported by RScF, grant number
    14-21-00053, and the Dynasty Foundation award.}}

\begin{document}

\maketitle

\tableofcontents

\section*{Introduction.}

Periodic cyclic homology $HP_\idot(A)$ of a unital associative flat
algebra $A$ over a commutative ring $k$ is defined as the homology
of an explicit complex --- or rather, a bicomplex
$CP_{\idot,\idot}(A)$. In one direction, the bicomplex is bounded:
one has $CP_{i,\idot}(A)=0$ for $i < 0$. However, there are no
restrictions in the the other direction. On the contrary,
$CP_{i,j}(A)$ is $2$-periodic with respect to the index
$j$. Therefore turning the bicomplex $CP_{\idot,\idot}(A)$ into a
complex requires a choice.

In the standard theory, what one considers is the product-total
complex of $CP_{\idot,\idot}(A)$. Traditionally, considering the
sum-total complex was assumed to be a wrong thing to do, since if
$k$ contains $\Q$, the sum-total complex is trivially acyclic (see
e.g.\ \cite[Section 5.1.2]{Lo}). Nevertheless, about ten years ago,
it has been suggested by M. Kontsevich (\cite{ko1},
\cite[2.32]{ko2}) that if $k$ is a field of positive characteristic
$p$, then taking the sum-total complex gives a meaningful and
interesting new theory.

\medskip

It seems that Kontsevich's suggestion was not taken seriously at the
time (in particular, I didn't follow it up in \cite{ka-1} and
\cite{ka0} --- in retrospect, a rather glaring omission). However,
recently a very similar phenomenon has appeared in the work of
A. Beilinson and B. Bhatt. They work in the commutative situation,
where periodic cyclic homology is known to be intimately related to
the de Rham cohomology. In particular, if $A$ is a commutative
algebra over a field $k$ of characteristic $0$ with a smooth
spectrum $X = \Spec A$, then $HP_\idot(A)$ is entirely expressible
in terms of de Rham cohomology groups $H^\hdot_{DR}(A) =
H^\hdot_{DR}(X)$. In the situations considered by Beilinson and
Bhatt, $k$ is a field $\F_q$ of positive characteristic $p$ (or more
generally, a truncated Witt vectors ring $W_n(\F_q)$, $n \geq 1$),
and $\Spec A$ is not smooth. So, instead of de Rham cohomology, one
works with {\em derived de Rham cohomology} introduced by L. Illusie
\cite{ill}. To define it for an affine variety $X = \Spec A$, one
replaces the algebra $A$ with a simplicial resolution $A_\idot$
whose terms are free commutative algebras, and takes de Rham complex
termwise. Then after passing to the usual standard complex of a
simplicial object, one obtains a bicomplex $\Omega^\hdot(A_\idot)$,
and here again, there is the issue of how to totalize it. In
Illusie's approach, one takes the product-total complex. If we
denote its cohomology simply by $H^\hdot_{DR}(A)$, then the standard
Hodge-to-de Rham spectral sequence induces a convergent spectral
sequence
$$
\Lambda^\hdot\Omega_\idot(A) \Rightarrow H^\hdot_{DR}(A),
$$
where $\Lambda^\hdot$ is an appropriate derived version of the
exterior power functors, and $\Omega_\idot(A)$ is the cotangent
complex of $A$. In particular, $H^\hdot_{DR}(A)$ does not depend on
the choice of a resolution $A_\idot$.

If one instead takes the sum-total complex of the bicomplex
$\Omega^\hdot(A_\idot)$ and denotes its homology by
$\overline{H}^\hdot_{DR}(A)$, then a different spectral sequence
becomes convergent --- namely, the spectral sequence
$$
H^\hdot_{DR}(A_\idot) \Rightarrow \overline{H}^\hdot_{DR}(A).
$$
If $k$ is a field of characteristic $0$, this shows that
$\overline{H}^\hdot_{DR}(A)$ is completely uninteresting. Indeed,
since the terms $A_i$ of the resolution $A_\idot$ are free algebras,
so that $\Spec A_i$ are just affine spaces, the de Rham cohomology
$H^\hdot_{DR}(A_i)$ is simply $k$ placed in degree $0$, irrespective
of $A_i$, and one deduces that $\overline{H}^\hdot_{DR}(A)$ is also
$k$ placed in degree $0$.

However, the situation changes drastically if $k$ is a field of
positive characteristic. In this case, affine spaces have lots of de
Rham cohomology --- in fact, one has the classic Cartier isomorphism
$$
H^i_{DR}(A_j) \cong \Omega^i(A_j^{(1)}),
$$
where ${(1)}$ indicates the Frobenius twist. With this in mind, what
we obtain is the derived version of the {\em conjugate spectral
  sequence}, and it reads as
$$
\Lambda^\hdot\Omega_\idot(A^{(1)}) \Rightarrow
\overline{H}^\hdot_{DR}(A).
$$
In particular, $\overline{H}^\hdot_{DR}(A)$ also does not depend on
the resolution $A_\idot$, and this justifies our notation. Beilinson
and Bhatt go on to apply this observation to some specific algebras
$A$ important in $p$-adic Hodge theory, where working with
$\overline{H}^\hdot_{DR}(-)$ instead of $H^\hdot_{DR}(-)$ allows one
to obtain much stronger results.

\medskip

The goal of the present paper is then the following. Motivated by
the work of Beilinson and Bhatt, we take up Kontsevich's original
suggestion and study ``wrong'' totalizations of the cyclic bicomplex
$CP_{\idot,\idot}(A)$. To make room for interesting applications, we
work with DG algebras $A_\idot$ instead of just associative
algebras. For any DG algebra $A_\idot$ over a commutative ring $k$,
we define a functorial {\em co-periodic cyclic complex}
$\bCP_\idot(A_\idot)$, and we call its homology
$\bHP_\idot(A_\idot)$ {\em co-periodic cyclic homology} of
$A_\idot$. If $A_\idot = A$ is concentrated in homological degree
$0$, then $\bCP_\idot(A)$ is literally the sum-total complex of the
periodic cyclic bicomplex $CP_{\idot,\idot}(A)$, and in the general
case, the definition requires only a minor modification. We then
extend the definition to small DG categories $A_\idot$, and we prove
the following:
\begin{itemize}
\item As soon as the base ring $k$ is Noetherian,
  $\bHP_\idot(A_\idot)$ is derived Morita-invariant, and moreover,
  gives an additive invariant of small DG categories in the sense of
  \cite{kel}.
\item If $k$ is a field of finite characteristic, and $A_\idot$ is a
  cohomologically bounded smooth DG algebra over $k$, then
  $\bHP_\idot(A_\idot) \cong HP_\idot(A_\idot)$.
\end{itemize}
Our actual statements, Theorem~\ref{morita.thm},
Theorem~\ref{loc.thm} and Theorem~\ref{comp.thm}, are slightly
stronger and more precise, but the above sums up the essential
points. We also prove the following somewhat more technical result.
\begin{itemize}
\item Assume that $k$ is a field of odd positive characteristic $p >
  2$. Then for any small DG category $A_\idot$, there exists a
  functorial {\em conjugate spectral sequence}
$$
HH_\idot(A_\idot)^{(1)}((u^{-1})) \Rightarrow \bHP_\idot(A_\idot),
$$
where $HH_\idot(A_\idot)$ is the Hochschild homology of $A_\idot$,
$(1)$ stands for the Frobenius twist, $u$ is a formal generator of
cohomological degree $2$, and $((u^{-1}))$ is the
shorthand for ``formal Laurent powers series in $u^{-1}$''.
\end{itemize}
The conjugate spectral sequence generalizes the usual commutative
spectral sequence in the same way that the Hodge-to-de Rham spectral
sequence is generalized to the ``Hochschild-to-cyclic'' spectral
sequence
$$
HH_\idot(A_\idot)((u)) \Rightarrow HP_\idot(A_\idot).
$$
Note that for the conjugate spectral sequence, we need to take power
series in $u^{-1}$ rather than $u$.

\medskip

The paper is organized as follows. Unfortunately, our definition of
co-periodic cyclic homology works by an explicit complex, and we
have not been able to find a more invariant treatment similar to
\cite{C}. We do rely on Connes' notion of a cyclic object and the
small category $\Lambda$ to encode all the relevant combinatoris,
but we are unable to use the machinery of derived categories. As a
next best choice, we equip everything in sight with a filtration,
and use filtered derived categories. This is still better than
writing down explicit maps of complexes, but it requires a lot of
preliminaries. These are contained in Section~\ref{filt.sec},
Section~\ref{mix.sec} and Section~\ref{cyc.sec}. Almost nothing in
these three sections is new, with a possible exception of some (but
not all) material in Subsection~\ref{t.subs},
Subsection~\ref{reso.subs} and Subsection~\ref{repr.subs}. However,
these are things that need to be spelled out in all the gory detail,
in order to fix notation and terminology and avoid
ambiguities. Section~\ref{filt.sec} is devoted to filtrations,
filtered derived categories, convergence of spectral sequences and
suchlike. The possibly new ingredient is a generalization of the
canonical truncation to filtered complexes given in
Subsection~\ref{t.subs}. Section~\ref{mix.sec} is devoted to mixed
complexes: we state the necessary definitions, discuss examples
arising from Tate cohomology of cyclic groups, and prove one simple
technical statement on what we call ``mixed resolutions'' that is
used later in the paper. Section~\ref{cyc.sec} contains the
definition of the co-periodic cyclic complex, and some immediate
observations on its properties. The input here is a cyclic object or
a complex of cyclic objects, so that the theory is essentialy
linear.

The technical heart of the paper are Section~\ref{edge.sec} and
Section~\ref{conj.sec}. We still work in the setting of complexes of
cyclic objects. In Section~\ref{edge.sec}, we start using our main
technical gadget, namely, a certain version $\Lambda_p$ of the
cyclic category $\Lambda$ that appeared in \cite{FT}. This is a
small category that comes equipped with two natural functors to
$\Lambda$ --- an ``edgewise subdivision'' functor $i:\Lambda_p \to
\Lambda$, and a projection $\pi:\Lambda_p \to \Lambda$ that
corresponds to the $p$-fold cover of a circle by itself. After
reminding classic results about edgewise subdivision and the functor
$i$, we give, in Proposition~\ref{edge.prop}, a filtered refinement
of the main classic result (this refinement seems to be new). Then
in Subsection~\ref{proj.subs}, we turn to the functor $\pi$ and
prove Lemma~\ref{proj.le}, a kind of a retraction statement that
allows to reduce questions about cyclic objects coming from the
projection $\pi$ to questions about their restriction to the
subcategory $\Delta^o \subset \Lambda$.

Section~\ref{conj.sec} is mostly devoted to the conjugate spectral
sequence: we start with a complex of cyclic objects, and we produce
a spectral sequence converging to its co-periodic cyclic
homology. We need to impose one condition on the complex that we
call ``tightness''; a version of it already appeared in \cite{ka1}
and before that, in \cite{ka-1}. Then in
Subsection~\ref{conve.subs}, we turn to the study of convergence of
this spectral sequence, and show that convergence questions for
cyclic objects can be reduced to questions about their restrictions
to $\Delta^o$ --- this uses Lemma~\ref{proj.le} in an essential way.

Finally, in Section~\ref{alg.sec}, we turn to our main subject ---
namely, DG algebras and DG categories. After giving a brief reminder
on DG algebras and DG categories, we define co-periodic cyclic
homology of a DG algebra and a small DG category, and prove our main
results. The main technical tool here is a very simple fact of
linear algebra --- namely, the computation of Tate cohomology of a
cyclic group $\Z/p\Z$ with coefficients in the $p$-th tensor power
$V^{\otimes_k p}$ of a vector space $V$ over a field $k$ of
characteristic $p$. This has already appeared in \cite{ka-1} and
more recently, in \cite{ka1}. In the present paper, in
Subsection~\ref{ten.subs}, we give a generalization of this
computation to complexes of vector spaces that seems to be
new. Together with the material of Section~\ref{conj.sec}, this
immediately gives the conjugate spectral sequence and yields
Theorem~\ref{morita.thm} and Theorem~\ref{loc.thm}. Then the same
computation is applied to a certain twisted version of Hochschild
homology, and together with the material of
Subsection~\ref{conve.subs}, allows us to prove Theorem~\ref{comp.thm},
our main comparison theorem. In the end, in
Subsection~\ref{derham.subs}, we explain briefly how our results
compare to results about de Rham cohomology, and in particular, to
results of Beilinson and Bhatt.

\subsection*{Acknowledgements.} I am very grateful to A. Beilinson,
B. Bhatt, A. Efimov, B. Feigin, V. Hinich, D. Kazhdan,
M. Kontsevich, A. Kuznetsov, N. Markarian, D. Orlov, D. Pomerleano,
and V. Vologodsky for useful remarks and interesting discussions on
the subject. I am especially grateful to Vadim Vologodsky for
convincing me that the main Tate cohomology computation should work
for complexes, and to David Kazhdan for encouraging me to revisit
the subject and to try to finally settle it.

\section{Filtrations.}\label{filt.sec}

\subsection{Filtered objects.}

Assume given an abelian category $\C$ satisfying $AB4$ and $AB4^*$
(that is, $\C$ has arbitrary products and sums, and both are
exact). In this paper, by a {\em filtration} $F^\hdot$ on an object
$E \in \C$ we will understand a collection of subobjects $F^jE
\subset E$, one for each integer $j$, such that $F^jE \subset F^iE$
whenever $j \geq i$. In other words, all our filtrations are
decreasing and indexed by all integers. A {\em filtered object} in
$\C$ is an object equipped with a filtration. A filtered object
$\langle E,F^\hdot \rangle$ is {\em concentrated in a filtered
  degree $i$} if $F^iE=E$ and $F^{i+1}E=0$. For any filtered object
$\langle E,F^\hdot \rangle$ and any integers $i \leq j$, we denote
$$
F^{[i,j]}E = F^iE/F^jE,
$$
and we denote $\gr^i_FE=E^{[i,i]}$, so that $\gr^\hdot_FE$ is the
associated graded quotient of the filtration $F^\hdot$. If $F$ in
$\gr^\hdot_F$ is clear from the context, we will drop it from
notation. For any $i \leq i' \leq j' \leq j$, we have a natural
commutative square
\begin{equation}\label{F.ij}
\begin{CD}
F^{[i,j]}E @<<< F^{[i',j]}E\\
@VVV @VVV\\
F^{[i,j]}E @<<< F^{[i',j']}E,
\end{CD}
\end{equation}
and this square is cartesian and cocartesian. We also have a natural
isomorphism
\begin{equation}\label{compl.iso}
\lim_{\overset{i}{\to}}\lim_{\overset{j}{\gets}}F^{[i,j]}E \cong
\lim_{\overset{j}{\gets}}\lim_{\overset{i}{\to}}F^{[i,j]}E
\end{equation}
and a natural commutative diagram
\begin{equation}\label{conv.sq}
\begin{CD}
E @>{p}>> \displaystyle\lim_{\overset{j}{\gets}}E/F^jE\\
@A{i}AA @AAA\\
\displaystyle\lim_{\overset{i}{\to}}F^iE @>>>
\displaystyle\lim_{\overset{j}{\gets}}
\displaystyle\lim_{\overset{i}{\to}}F^{[i,j]}E.
\end{CD}
\end{equation}

\begin{defn}\label{compl.def}
A filtration $F^\hdot$ on an object $E \in \C$ is {\em convergent},
or equvalently, the object $E$ is {\em convergent with respect to
  $F^\hdot$}, if all the maps in the diagram \eqref{conv.sq} are
isomorphisms. The {\em completion} $\wh{E}$ of a filtered object
$\langle E,F^\hdot \rangle$ in $\C$ is the object \eqref{compl.iso}.
\end{defn}

\begin{exa}\label{sum.exa}
Assume given a graded object $E = \bigoplus_{i \in \Z}E_i \in \C$,
and filter it by
$$
F^iE = \bigoplus_{j \geq i}E_j \subset E.
$$
Then in general, this filtration is {\em not} convergent. The
completion $\wh{E}$ can be described as
$$
\wh{E} = \prod_{j \geq i}E_j \oplus \bigoplus_{j < i}E_j,
$$
where $i$ is an arbitrary fixed integer (the right-hand side does
not depend on the choice of $i$). If all the components $E_i$ are
identified with the same object $E' \in \C$, then it is convenient
to use shorthand notation
$$
E = E'[u,u^{-1}],
$$
that is, Laurent polynomials in one formal variable $u$ of grading
degree $1$ with coefficients in $E'$. Then we have $\wh{E}=E'((u))$,
the formal Laurent power series in the same variable $u$ with the
same coefficients.
\end{exa}

In general, our terminology is misleading, since there are no
natural maps between a filtered object $E$ and its completion
$\wh{E}$. However, say that a filtration $F^\hdot$ on $E$ is {\em
  exhaustive} if the map $i$ in \eqref{conv.sq} is an
isomorphism. Then $p$ gives a natural map
$$
E \to \wh{E} = \lim_{\overset{j}{\gets}}E/F^jE.
$$
Alternatively, say that a filtration $F^\hdot$ is {\em bounded
  below} if $F^iE=0$ for $i \gg 0$. Then $p$ is an isomorphism, and
$i$ gives a natural map
$$
\wh{E} = \lim_{\overset{i}{\to}}E \to E.
$$
For any filtered object $\langle E,F^\hdot \rangle$, the filtration
$F^\hdot$ induces a filtration on the completion $\wh{E}$, and this
filtration is convergent. Two different filtrations on the same
object can have the same completion; here is one standard
situations when this happens.

\begin{defn}
Two filtrations $F^\hdot_1$, $F^\hdot_2$ on an object $E \in \C$ are
{\em commensurable} if for any integer $i$, there exist integers
$i_1$, $i_2$ such that $i_1 \leq i \leq i_2$ and
$$
F^{i_2}_1E \subset F^i_2E \subset F^{i_1}_1E, \qquad F^{i_2}_2E
\subset F^i_1E \subset F^{i_1}_2E.
$$
\end{defn}

Then for any two commensurable filtrations $F_1$, $F_2$ and for any
two integers $i \leq j$, $F^{[i,j]}_1E$ is concentrated in a finite
range of filtered degrees with respect to $F_1$, so that it is
complete, and analogously, $F^{[i,j]}_2E$ is concentrated in a
finite range of filtered degress with respect to $F_2$. The
completions of $E$ with respect to $F_1$ and $F_2$ are canonically
identified.

For example, for any integers $n \geq 1$, $m$, and filtration $F^\hdot$,
define the {\em $n$-th rescaling} $F^\hdot_{[n]}$ and the {\em shift
  by $m$} $F^\hdot_{m}$ of the filtration $F^\hdot$ by
\begin{equation}\label{scale}
F^i_{[n]}E = F^{in}E, \qquad F^i_mE = F^{i+m}E, \qquad i \in \Z.
\end{equation}
Then all rescalings and shifts of a given filtration are mutually
commensurable. For any object $E$ equipped with a filtration
$F^\hdot$, and for any integer $n \geq 1$, we will denote by
$E^{[n]}$ the filtered object $\br{E,F^\hdot_{[n]}}$.

\subsection{Filtered complexes.}\label{DF.subs}

A {\em filtered complex} in $\C$ is a filtered object in the abelian
category $C_\idot(\C)$ of unbounded chain complexes in $\C$. We
denote the category of filtered complexes by $CF_\idot(\C)$. A map
$f:E_\idot \to E'_\idot$ in $CF_\idot(\C)$ is a {\em filtered
  quasiisomorphism} if the corresponding map $\gr^i(f):\gr^iE_\idot
\to \gr^iE'_\idot$ is a quasiisomorphism for any integer $i$. Just
like inverting quasiisomorphisms in $C_\idot(\C)$ gives the derived
category $\D(\C)$, inverting filtered quasiisomorphisms in
$CF_\idot(\C)$ gives the triangulated {\em filtered derived
  category} $\DF(\C)$, see e.g.\ \cite{BBD}. For any filtered
quasiisomorphism $f:\langle E_\idot,F^\hdot \rangle \to \langle
E'_\idot,F^\hdot \rangle$ and any integers $i \leq j$, the induced
map
$$
F^{[i,j]}(f):F^{[i,j]}E_\idot \to F^{[i,j]}E'_\idot
$$
is also a quasiisomorphism, so that for any object $E \in \DF(\C)$,
we have well-defined objects $F^{[i,j]}E \in \D(\C)$ and commutative
squares \eqref{F.ij}. On the other hand, the map $f:E_\idot \to
E'_\idot$ is {\em not} required to be a quasiisomorphism, and
neither are the maps $f:F^iE_\idot \to F^iE'_\idot$. Thus for
an arbitrary filtered complex $\langle E_\idot,F^\hdot \rangle$, the
object $E \in \D(\C)$ corresponding to $E_\idot$ cannot in general
be recovered from the object in $\DF(\C)$ corresponding to $\langle
E_\idot,F^\hdot \rangle$. However, every object $E \in \DF(\C)$ can be
represented by a complex $E_\idot$ in $\C$ equipped with a
convergent filtration $F^\hdot$, and such a representative is unique
up to a quasiisomorphism. Thus one can think of objects in $\DF(\C)$
as corresponding to convergent filtered complexes. In particular,
for any $E \in \DF(\C)$, we have a well-defined completion $\wh{E}
\in \D(\C)$ (explicitly, it is still given by \eqref{compl.iso},
where the limits have to be replaced with their derived
functors). Sending $E$ to $\wh{E}$ gives a completion functor
\begin{equation}\label{df.compl}
\DF(\C) \to \D(\C).
\end{equation}
The filtration $F^\hdot$ on $E$ generates the standard spectral
sequence
\begin{equation}\label{spectr}
H_\idot(\gr^\hdot(E)) \Rightarrow H_\idot(\wh{E}),
\end{equation}
where $H_\idot(-)$ stands for homology objects, and the specral
sequence converges in the usual sense: its $E_\infty$-term is
naturally identified with the associated graded quotient of the
right-hand side with respect to the natural filtration.

The two most standard examples of a filtration on a complex are the
following ones.

\begin{exa}\label{stu.can.exa}
For any complex $E_\idot$ in $\C$, the {\em stupid filtration}
$F^\hdot E_\idot$ and the {\em canonical filtration} $\tau^\hdot
E_\idot$ are given by
$$
F^iE_j = \begin{cases} E_j, &\quad j +i \leq 0,\\0, &\quad j +i > 0,
\end{cases}, \qquad
\tau^iE_j = \begin{cases} 0, &\quad j > i,\\ \Ker d, &\quad
  j=i,\\ E^j, &\quad j < i.
\end{cases}
$$
Both filtrations are convergent.
\end{exa}

\begin{remark}
Our $\tau^i$ is usually denoted $\tau_{\leq -i}$; we change notation
to keep all filtrations decreasing and simplify the statements.
\end{remark}

Canonical filtration is called canonical because it is compatible
with quasiisomorphisms: any quasiisomorphism $f:E_\idot \to
E'_\idot$ in $C_\idot(\C)$ is a filtered quasiisomorphism with
respect to the canonical filtration (the spectral sequence
\eqref{spectr} trivially degenerates for dimension
reasons). Therefore sending $E_\idot$ to $\langle E_\idot,\tau^\hdot
\rangle$ descends to a functor $\D(\C) \to \DF(\C)$. This functor is
fully faithful. For any integer $i$, the associated graded quotient
$\gr^i_\tau E_\idot$ with respect to the canonical filtration is
naturally quasiisomorphic to the $i$-th homology object
$H_i(E_\idot)$ of the complex $E_\idot$. To make this
quasiisomorphism more explicit, one introduces another functorial
decreasing filtration $\beta^\hdot$ on $E_\idot$ given by
\begin{equation}\label{beta}
\beta^iE_j = \begin{cases} 0, &\quad j > i,\\ \Im d, &\quad
  j=i,\\ E_j, &\quad j < i.
\end{cases}
\end{equation}
Then $\tau^{i+1}E_\idot \subset \beta^iE_\idot \subset
\tau^iE_\idot$, the embedding $\tau^{i+1}E_\idot \to \beta^iE_\idot$
is a quasiisomorphism, and the quotient
$\tau^iE_\idot/\beta^iE_\idot$ is isomorphic to the homology object
$H_i(E_\idot)$ placed in homological degree $i$.

Stupid filtration is sufficiently functorial to define a functor
$C_\idot(\C) \to \DF(\C)$, $E_\idot \mapsto \langle E_\idot,F^\hdot
\rangle$, but the functor does not factor through $\D(\C)$:
conversely, it is already a fully faithful embedding on the level of
the category $C_\idot(\C)$, and it identifies $C_\idot(\C)$ with the
heart of a natural $t$-structure on $\DF(\C)$ (see
e.g. \cite[Appendix]{beil}).

For a more interesting example of a filtered complex, denote by
$C_{\idot,\idot}(\C)$ the category of chain bicomplexes in $\C$,
again unbounded in either coordinate. We will call the first
coordinate the {\em vertical direction}, and the second coordinate
the {\em horizontal direction}. Then any $E_{\idot,\idot} \in
C_{\idot,\idot}(\C)$ has at least four natural filtrations --
namely, the stupid and the canonical filtration in either of the two
directions. We also have the two totalization functors from
$C_{\idot,\idot}(\C)$ to $C_\idot(\C)$, namely, the sum-total and
the product-total complex functors
$$
\tot,\Tot:C_{\idot,\idot}(\C) \to C_\idot(\C),
$$
with a natural map $\tot \to \Tot$ between them. The four
filtrations on $E_{\idot,\idot}$ induce filtrations on
$\tot(E_{\idot,\idot})$ and $\Tot(E_{\idot,\idot})$. Sometimes, some
of these filtrations are automatically convergent. For example,
assume that $E_{\idot,\idot}$ is a second-quadrant bicomplex -- that
is, $E_{i,j}=0$ unless $i \leq 0 \leq j$. Then
$\Tot(E_{\idot,\idot})$ is convergent with respect to the stupid
filtration in the vertical direction and the canonical filtration in
the horizontal direction, while $\tot(E_{\idot,\idot})$ is
convergent with respect to the stupid filtration in the horizontal
direction and the canonical filtration in the vertical
direction. Thus the spectral sequences \eqref{spectr} for
$\tot(E_{\idot,\idot})$ and $\Tot(E_{\idot,\idot})$ are different,
and the natural map $\tot(E_{\idot,\idot}) \to
\Tot(E_{\idot,\idot})$ need not be a quasiisomorphism. Here is the
standard example of such a situation.

\begin{lemma}\label{snake.lim.le}
Assume given objects $E_n \in \C$ and morphisms $e_n:E_n \to
E_{n+1}$, $n \geq 0$, and consider the bicomplex $E_{\idot,\idot}$
given by $E_{-n,n}=E_{-(n+1),n}=E_n$, $E_{i,j}=0$ otherwise, with the
non-trivial horizontal differentials given by the identity maps, and
the vertical differentials given by the maps $e_i$. Assume that
either the category $\C$ satisfies AB5 (filtered direct limits are
exact), or all the maps $e_i$ are injective. Then
$\Tot(E_{\idot,\idot})$ is quasiisomorphic to $E_0$ placed in degree
$0$, while $\tot(E_{\idot,\idot})$ is quasiisomorphic to the cone of
the natural map
$$
E_0 \to \lim_{\overset{i}{\to}}E_i,
$$
where the limit is take with respect to the maps $e_i$. The natural
map from $\tot(E_{\idot,\idot})$ to $\Tot(E_{\idot,\idot})$ fits
into a distinguished triangle
$$
\begin{CD}
\tot(E_{\idot,\idot}) @>>> \Tot(E_{\idot,\idot}) @>>>
\displaystyle\lim_{\overset{i}{\to}}E_i @>>>.
\end{CD}
$$
\end{lemma}

\begin{remark}
Explicitly, the bicomplex $E_{\idot,\idot}$ is given by the
snake diagram
$$
\begin{CD}
... @>{\id}>> E_3\\
@.@AA{e_2}A\\
@.E_2 @>{\id}>> E_2\\
@.@.@AA{e_1}A\\
@.@.E_1 @>{\id}>> E_1\\
@.@.@.@AA{e_0}A\\
@.@.@.E_0.
\end{CD}
$$
\end{remark}

\proof{} The first claim is clear from \eqref{spectr}. For second,
compute the limit by the telescope construction. The last claim is
clear.
\endproof

\subsection{Truncations and functoriality.}\label{t.subs}

For a general bicomplex $E_{\idot,\idot}$, neither
$\tot(E_{\idot,\idot})$ nor $\Tot(E_{\idot,\idot})$ is convergent
with respect to either of the filtrations -- the situation for
$\tot(E_{\idot,\idot})$ is described in Example~\ref{sum.exa}, and
$\Tot(E_{\idot,\idot})$ is dual. However, the completion of
$\tot(E_{\idot,\idot})$ with respect to the stupid filtration in the
vertical direction is isomorphic to its completion with respect to
the canonical filtration in the horizontal direction. We will need a
slight generalization of this fact.

First of all, fix the identification $C_{\idot,\idot}(\C) \cong
C_\idot(C_\idot(\C))$ so that the outer index in
$C_\idot(C_\idot(\C))$ corresponds to the horizontal direction in
$C_{\idot,\idot}(\C)$. Since $C_\idot(\C)$ is an abelian category in
its own right, this defines the truncations functors $\tau^\hdot$
and their counterparts $\beta^\hdot$ of \eqref{beta} on
$C_{\idot,\idot}(\C) \cong C_\idot(C_\idot(\C))$. Now denote by
$$
\tot^f:C_{\idot,\idot}(\C) \to CF_\idot(\C)
$$
the functor sending a bicomplex $E_{\idot,\idot}$ to
$\tot(E_{\idot,\idot})$ with the filtration induced by the stupid
filtration in the vertical direction on $E_\idot$. Say that a map
$f:E_{\idot,\idot} \to E'_{\idot,\idot}$ is a {\em horizontal
  quasiisomorphism} if $f:E_{i,\idot} \to E'_{i,\idot}$ is a
quasiisomorphism for any integer $i$. Then inverting horizontal
quasiisomorphisms in $C_{\idot,\idot}(\C) \cong
C_\idot(C_\idot(\C))$ gives the derived category $\D(C_\idot(\C))$,
and $\tot^f$ descends to a functor
$$
\D(C_\idot(\C)) \to \DF(\C).
$$
It is well-known (see e.g. \cite[Appendix]{beil}) that this functor
is an equivalence of categories. Thus in particular, the truncation
functors $\tau^\hdot$ of Example~\ref{stu.can.exa} are well-defined
on the filtered derived category $\DF(\C)$.

It turns out that in fact both $\tau^\hdot$ and the related
truncation functors $\beta^\hdot$ of \eqref{beta} can be defined
already on the category $CF_\idot(\C)$. Namely, for any filtered
complex $\langle E_\idot,F^\hdot \rangle$ with differential $d$, and
any integers $n$, $i$, let
\begin{equation}\label{tau.beta.F}
\begin{aligned}
\tau^nE_i &= d^{-1}(F^{n+1-i}E_{i-1}) \cap F^{n-i}E_i \subset E_i,\\
\beta^nE_i &= F^{n+1-i}E_i + d(F^{n-i}E_{i+1}) \subset E_i.
\end{aligned}
\end{equation}
Then for any $n$, $\tau^nE_\idot$ and $\beta^nE_\idot$ are
subcomplexes in $E_\idot$, and we have $\tau^{n+1}E_\idot \subset
\beta^nE_\idot \subset \tau^nE_\idot$. Equip $\tau^nE_\idot$ and
$\beta^nE_\idot$ with the filtrations induced by $F^\hdot$, and
denote
\begin{equation}\label{h.i}
H_n(E_\idot) = \tau^nE_\idot/\beta^nE_\idot \subset CF_\idot(\C),
\end{equation}
again equipped with the induced filtration. For any bicomplex
$E_{\idot,\idot}$ and any integer $n$, we then obviously have
$$
\tau^n\tot^f(E_{\idot,\idot}) \cong \tot^f(\tau^nE_{\idot,\idot}),
\qquad \beta^n\tot^f(E_{\idot,\idot}) \cong
\tot^f(\beta^nE_{\idot,\idot}),
$$
and $H_n(\tot^f(E_{\idot,\idot}))$ coincides with
$\tot^f(H_n(E_{\idot,\idot}))$, where $H_n(E_{\idot,\idot}) \in
C_\idot(\C)$ is the $n$-th homology object of $E_{\idot,\idot} \in
C_{\idot,\idot}(\C) \cong C_\idot(C_\idot(\C))$.

\begin{lemma}\label{tau.le}
For any filtered complex $\langle E_\idot,F^\hdot \rangle$ in $\C$
and any integer $n$, $i$, we have
$$
\gr^i_F\tau^nE_\idot \cong \tau^{n-i}\gr^i_FE_\idot, \qquad
\gr^i_F\beta^nE_\idot \cong \beta^{n-i}\gr^i_FE_\idot.
$$
Moreover, the embedding $\tau^{n+1}E_\idot \to \beta^nE_\idot$ is a
filtered quasiisomorphism, and we have $\gr^i_FH_n(E_\idot) \cong
H_{n-i}(\gr^i_FE_\idot)[n]$, so that the filtration $F^\hdot$ on
$H_n(E_\idot)$ is the stupid filtration shifted by $n$, as in
\eqref{scale}. Finally, for any integer $i$, the filtrations
$\tau^\hdot$, $\beta^\hdot$ and $F^\hdot$ on $E_i$ are
commensurable.
\end{lemma}

\proof{} Clear. \endproof

We will also need one fact about functoriality properties of the
filtrations $\tau^\hdot$, $\beta^\hdot$. Assume given abelian
categories $\C$, $\C'$. Any additive functor $\phi:\C \to \C'$
naturally extends to a functor $\phi_\idot:\C \to \C'$ by applying
it termwise, that is, setting $\phi(E_\idot)_i = \phi(E_i)$ for any
complex $E_\idot$ in $\C$ and any integer $i$. For any additive
functor $\phi_\idot:\C \to C_\idot(\C')$, we define its extension
$\phi_\idot:C_\idot(\C) \to C_\idot(\C')$ as the composition
$$
\begin{CD}
C_\idot(\C) @>{\phi_{\idot,\idot}}>> C_\idot(C_\idot(\C')) \cong
C_{\idot,\idot}(\C') @>{\tot}>> C_\idot(\C').
\end{CD}
$$
That is, we apply $\phi_\idot$ termwise, and then take the sum-total
complex $\tot(-)$. Explicitly, for any complex $E_\idot$ in
$\C$, the complex $\phi_\idot(E_\idot)$ in $\C'$ has terms
\begin{equation}\label{phi.c}
\phi_\idot(E_\idot)_i = \bigoplus_j \phi_{i-j}(E_j).
\end{equation}
We also define the functor
\begin{equation}\label{phi.f}
\phi^f = \tot^f \circ \phi_{\idot,\idot}:C_\idot(\C) \to CF_\idot(\C').
\end{equation}
Moreover, say that a filtration $F^\hdot$ on a complex $E_\idot$ is
{\em termwise-split} if for any $i$, $j$, the embedding $F^iE_j \to
E_j$ is a split injection, and denote by $\bCF_\idot(\C) \subset
CF_\idot(\C)$ the full subcategory spanned by complexes with
termwise-split filtrations. Then since any functor sends split
injections to split injections, we can further extend the functor
$\phi^f$ to a functor
\begin{equation}\label{phi.f.b}
\phi^f:\bCF_\idot(\C) \to \bCF_\idot(\C')
\end{equation}
by setting
$$
F^n\phi^f(E_\idot)_i = \bigoplus_j\phi_{i-j}(F^{n+j}E_j) \subset
\phi_\idot(E_\idot)_i,
$$
where we recall that $\phi_\idot(E_\idot)_i$ is explicitly given by
\eqref{phi.c}. Note that if a complex $E_\idot$ is placed in
filtered degree $0$, then the filtration is tautologically termwise
split, and \eqref{phi.f.b} coincides with \eqref{phi.f}. If moreover
$E_\idot$ is concentrated in homological degree $0$, then
$\phi^f(E_\idot)$ is just $\phi_\idot(E^0)$ equipped with the stupid
filtration.

\begin{lemma}\label{tau.func.le}
Assume given an additive functor $\phi_\idot:\C \to
C_\idot(\C')$. Then for any integer $n$ and any complex $E_\idot$ in
$\C$ equipped with a termwise-split filtration $F^\hdot$, the
filtrations on the truncations $\tau^nE_\idot,\beta^nE_\idot \subset
E_\idot$ of \eqref{tau.beta.F} are termwise-split, and the natural
maps $\phi^f(\tau^nE_\idot) \to \phi^f(E_\idot)$,
$\phi^f(\beta^nE_\idot) \to \phi^f(E_\idot)$ factor through natural
maps
$$
\phi^f(\tau^nE_\idot) \to \tau^n(\phi^f(E_\idot)) \subset
\phi^f(E_\idot), \
\phi^f(\beta^nE_\idot) \to \beta^n(\phi^f(E_\idot)) \subset
\phi^f(E_\idot).
$$
Moreover, if the functor $\phi_\idot$ is right-exact, the map
$\phi^f(\tau^nE_\idot) \to \tau^n\phi^f(E_\idot)$ is surjective.
\end{lemma}

\proof{} By Lemma~\ref{tau.le}, for any integers $i$, $n$, we have
$(\tau^nE_\idot)_i = F^{n-i}(\tau^nE_\idot)_i$,
$F^{n+1-i}(\tau_nE_\idot)_i = F^{n+1-i}E_i$, and a splitting of
the projection
$$
F^{n-i}(\tau^nE_\idot)_i \to \gr^{n-i}_F(\tau^nE_\idot)_i \cong
(\tau^i\gr^{n-i}_FE_\idot)_i
$$
is induced by a splitting of the projection $F^{n-i}E_i \to
\gr^{n-i}_FE_i$. Therefore the filtration on $\tau^nE_\idot$ induced
by $F^\hdot$ is termwise-split. The same argument shows that the
filtration on $\beta^nE_\idot$ is also termwise-split. Then applying
the additive functor $\phi$ commutes with taking associated graded
quotients with respect to $F^\hdot$, and by Lemma~\ref{tau.le}, it
suffices to prove the remaining claims when $E_\idot$ is
concentrated in a single filtered degree. Then up to a shift,
$\tau^\hdot$ is the canonical filtration, $\beta^\hdot$ is given by
\eqref{beta}, and the claims are obvious.
\endproof

\section{Mixed complexes.}\label{mix.sec}

\subsection{Mixed complexes and expansions.}\label{mix.gen.subs}

As in Section~\ref{filt.sec}, we fix an abelian category $\C$
satisfying  $AB4$ and $AB4^*$.

\begin{defn}\label{mx.def}
A {\em mixed complex} $\langle V_\idot,B \rangle$ in $\C$ is a
complex $V_\idot$ equipped with a map $B:V_\idot \to V_\idot[-1]$
such that $B^2=0$. A {\em map of mixed complexes} $f:\langle
V_\idot,B \rangle \to \langle V'_\idot,B \rangle$ is a map of
complexes $f:V_\idot \to V'_\idot$ that commutes with $B$.
\end{defn}

For every map of mixed complexes, its cone is also a mixed
complex. For every integer $n$, the canonical truncation
$\tau^nV_\idot$ of a mixed complex $\br{V_\idot,B}$ is a mixed
complex, and the natural embedding $\tau^nV_\idot \to V_\idot$ is a
map of mixed complexes.

\begin{defn}\label{exp}
The {\em expansion} $\Exp(\langle V_\idot,B \rangle)$ of a mixed
complex $\langle V_\idot,B \rangle)$ is the complex isomorphic to
$V_\idot[u^{-1}]$ as a graded object in $\C$, with $u$ being a
formal variable of cohomological degree $2$, and equipped with the
differential $d+Bu$, where $d:V_\idot \to V_{\idot-1}$ is the
differential in the complex $V_\idot$.
\end{defn}

We will write $\Exp(V_\idot)$ when $B$ is clear from the
definition. The meaning of the shorthand notation $V_\idot[u^{-1}]$
is the same as in Example~\ref{sum.exa}; as usual, polynomials in
$u^{-1}$ are understood as a module over polynomials in $u$ via the
identification $k[u^{-1}] = k[u,u^{-1}]/uk[u]$. Equivalently,
$\Exp(V_\idot)$ is the sum-total complex of the bicomplex
\begin{equation}\label{exp.bi}
\begin{CD}
@>{B}>> V_\idot[3] @>{B}>> V_\idot[2] @>{B}>> V_\idot[1] @>{B}>> V_\idot,
\end{CD}
\end{equation}
and $u:\Exp(V_\idot) \to \Exp(V_\idot)[2]$ is induced by the obvious
map shifting the bicomplex by $1$ in either direction. By
definition, we have a short exact sequence of complexes
$$
\begin{CD}
0 @>>> V_\idot @>>> \Exp(V_\idot) @>{u}>> \Exp(V_\idot)[2] @>>> 0.
\end{CD}
$$
The same shorthand notation is used in the following definition.

\begin{defn}\label{per}
The {\em periodic\/} resp. {\em co-periodic\/} resp. {\em polynomial
  periodic\/} expansions $\Per(V_\idot)$ resp. $\cPer(V_\idot)$
resp. $\per(V_\idot)$ of a mixed complex $\langle V_\idot,B \rangle$
are the complexes in $\C$ given by
$$
\Per(V_\idot) = V_\idot((u)), \qquad \cPer(V_\idot)((u^{-1})),
\qquad \per(V_\idot) = V_\idot[u,u^{-1}]
$$
as graded objects, each equipped with the differential $d+Bu$.
\end{defn}

By definition, for every mixed complex $\br{V_\idot, B}$, we have
natural maps
\begin{equation}\label{per.Per}
\begin{CD}
\Per(V_\idot) @<{l}<< \per(V_\idot) @>{r}>> \cPer(V_\idot).
\end{CD}
\end{equation}
In general, neither of these maps is an isomorphism, or even a
quasiisomorphism. We also have
\begin{equation}\label{per.exp}
\Per(V_\idot) = \lim_{\overset{u}{\gets}}\Exp(V_\idot).
\end{equation}
All the periodic expansions of Definition~\ref{per} are obviously
exact functors on the category of mixed complexes. In particular,
they all commute with finite direct sums. The polynomial periodic
expansion functor $\per(-)$ also commutes with infinite direct
sums. A filtration $F^\hdot$ on the complex $V_\idot$ preserved by
the differential $B$ induces filtrations on $\Per(V_\idot)$,
$\cPer(V_\idot)$, $\per(V_\idot)$ by setting
\begin{equation}\label{mixed.filt}
\begin{aligned}
F^i\Per(V_\idot) &= \Per(F^iV_\idot),\\
F^i\cPer(V_\idot) &= \cPer(F^iV_\idot),\\
F^i\per(V_\idot) &= \per(F^iV_\idot)
\end{aligned}
\end{equation}
for any integer $i$. The corresponding associated graded quotients
are given by
\begin{equation}\label{mixed.filt.gr}
\begin{aligned}
\gr^F_i\Per(V_\idot) &= \Per(\gr^F_iV_\idot),\\
\gr^F_i\cPer(V_\idot) &= \cPer(\gr^F_iV_\idot),\\
\gr^F_i\per(V_\idot) &= \per(\gr^F_iV_\idot),
\end{aligned}
\end{equation}
again for any integer $i$.

\begin{defn}\label{mix.qis}
A mixed complex $\br{V_\idot,B}$ with the differential
$d:V_{\idot+1} \to V_\idot$ is {\em contractible} if there exists a
map $h:V_\idot \to V_{\idot+1}$ such that $dh+hd=\id$ and $Bh=hB$. A
mixed complex is {\em acyclic} if it is acyclic with respect to the
differential $d$, and it is {\em strongly acyclic} if it is finite
extension of contractible mixed complexes.  A map of mixed complexes
$f:V_\idot \to V'_\idot$ is a {\em quasiisomorphism} resp.\ {\em
  strong quasiisomorphism} if if its cone is acyclic resp.\ strongly
acyclic.
\end{defn}

\begin{defn}\label{bnd.l}
A mixed complex $\br{V_\idot,B}$ is {\em bounded from below} if for
some integer $n$, $\tau^nV_\idot$ is acylic, and it is {\em strongly
  bounded from above} if for some integer $n$, the map
$\tau^nV_\idot \to V_\idot$ is a strong quasiisomorphism.
\end{defn}

\begin{lemma}\label{mix.qis.le}
If a map $f:V_\idot \to V'_\idot$ of mixed complexes is a
quasiisomorphism, then the map $\Per(f):\Per(V_\idot) \to
\Per(V'_\idot)$ is a quasiisomorphism. If $f$ is a strong
quasiisomorphism, then the maps $\per(f)$ and $\cPer(f)$ are also
quasiisomorphisms.
\end{lemma}

\proof{} For any mixed complex $\langle V_\idot,B \rangle$, let
$\per_{\idot,\idot}(V_\idot)$ be the bicomplex
\begin{equation}\label{per.bi}
\begin{CD}
@>{B}>> V_\idot[1] @>{B}>> V_\idot @>{B}>> V_\idot[-1] @>{B}>>,
\end{CD}
\end{equation}
a periodic version of the bicomplex \eqref{exp.bi}. Then by
definition, we have $\per(V_\idot) =
\tot(\per_{\idot,\idot}(V_\idot))$, and $\Per(V_\idot)$ is its
completion with respect to the stupid filtration in the horizontal
direction. Applying the spectral sequence \eqref{spectr}, we get the
first claim. For the second claim, it suffices to check that for a
strongly acyclic mixed complex $V_\idot$, $\per(V_\idot)$ and
$\cPer(V_\idot)$ are acyclic. This is obvious: $\per(-)$ and
$\cPer(-)$ are exact functors, and they send contractible mixed
complexes to contractible complexes.
\endproof

\begin{corr}\label{b.be}
Assume given a mixed complex $\langle V_\idot,B \rangle$.
\begin{enumerate}
\item If $V_\idot$ is bounded from below in the sense of
  Definition~\ref{bnd.l}, then the natural map $l:\per(V_\idot) \to
  \Per(V_\idot)$ is a quasiisomorphism.
\item If $V_\idot$ is strongly bounded from above in the sense of
  Definition~\ref{bnd.l}, then the natural map $r:\per(V_\idot) \to
  \cPer(V_\idot)$ is a quasiisomorphism.
\end{enumerate}
\end{corr}

\proof{} By Lemma~\ref{mix.qis.le}, it suffices to prove that if we
have an integer $n$ such that $V_i=0$ for $i \geq n$ resp.\ $i \leq
n$, then $l$ resp.\ $r$ is quasiisomorphism. In fact, under these
assumptions it is even an isomorphism.
\endproof

\subsection{Cyclic groups.}\label{cycl.grp.subs}

To obtain a useful example of a mixed complex, fix an integer $n
\geq 1$, and let $C = \Z/n\Z$ be the cyclic group of order $n$, with
the generator $\sigma \in C$. If we denote $\K_1=\K_0=\Z[C]$, then
we have a natural exact sequence
\begin{equation}\label{K.eq}
\begin{CD}
0 @>>> \Z @>{\kappa_1}>> \K_1 @>{\id - \sigma}>> \K_0 @>{\kappa_0}>>
\Z @>>> 0
\end{CD}
\end{equation}
of $\Z[C]$-modules, where the action of $C$ on $\Z$ is
trivial. We thus have a complex $\K_\idot$ of $\Z[C]$-modules of
length $2$, with homology in degrees $0$ and $1$ identified with
$\Z$. Denote
$$
B = \kappa_1 \circ \kappa_0:\K_0 \to \K_1.
$$
Then $B \circ (\id -\sigma) = (\id - \sigma) \circ B = 0$, so that
$B$ defines a map of complexes
\begin{equation}\label{B.K.grp}
B:\K_\idot \to \K_\idot[-1]
\end{equation}
turning $\K_\idot$ into a mixed complex. Explicitly, we have
$$
B = \id + \sigma + \dots + \sigma^{n-1}.
$$
Moreover, fix a ring $R$. Then the category $\C$ of left
$R[C]$-modules is an abelian category satisfying $AB4$ and $AB4^*$,
and for any $R[C]$-module $E$, we have a natural mixed complex
\begin{equation}\label{K.grp}
\K_\idot(E) = \K_\idot \otimes E
\end{equation}
in $\C$, with the differential $B$ induced by the map
\eqref{B.K.grp}. The expansion $\Exp(\K_\idot(E))$ is then the
standard periodic free resolution of $E$. Since the complex
$\K_\idot(E)$ is concentrated in a finite number of degrees, the
maps \eqref{per.Per} for $\K_\idot(E)$ are quasiisomorphisms -- in
fact, isomorphisms -- and the complex $\per(\K_\idot(E)) \cong
\Per(\K_\idot(E)) \cong \cPer(\K_\idot(E))$ is acyclyc.

Taking coinvariants with respect to $C$, -- or equivalently, with
respect to the generator $\sigma \in C$ -- we obtain a mixed
complex $\K_\idot(E)_\sigma$ in the category of $R$-modules. Its
expansion
\begin{equation}\label{c.eq}
C_\idot(C,E) = \Exp(\K_\idot(E)_\sigma) \cong \Exp(\K_\idot(E))_\sigma
\end{equation}
is the standard periodic complex computing the homology groups
$H_\idot(C,E)$, and the periodic expansion
\begin{equation}\label{tate.eq}
\vC_\idot(C,E) = \Per(\K_\idot(E)_\sigma) \cong
\Per(\K_\idot(E))_\sigma
\end{equation}
computes the Tate homology groups $\vH_\idot(C,E)$. All the
non-trivial terms of these complexes are isomorphic to $E$ as an
$R$-modules, and the differentials are $\id - \sigma$ resp.\ $\id +
\sigma + \dots + \sigma^{n-1}$ in odd resp.\ even degrees.

More generally, let $E_\idot$ be a complex of $R[C]$-modules. Then
$\K_\idot(E_\idot) = \K_\idot \otimes E_\idot$ equipped with a map
\eqref{B.K.grp} is still a mixed complex of $R[C]$-modules, and
taking coinvariants, we obtain a mixed complex
$\K_\idot(E_\idot)_\sigma$ of $R$-modules. Since
$\K_\idot(E_\idot)_\sigma$ is no longer concentrated in a finite
range of degrees, the maps \eqref{per.Per} need not be
quasiisomorphisms anymore. However, we can still define the Tate
homology complex $\vC_\idot(C,E_\idot)$ by \eqref{tate.eq}, and for
any quasiisomorphism $E_\idot \to E'_\idot$ of complexes of
$R[C]$-modules, the corresponding map
$$
\vC_\idot(C,E_\idot) \to \vC_\idot(C,E'_\idot)
$$
is a quasiisomorphism by Lemma~\ref{mix.qis.le}. Taking other
periodic expansions of Definition~\ref{per}, we obtain two more
versions of the Tate cohomology complex that we denote by
\begin{equation}\label{tate.bis}
\vbC_\idot(C,E_\idot) = \cPer(\K_\idot(E_\idot)_\sigma), \qquad
\wC_\idot(C,E_\idot) = \per(\K_\idot(E_\idot)_\sigma).
\end{equation}
These are not in general invariant under quasiisomorphisms. To get
an invariance result, we proceed as in Definition~\ref{mix.qis}.

\begin{defn}\label{strong.qis}
A complex $E_\idot$ in an abelian category $\C$ is {\em strongly
  acyclic} if it is a finite extension of contractible complexes. A
map $f:E_\idot \to E'_\idot$ is a {\em strong quasiisomorphism} if
its cone is strongly acylic.
\end{defn}

\begin{defn}\label{strong.bnd}
A complex $E_\idot$ in an abelian category $\C$ is {\em bounded from
  below} resp.\ {\em strongly bounded from below} if for some
integer $n$, $\tau^nE_\idot$ is acylic resp.\ strongly acyclic, and
it is {\em bounded from above} resp.\ {\em strongly bounded from
  above} if for some integer $n$, the map $\tau^nE_\idot \to
E_\idot$ is a quasiisomorphism resp.\ strong quasiisomorphism.
\end{defn}

\begin{lemma}\label{tate.bnd.le}
\begin{enumerate}
\item For every strong quasiisomorphism $f$ of complexes of
  $R[C]$-modules, the maps $\per(f)$ and $\cPer(f)$ are
  quasiisomorphisms.
\item If a complex $E_\idot$ of $R[C]$-modules is bounded from below
  in the sense of Definition~\ref{strong.bnd}, then the map
  $l:\wt{C}_\idot(C,E_\idot) \to \vC_\idot(C,E_\idot)$ is a
  quasiisomorphism. If $E_\idot$ is strongly bounded from above,
  then the map $r:\wt{C}_\idot(C,E_\idot) \to \vbC_\idot(C,E_\idot)$
  is a quasiisomorphism.
\end{enumerate}
\end{lemma}

\proof{} Since $\K(E)_\sigma$ is an exact functor, the claims
immediately follow from Lemma~\ref{mix.qis.le}
resp.\ Corollary~\ref{b.be}.
\endproof

\begin{lemma}\label{tate.compl.le}
For any complex $E_\idot$ of $R[C]$-modules, equip the Tate complex
$\wC_\idot(C,E_\idot)$ with the filtration induced by the stupid
resp. canonical filtration on $E_\idot$. Then its completion
coincides with the complex $\vbC_\idot(C,E_\idot)$
resp. $\vC_\idot(C,E_\idot)$.
\end{lemma}

\proof{} Clear. \endproof

We will also need a version of homology complexes twisted by a
sign. Namely, let $\sigma^\dg = (-1)^{|C|+1}\sigma\in R[C]$, where
$|C|$ is the order of the cyclic group $C$. Then $\sigma^\dg$ is an
invertible element of order $|C|$, so that any $R[C]$-module $E$
gives another $R[C]$-module $E^\dg$ with the generator of $C$ acting
by $\sigma^\dg$. Alternatively, $E^\dg = E \otimes I$, where $I$ is
the one-dimensional sign representation: $I = \Z$ as a $\Z$-module,
and $\sigma$ acts by $(-1)^{|C|+1}$. Then taking coinvariants with
respect to $\sigma^\dg$, we obtain a mixed complex
$\K(E)_{\sigma^\dg}$ and its expansions
\begin{equation}\label{c.tw}
\begin{aligned}
C^\dg_\idot(C,E) &= \Exp(\K(E)_{\sigma^\dg}) \cong
\Exp(\K(E))_{\sigma^\dg}, \\
\vC^\dg_\idot(C,E) &= \Per(\K(E)_{\sigma^\dg}) \cong
\Per(\K(E))_{\sigma^\dg},
\end{aligned}
\end{equation}
the twisted versions of \eqref{c.eq} and
\eqref{tate.eq}. Explicitly, all the terms in $\vC^\dg_\idot(C,E)$
are isomorphic to $E$ as $R$-modules, with the alternating
differentials $\id - \sigma^\dg$, $\id + \sigma^\dg + \dots _
(\sigma^\dg)^{n-1}$. More generally, for any complex $E_\idot$ of
$R[C]$-modules, we have the mixed complex $\K(E_\idot)_{\sigma^\dg}$
and its expansions
\begin{equation}\label{tate.tw}
\begin{aligned}
C_\idot^\dg(C,E_\idot) &= \Exp(\K_\idot(E_\idot)_{\sigma^\dg}), \qquad
\vC_\idot^\dg(C,E_\idot) = \cPer(\K_\idot(E_\idot)_{\sigma^\dg}),\\
\vbC_\idot^\dg(C,E_\idot) &=
\cPer(\K_\idot(E_\idot)_{\sigma^\dg}), \qquad \wC_\idot^\dg(C,E_\idot)
= \per(\K_\idot(E_\idot)_{\sigma^\dg}),
\end{aligned}
\end{equation}
the twisted versions of \eqref{c.eq}, \eqref{tate.eq} and
\eqref{tate.bis}. We of course have
$$
C_\idot^\dg(C,E_\idot) \cong C_\idot(C,E_\idot \otimes I),
$$
and similarly for the periodic complexes, so that the general
properties of periodic complexes also hold for their twisted
versions. To wit, $\vC^\dg_\idot(C,-)$ sends quasiisomorphisms to
quasiisomorphisms, and we have the maps \eqref{per.Per},
Lemma~\ref{tate.bnd.le} and Lemma~\ref{tate.compl.le}.

Finally, let us note that one can modify the constructions of this
Subsection in the following way. Take a larger finite cyclic group
$C'$ containing $C$, and denote by $\K_\idot'$ the complex
\eqref{K.eq} for the group $C'$ considered as a complex of
$R[C]$-modules by restriction with respect to the embedding $C
\subset C'$. Then $\K_\idot'$ is also a mixed complex, and for any
$R[C]$-module $E$, we can form a mixed complex $\K'_\idot(E)_\sigma
= (\K' \otimes E)_\sigma$ and its expansions
\begin{equation}\label{c.cpr}
C_\idot(C,C',E) = \Exp(\K'_\idot(E)_\sigma), \qquad
\vC_\idot(C,C',E) = \Per(\K'_\idot(E)_\sigma).
\end{equation}
However, both $\K_\idot$ and $\K'_\idot$ are complexes of projective
$\Z[C]$-modules, and it is easy to see that there a decomposition
$$
\K'_\idot \cong \K_\idot \oplus K_\idot
$$
of mixed complexes of $\Z[C]$-modules with $K_\idot$ contractible in
the sense of Definition~\ref{mix.qis}. Therefore the complexes
\eqref{c.cpr} are canonically chain-homotopy equivalent to their
counterparts \eqref{c.eq}, \eqref{tate.eq}. Analogously, for any
complex $E_\idot$ of $R[C]$-modules, one can define complexes
$\vC_\idot(C,C',E_\idot)$, $\vbC_\idot(C,C',E_\idot)$,
$\wC_\idot(C,C',E_\idot)$ and their twisted versions as in
\eqref{tate.tw}; however, all these complexes ara canonically
chain-homotopy equivalent to their counterparts without the group
$C'$.

\subsection{Mixed resolutions.}\label{reso.subs}

Now return to the general situation of
Subsection~\ref{mix.gen.subs}: $\C$ is again an arbitrary abelian
category satisfying $AB4$ and $AB4^*$. We note that while every map
of mixed complexes induces a map of their expansions, the converse
is not true -- there are useful $u$-equivariant maps of expansions
that do not come from maps of mixed complexes. We will need one
example of this kind.

\begin{defn}\label{mix.reso}
A {\em mixed resolution} of an object $M \in \C$ is a mixed complex
$\langle M_\idot,B \rangle$ in $\C$, with homology objects
$H_\idot(M_\idot)$, equipped with a map $a:M_0 \to M$ such that
\begin{enumerate}
\item $M_i = 0$ for $i < 0$, and $H_i(M_\idot)=0$ for $i \geq 2$,
\item $a:M_0 \to M$ induces an isomorphism $H_0(M_\idot) \cong M$,
  and
\item $B:M_\idot \to M_\idot[-1]$ induces an isomorphism
  $H_1(M_\idot) \cong H_0(M_\idot)$.
\end{enumerate}
The {\em cohomology class}
$$
\alpha \in \Ext^2(M,M) \cong \Ext^2(H_0(M_\idot),H_1(M_\idot))
$$
of a mixed resolution $M_\idot$ is the class it represents by Yoneda
as a complex.
\end{defn}

\begin{lemma}\label{mix.map.le}
Assume given two mixed resolutions $P_\idot$, $M_\idot$ of the same
object $M \in \C$ and with the same cohomology class $\alpha \in
\Ext^2(M,M)$, and assume that $P_i$ is a projective object in $\C$
for any $i \geq 0$. There there exists a $u$-equivariant map
$$
\nu:\Exp(P_\idot) \to \Exp(M_\idot)
$$
that induces the given isomorphism $H_0(P_\idot) \cong H_0(M_\idot)
\cong M$ in degree $0$.
\end{lemma}

\proof{} Since $P_0$ is projective, the isomorphism $H_0(P_\idot)
\to H_0(M_\idot)$ can be lifted to a map $\nu_0:P_0 \to M_0$, and
since both mixed resolutions have the same cohomology class, and
$P_i$ is projective for any $i \geq 0$, the map $\nu_0$ can be
extended to a map of complexes $\nu_0:P_\idot \to M_\idot$ that
induces the given isomorphism $H_1(P_\idot) \cong H_0(P_\idot) \cong
H_0(M_\idot) \cong H_1(M_\idot)$ on homology in degree $1$. For any
mixed complex $V_\idot$ and any $i \geq 0$, denote by
$F_i\Exp(V_\idot) \subset \Exp(V_\idot)$ be $(-i)$-th term of the
stupid filtration of the bicomplex \eqref{exp.bi} in the horizontal
direction. Then by definition, $u$ sends $F_i\Exp(V_\idot)$ into
$F_{i-1}\Exp(V_\idot)$, and we have $\gr^i_F\Exp(V_\idot) \cong
V_\idot[i]$. By induction, we may assume that we are given a
$u$-equivariant map of complexes
$$
\nu_i:F_i\Exp(P_\idot) \to F_i\Exp(M_\idot),
$$
and we need to extend it to a $u$-equivariant map $\nu_{i+1}$. First
extend it to a $u$-equivarient graded map $\wt{\nu}_{i+1}$ by
setting
$$
\wt{\nu}_{i+1} = u^{-1}\nu_iu
$$
on $u^{-(i+1)}P_\idot \subset F_i\Exp(P_\idot)$. Then while
$\wt{\nu}_{i+1}$ does not necessarly commute with the differential
$d+uB$, we know by induction that the commutator
$$
(d+uB)\wt{\nu}_{i+1}+\wt{\nu}_{i+1}(d+uB):F_{i+1}\Exp(P_\idot) \to
F_{i+1}\Exp(M_\idot)
$$
is divisible by $u^i$. That is, we have
$(d+uB)\wt{\nu}_{i+1}+\wt{\nu}_{i+1}(d+uB) = u^i\tau_i$ for some map
$\tau_i:P_\idot \to M_\idot$. This is a map of complexes, and
moreover, it is equal to $0$ on homology: if $i=1$, this follows
from our construction of the map $\nu_0$, and for $i \geq 2$, this
follows directly from Definition~\ref{mix.reso}~\thetag{i}. Since
$P_\idot$ is a complex of projective objects, $\tau_i$ must be
homotopic to $0$ -- that is, $\tau_i = dh_i+h_id$ for some graded
map $h_i:P_\idot \to M_{\idot+1}$. To finish the proof, take $\nu_i
= \wt{\nu}_i + u^ih_i$.
\endproof

\section{Cyclic homology.}\label{cyc.sec}

\subsection{Cyclic complexes.}\label{cp.subs}

We will use the same notation  and conventions as in \cite[Section
  1]{ka1}; in particular, for any small category $I$ and ring $R$,
we denote by $\Fun(I,R)$ the abelian category of functors from $I$
to the category of $R$-modules, and we denote by $\D(I,R)$ its
derived category. For any functor $\gamma:I' \to I$ between small
categories, we denote by $\gamma^*:\Fun(I',R) \to \Fun(I,R)$ the
pullback functors, and we denote by $\gamma_!,\gamma_*:\Fun(I,R)
\to \Fun(I',R)$ its left and right-adjoint. If $I'$ is the point
category, and $\gamma$ is the tautological projection, then by
definition, for any $E \in \D(I,R)$,
$$
L^\hdot\gamma_!E = H_\idot(I,E), \qquad R^\hdot\gamma_*E =
H^\hdot(I,E)
$$
are the homology resp. cohomology groups of the category $I$ with
coefficients in $E$.

To study cyclic homology, we will use A. Connes' cyclic category
$\Lambda$ of \cite{C} and its $l$-fold covers $\Lambda_l$, $l \geq
1$, of \cite[Appendix, A2]{FT}. The exact definitions can be found
for example in \cite[Section 1]{ka1}. We will need to know that for
any $l \geq 1$, objects in $\Lambda_l$ are numbered by non-negative
integers, with $[n] \in \Lambda_l$ being the object corresponding to
$n \geq 1$, and that we have natural functors
\begin{equation}\label{i.pi}
i_l,\pi_l:\Lambda_l \to \Lambda
\end{equation}
given by $i_l([n])=[nl]$, $\pi_l([n])=[n]$ on the level of
objects. If $l=1$, both functors are equivalences. Objects $[n] \in
\Lambda$ can be geometrically thought of as cellular decompositions
of the circle $S^1$; $0$-cells of such a decomposition are called
{\em vertices}. The set of vertices corresponding to $[n] \in
\Lambda$ is denoted $V([n])$, and this gives a functor $V$ from
$\Lambda$ to the category of finite sets. For any morphism $f:[n']
\to [n]$ and any vertex $v \in V([n])$, the set $f^{-1}(v) \subset
V([n'])$ carries a natural total order. For any $n,l \geq 1$, the
automorphism group $\Aut([n])$ of the object $[n] \in \Lambda_l$ is
naturally identified with the cyclic group $\Z/nl\Z$. Moreover, for
any $l \geq 1$, we have a natural functor
\begin{equation}\label{j.l}
j_l:\Delta^o \to \Lambda_l,
\end{equation}
where $\Delta^o$ is the opposite to the category of finite non-empty
totally ordered sets. Contrary to the standard usage, we denote by
$[n] \in \Delta$ the set with $n$ elements, so that
$j_l([n])=[n]$. For $l=1$, we simplify notation by writing $j=j_1$,
and for any $l \geq 1$, we have $\pi_l \circ j_l \cong j$. The
functor $j$ is an equivalence between $\Delta^o$ and the category of
objects $[n] \in \Lambda$ with a distinguished vertex $v \in
V([n])$.

Recall that for any simplicial $R$-module $M \in \Fun(\Delta^o,R)$,
the homology $H_\idot(\Delta^o,M)$ can be computed by the {\em
  standard complex} $CH_\idot(M)$ with terms $CH_i(M) = M([i+1])$,
$i \geq 0$, and the differential given by
\begin{equation}\label{d.del}
d = \sum_{0 \leq j \leq i}(-1)^jd^i_j,
\end{equation}
where $d^i_\idot:[i] \to [i+1]$ are the standard face maps. We can
also consider the complex $CH'_\idot(M)$ with the same terms
$CH_i'(M) = M([i+1])$ and with the differential
\begin{equation}\label{d.del.bis}
d' = \sum_{0 \leq j < i}(-1)^jd^i_j.
\end{equation}
This complex is canonically contractible. More generally, if we
extend the functors $CH_\idot(-)$ and $CH'_\idot(-)$ to complexes
as in \eqref{phi.c}, then for any complex $M_\idot$ of simplicial
$R$-modules, the complex $CH'_\idot(M_\idot)$ is acyclic, and the
complex $CH_\idot(M_\idot)$ computes the homology
$H_\idot(\Delta^o,M_\idot)$ of the category $\Delta^o$ with
coefficients in $M_\idot$. This homology is denoted
$HH_\idot(M_\idot)$.

Now, it is well-known and easy to check that if $M = j_l^*E$ for for
$l \geq 1$ and $E \in \Fun(\Lambda_l,R)$, then $d \circ (\id -
\sigma^\dg) = (\id - \sigma^\dg) \circ d'$, where in every degree
$i$, $\sigma$ is the generator of the cyclic group $Z/l(i+1)\Z \cong
\Aut([i+1])$, and $\sigma^\dg$ is its twist by the sign, as in
Subsection~\ref{cycl.grp.subs}. Therefore we have a natural map of
complexes
\begin{equation}\label{id.si.eq}
\id - \sigma^\dg:CH'_\idot(j_l^*E) \to CH_\idot(j_l^*E).
\end{equation}
In particular, if we denote
\begin{equation}\label{cc}
cc_i(E) = CH_i(j_l^*E)_{\sigma^\dg} = E([i+1])_{\sigma^\dg} =
\Coker(\id - \sigma^\dg), \qquad i \geq 0,
\end{equation}
then the differential $d$ of \eqref{d.del} descends to a
well-defined differential $d:cc_\idot(E) \to cc_\idot(E)$ turning
$cc_\idot(E)$ into a functorial complex. We call it the {\em reduced
  cyclic complex} of the object $E \in \Fun(\Lambda_l,R)$.

We note right away that if $l$ is odd, then for any $n \geq 1$,
$l(n+1)$ has the same parity as $(n+1)$, so that the sign twist in
the definition of $\sigma^\dg$ for the objects $[n] \in \Lambda_l$,
$[n] \in \Lambda$ is the same. Then for any object $E$ in
$\Fun(\Lambda_l,R)$, \eqref{cc} provides a functorial isomorphism
\begin{equation}\label{cc.pi.0}
cc_\idot(\pi_{l!}E) \cong cc_\idot(E),
\end{equation}
where $\pi_l$ is the functor \eqref{i.pi}.

To define the usual cyclic complex, we recall from \cite[Section
  1]{ka1} that the exact sequences \eqref{K.eq} for the groups
$\Aut([n])$, $[n] \in \Lambda_l$ fit together into a single exact
sequence
\begin{equation}\label{K.la.eq}
\begin{CD}
0 @>>> \Z @>{\kappa_1}>> \K_1 @>>> \K_0 @>{\kappa_0}>>
\Z @>>> 0
\end{CD}
\end{equation}
in the category $\Fun(\Lambda_l,\Z)$, with $\Z$ being the constant
functor with value $\Z$. Then for any $E \in \Fun(\Lambda_l,R)$, we
have a natural length-$2$ complex $\K_\idot(E) = \K_\idot \otimes E$
whose homology objects are $E$ in degree $0$ and $1$. Moreover,
equipping $\K_\idot(E)$ with the differential $B=\kappa_1 \circ
\kappa_0$ turns it into a mixed complex in
$\Fun(\Lambda_l,R)$. Applying the functor $cc_\idot$ and taking the
sum-total complex, we obtain a functorial mixed complex
$$
CH_\idot(E) = \tot(cc_\idot(\K_\idot(E)))
$$
in the category of $R$-modules.

By definition, $CH_\idot(E)$ is the cone of a map $cc_\idot(\K_1(E))
\to cc_\idot(\K_0(E))$ induced by the map $\K_1 \to
\K_0$. Explicitly, we have natural identifications
\begin{equation}\label{K.ch}
cc_\idot(\K_0(E)) \cong CH_\idot(j_l^*E), \qquad cc_\idot(\K_1(E))
\cong CH'_\idot(j_l^*E),
\end{equation}
and the map is exactly the map \eqref{id.si.eq}. Thus we have a
natural injective map
\begin{equation}\label{ch.ch}
CH_\idot(j_l^*E) \to CH_\idot(E),
\end{equation}
and the quotient $CH'_\idot(j^*_lE)$ is an acyclic complex, so that
\eqref{ch.ch} is a quasiisomorphism. The homology of the complex
$CH_\idot(E)$ is denoted $HH_\idot(E)$, and it coincides with
$HH_\idot(j_l^*E)$. The augmentation map $\kappa_0:\K_\idot(E) \to
E$ induces a natural map
\begin{equation}\label{gamma}
\gamma:CH_\idot(E) \to cc_\idot(E).
\end{equation}
In terms of the identifications \eqref{K.ch}, this map is given by
the projections from $CH_i(j_l^*E)$, $i \geq 0$, to their quotients
$cc_i(E)$ of \eqref{cc}.

\begin{defn}
For any $l \geq 1$ and any $E \in \Fun(\Lambda_l,R)$, the {\em
  cyclic complex} $CC_\idot(E)$ is given by
$$
CC_\idot(E) = \Exp(CH_\idot(E)),
$$
where $\Exp$ is the expansion of Definition~\ref{exp}. The homology
of the complex $CC_\idot(E)$ is denoted $HC_\idot(E)$.
\end{defn}

Explicitly, $CC_\idot(E)$ is the total complex of a bicomplex whose
even-numbered columns are all isomorphic to $CH_\idot(j_l^*E)$, and
whose odd-num\-be\-red columns are acyclic. For $l=1$, $CC_\idot(E)$
is the standard cyclic complex of a cyclic object $E \in
\Fun(\Lambda,R)$ (see e.g. \cite[Section 2.1.2]{Lo}). By definition,
we have the embedding $CH_\idot(E) \to CC_\idot(E)$, and it fits
into a short exact sequence of complexes
$$
\begin{CD}
0 @>>> CH_\idot(E) @>>> CC_\idot(E) @>{u}>> CC_\idot(E)[2] @>>> 0,
\end{CD}
$$
where $u$ is the periodicity map of Definition~\ref{exp}. The map
$\gamma$ of \eqref{gamma} extends to a natural augmentation map
\begin{equation}\label{aug.cc}
\alpha:CC_\idot(E_\idot) \to cc_\idot(E_\idot).
\end{equation}
If for any $[n] \in \Lambda_l$, $H_i(\Z/nl\Z,E([n])=0$ for $i \geq
1$, then this map is a quasiisomorphism (this is obvious from the
second spectral sequence for the bicomplex that gives
$CC_\idot(E)$).

More generally, we extend the functor $cc_\idot(-)$ to complexes as
in \eqref{phi.c}. Then for every complex $E_\idot$ in
$\Fun(\Lambda_l,R)$, we denote by $\K_\idot(E_\idot)$ the cone of
the natural map $\K_1 \otimes E_\idot \to \K_0 \otimes E_\idot$, and
we denote
\begin{equation}\label{ch.cc}
CH_\idot(E_\idot) = cc_\idot(\K_\idot(E_\idot)), \qquad
CC_\idot(E_\idot) = \Exp(CH_\idot(E_\idot)).
\end{equation}
We also denote by $HH_\idot(E_\idot)$ resp.\ $HC_\idot(E_\idot)$ the
homology of $CH_\idot(E_\idot)$ resp.\ $CC_\idot(E_\idot)$.
Alternatively, $CH_\idot(E_\idot)$ is the cone of the natural map
$$
\id - \sigma^\dg:CH'_\idot(j_l^*E_\idot) \to
CH_\idot(j_l^*E_\idot),
$$
and in terms of the identifications \eqref{K.ch}, this map is the
natural map induced by the map $\K_1 \to \K_0$. In particular, we
have $HH_\idot(E_\idot) \cong HH_\idot(j_l^*E_\idot)$, and for any
quasiisomorphism $E_\idot \to E'_\idot$, the corresponding map
\begin{equation}\label{ch.qi}
CH_\idot(E_\idot) \to CH_\idot(E'_\idot)
\end{equation}
of mixed complexes is a quasiisomorphism in the sense of
Definition~\ref{mix.qis}. We note that $\Exp$ commutes with the
totalization $\tot$, so that we have $CC_\idot(E_\idot) \cong
cc_\idot(\Exp(\K_\idot(E_\idot)))$. We also have natural
augmentation maps $\gamma$, $\alpha$ of \eqref{gamma}
resp. \eqref{aug.cc}.

\subsection{Periodic complexes -- definitions.}\label{per.co.subs}

We now assume given a complex $E_\idot$ in the category
$\Fun(\Lambda_l,R)$, take the mixed complex $CH_\idot(E_\idot)$, and
consider the periodic expansions of Definition~\ref{per}.

\begin{defn}\label{cp.def}
Assume given a complex $E_\idot$ in the category
$\Fun(\Lambda_l,R)$. Then the {\em periodic cyclic complex}
$CP_\idot(E_\idot)$ and the {\em co-periodic cyclic complex}
$\bCP_\idot(E_\idot)$ are given by
\begin{equation}\label{cp.eq}
CP_\idot(E_\idot) = \Per(CH_\idot(E_\idot)), \qquad
\bCP_\idot(E_\idot) = \cPer(CH_\idot(E_\idot)).
\end{equation}
The homology of the complexes $CP_\idot(E_\idot)$ resp.\
$\bCP_\idot(E_\idot)$ are denoted by $HP_\idot(E_\idot)$
resp.\ $\bHP_\idot(E_\idot)$.
\end{defn}

The complexes $CP_\idot(E_\idot)$ and $\bCP_\idot(E_\idot)$ and
their homology are the main objects of study in this
paper. Unfortunately, there is no map from one complex to the other
one. To be able to compare the two, we have to consider an
additional functorial complex.

\begin{defn}\label{cp.bis.def}
  For any complex $E_\idot$ in $\Fun(\Lambda_l,R)$, the {\em
    polynomial periodic cyclic complex} $cp_\idot(E_\idot)$ is given
  by
$$
cp_\idot(E_\idot) = \per(CH_\idot(E_\idot)) \cong
cc_\idot(\per(\K_\idot(E_\idot))).
$$
\end{defn}

Then for any complex $E_\idot$ in $\Fun(\Lambda_l,R)$, the natural
maps \eqref{per.Per} induce natural maps
\begin{equation}\label{CP.cp}
\begin{CD}
CP_\idot(E_\idot) @<<< cp_\idot(E_\idot) @>>> \bCP_\idot(E_\idot)
\end{CD}
\end{equation}
from the polynomial periodic cyclic complex $cp_\idot(E_\idot)$ to
the periodic and co-periodic cyclic complexes. In general, neither
of these maps is a quasiisomorphism. Moreover, for any
quasiisomorphism $f:E_\idot \to E'_\idot$, we can consider the
corresponding maps
\begin{equation}\label{cp.qis}
CP_\idot(E_\idot) \to CP_\idot(E'_\idot), \ 
\bCP_\idot(E_\idot) \to \bCP_\idot(E'_\idot), \ 
cp_\idot(E_\idot) \to cp_\idot(E'_\idot).
\end{equation}
induced by the map \eqref{ch.qi}. Then the first of these maps is a
quasiisomorphism by Lemma~\ref{mix.qis.le}, so that $CP_\idot(-)$
descends to a functor on the derived category
$\D(\Lambda_l,R)$. However, this need not be true for the other two
maps, so that neither $cp_\idot(-)$ nor $\bCP_\idot(-)$ are defined
on the level of the derived category, and one cannot study them by
the standard techniques of homological algebra.

\medskip

One way out of this difficulty is to notice that for any abelian
category $\C$, complexes in $\C$ strongly acyclic in the sense of
Definition~\ref{strong.qis} form a triangulated subcategory in the
homotopy category of chain complexes. One can then invert strong
quasiisomorphisms in $C_\idot(\C)$ and obtain the so-called {\em
  absolute derived category} $\D_{abs}(\C)$. This a triangulated
category introduced and studied extensively by L. Positselski
\cite{posic}. If a map $E_\idot \to E'_\idot$ is a strong
quasiisomorphism, then all the maps \eqref{cp.qis} are
quasiisomorphisms, so that $cp_\idot(-)$ and $\bCP_\idot(-)$ do
descend to the co-derived category $\D_{abs}(\Lambda_l,R)$ of the
abelian category $\Fun(\Lambda_l,R)$. However, absolute derived
categories are at present not well enough understood for our
purposes. So, for most of our results, we adopt a more conventional
approach -- we equip everything in sight with a filtration, and use
filtered derived categories described in Subsection~\ref{DF.subs}.

\medskip

We recall that the cyclic complex functor $cc_\idot(-)$ can be
promoted to a functor $cc^f$ of \eqref{phi.f.b} from the category
of filtered complexes in $\Fun(\Lambda_l,R)$ with a termwise-split
filtration to the category of filtered complexes of $R$-modules. We
will need to modify this slightly.

\begin{defn}\label{stand.filt.def}
For any filtered complex $\langle E_\idot,F^\hdot \rangle$ in
$\Fun(\Lambda_l,R)$ with term\-wise-split filtration, the {\em
  standard filtration} $F^\hdot$ on the complex $cc_\idot(E_\idot)$
is given by
\begin{equation}\label{cc.fi}
  F^n cc_\idot(E_\idot) = F^{n-1}cc^f(E_\idot),
\end{equation}
where the filtered complex $cc^f(E_\idot)$ is as in \eqref{phi.f}.
\end{defn}

In other words, we shift the filtration on $cc^f(E_\idot)$ by $1$ in
the sense of \eqref{scale}. From now on, for any filtered complex
$E_\idot$ in $\Fun(\Lambda_l,R)$ with a termwise-split filtration,
we will write simply $cc_\idot(E_\idot)$ to mean $cc_\idot(E_\idot)$
equipped with the standard filtration \eqref{cc.fi}.

For any filtered complex $E_\idot$ with termwise-split filtration,
the filtration on $E_\idot$ induces a termwise-split filtration on
$\K_\idot(E_\idot) = \K_\idot \otimes E_\idot$, and therefore the
standard filtration \eqref{cc.fi} induces a filtration on
$CH_\idot(E_\idot) = cc_\idot(\K_\idot(E_\idot))$. This filtration
is preserved by the differential $B$ in the mixed complex
$CH_\idot(E_\idot)$, thus induces filtrations \eqref{mixed.filt} on
$CP_\idot(E)$, $\bCP_\idot(E)$ and $cp_\idot(E)$. We will also call
them the {\em standard filtrations}.

Spelling out Definition~\ref{stand.filt.def}, we see that for any
complex $E_\idot$ in $\Fun(\Lambda_l,R)$ equipped with a
termwise-split filtration $F^\hdot$, the associated graded quotient
of the mixed complex $CH_\idot(E_\idot)$ with respect to the
standard filtration is given by
\begin{equation}\label{gr.stand.ch}
\gr_F^iCH_\idot(E_\idot) \cong
\bigoplus_{m \geq 1}\K_\idot(\gr^{i-m}_FE_\idot([m]))_{\sigma^\dg}[m-1]
\end{equation}
for any integer $i$, where $\K_\idot(-)$ is the mixed complex
\eqref{K.grp}, and $[m-1]$ in the right-hands side stands for the
homological shift. The associated graded quotients of the periodic
comlexes are then given by \eqref{mixed.filt.gr}. In particular,
since the polynomial expansion functor commutes with arbitrary
direct sums, we have
\begin{equation}\label{gr.stand.cp}
\gr_F^icp_\idot(E_\idot) = \bigoplus_{m \geq
  1}\wC^\dg_\idot(\Z/lm\Z,\gr^{i-m}_FE_\idot([m]))[m-1],
\end{equation}
where $\wC^\dg_\idot(-)$ is the twisted Tate homology complex of
\eqref{tate.tw} with respect to the cyclic group $\Z/ml\Z =
\Aut([m])$. Moreover, if $E_\idot$ is concentrated in a finite range
of filtered degrees, then the direct sum in \eqref{gr.stand.ch} is
in fact finite, so that we also have
\begin{equation}\label{gr.stand.CP}
\begin{aligned}
\gr_F^iCP_\idot(E_\idot) &= \bigoplus_{m \geq
  1}\vC^\dg_\idot(\Z/lm\Z,\gr^{i-m}_FE_\idot([m]))[m-1],\\
\gr_F^i\bCP_\idot(E_\idot) &= \bigoplus_{m \geq
  1}\vbC^\dg_\idot(\Z/lm\Z,\gr^{i-m}_FE_\idot([m]))[m-1],
\end{aligned}
\end{equation}
where $\vC^\dg_\idot(-)$ and $\vbC^\dg_\idot(-)$ are again the
complexes of \eqref{tate.tw}.

\subsection{Periodic complexes -- first properties.}

It turns out that standard filtrations are already useful when we
take a complex $E_\idot$ in $\Fun(\Lambda_l,R)$ and treat it as a
filtered complex concentrated in filtered degree $0$. Then by
definition, the standard filtrations on $cp_\idot(E_\idot)$,
$CP_\idot(E_\idot)$, $\bCP_\idot(E_\idot)$ are bounded below
(actually, already the term $F^0$ vanishes). However, while the
standard filtration on $cp_\idot(E_\idot)$ is exhaustive, the same
need not be true for $CP_\idot(E_\idot)$ or
$\bCP_\idot(E_\idot)$. In fact, let us introduce the following.

\begin{defn}\label{restr.def}
  For any complex $E_\idot$ in $\Fun(\Lambda_l,R)$, the {\em
    restricted} periodic resp.\ co-periodic complexes
  $CP^f_\idot(E_\idot)$, $\bCP^f_\idot(E_\idot)$ are given by
$$
CP^f_\idot(E_\idot) = cc_\idot(\Per(\K_\idot(E_\idot))),
\quad
\bCP^f_\idot(E_\idot) = cc_\idot(\cPer(\K_\idot(E_\idot))).
$$
\end{defn}

Then $CP^f_\idot(E_\idot)$, $\bCP^f_\idot(E_\idot)$ are precisely
the completions of $CP_\idot(E_\idot)$ resp. $\bCP_\idot(E_\idot)$
with respect to the standard filtrations in the sense of
Definition~\ref{compl.def}. Alternatively, up to a shift of
filtration, they coincide with the extensions \eqref{phi.f} of the
functors $CP_\idot(-)$, $\bCP_\idot(-)$ from $\Fun(\Lambda_l,R)$ to
complexes of $R$-modules, so that the notation is consistent. Since
the filtration on $E_\idot$ is bounded below, the standard
filtrations are also bounded below, and we have natural functorial
maps
\begin{equation}\label{cpf.cp}
CP^f_\idot(E_\idot) \to CP_\idot(E_\idot), \qquad
\bCP^f_\idot(E_\idot) \to \bCP_\idot(E_\idot).
\end{equation}
Together with the natural maps
$$
\begin{CD}
CP^f_\idot(E_\idot) @<<< cp_\idot(E_\idot) @>>>
\bCP^f_\idot(E_\idot)
\end{CD}
$$
obtained by applying $cc_\idot$ to the maps \eqref{per.Per}, the
maps \eqref{cpf.cp} fit into the following diagram
\begin{equation}\label{5.dia}
\begin{CD}
cp_\idot(E_\idot) @>{r}>> \bCP^f_\idot(E_\idot) @>{R}>>
\bCP_\idot(E_\idot)\\
@|\\
cp_\idot(E_\idot) @>{l}>> CP^f_\idot(E_\idot) @>{L}>>
CP_\idot(E_\idot),
\end{CD}
\end{equation}
a refinement of \eqref{CP.cp}. All the maps in this diagram are
injective maps of complexes, and in general, none of them are
quasiisomorphisms.

By \eqref{gr.stand.CP} and Lemma~\ref{mix.qis.le}, for any
quasiisomorphism $E_\idot \to E'_\idot$ of complexes in
$\Fun(\Lambda_l,R)$, the first map in \eqref{cp.qis} is a filtered
quasiisomorphism with respect to the standard filtration, so that
the corresponding map
$$
CP^f_\idot(E_\idot) \to CP^f_\idot(E'_\idot)
$$
is a quasiisomorphism. Therefore just as $CP_\idot(-)$, the complex
$CP^f_\idot(-)$ is defined on the level of the derived category
$\D(\Lambda_l,R)$, and so is the map $L$
of \eqref{5.dia}. Note that the map need not be an isomorphism even
in the simplest examples such as the following.

\begin{lemma}\label{cp.const.le}
  Let $E_\idot$ be the constant functor $R$ placed in homological
  degree $0$. Then $CP_\idot(R)$ is quasiisomorphic to the free
  rank-$1$ module over the Laurent power series ring $R((u))$, while
  $CP^f_\idot(R)$ is quasiisomorphic to $(R \lotimes_{\Z}
  (\Q/\Z))((u))[1]$. The map $L$ of \eqref{5.dia} fits into a
  distinguished triangle
$$
\begin{CD}
CP^f_\idot(R) @>{L}>> CP_\idot(R) @>>> CP_\idot(R \otimes_{\Z} \Q) @>>>
\end{CD}
$$
\end{lemma}

\proof{} Let $E_i=R$, $i \geq 0$, let $e_i:E_i \to E_{i+1}$ be
given by
$$
e_i = \begin{cases} (i+1)\id, &\quad i=2n,\\
2\id, &\quad i=2n+1,
\end{cases}
$$
and consider the bicomplex $E_{\idot,\idot}$ of
Lemma~\ref{snake.lim.le}. Then one immediately checks that the
periodic bicomplex $\per_{\idot,\idot}(CH_\idot(R))$ of
\eqref{per.bi} coincides on the nose with $E_{\idot,\idot}((u))$. To
finish the proof, apply Lemma~\ref{snake.lim.le}.
\endproof

As for the other functorial complexes of \eqref{5.dia}, then in the
simple situation such as that of Lemma~\ref{cp.const.le}, they
contain no new information. Indeed, one checks easily that for any
complex $E_\idot$ concentrated in a finite range of homological
degrees, the maps $r$, $l$, and $R$ are isomorphisms, so that the
only non-trivial invariants of $E_\idot$ are $CP_\idot(E_\idot)$ and
$CP^f_\idot(E_\idot) \cong cp_\idot(E_\idot) \cong \bCP^f(E_\idot)
\cong \bCP_\idot(E_\idot)$.

For unbounded complexes, the situation is more difficult, and
neither of the maps is even a quasiisomorphism without additional
assumptions. The only obvious general fact is the following result.

\begin{lemma}\label{cpf.Q.le}
  For any complex $E_\idot$ in $\Fun(\Lambda_l,R)$, we have
  quasiisomorphisms $CP^f_\idot(E_\idot) \otimes \Q \cong
  cp_\idot(E_\idot) \otimes \Q \cong \bCP_\idot^f(E_\idot) \otimes \Q
  \cong 0$.
\end{lemma}

\proof{} Since all the complexes are complete with respect to the
standard filtration, it suffices to check that all the associated
graded quotients \eqref{gr.stand.cp}, \eqref{gr.stand.CP} are
acyclic. But once we tensor with $\Q$, the differential in the Tate
homology complex becomes chain-homotopic to $0$.
\endproof

To analyse $\bCP_\idot(E_\idot)$, we can also use the standard
filtrations, but in a different way: take a complex $E_\idot$ in
$\Fun(\Lambda_l,R)$, and equip it with the stupid filtration instead
of the trivial one (it is automatically termwise-split). Then we
have the following result.

\begin{lemma}\label{stu.cp.le}
Assume that a complex $E_\idot$ in $\Fun(\Lambda_l,R)$ is equipped
with a stupid filtration. Then the map
$$
R \circ r:cp_\idot(E_\idot) \to \bCP_\idot(E_\idot)
$$
is a filtered quasiisomorphism with respect to the standard
filtrations, and $\bCP_\idot(E_\idot)$ is complete.
\end{lemma}

\proof{} In the case of the stupid filtration, \eqref{gr.stand.ch}
provides an isomorphism
\begin{equation}\label{gr.stu.cp}
\gr_F^iCH_\idot(E_\idot) \cong \bigoplus_{m \geq
  1}\K_\idot(E^{i-m}([m]))_{\sigma^\dg}[i-1]
\end{equation}
for any integer $i$. In particular, for any $i$, $\gr_F^i
CH_\idot(E_\idot)$ is concentrated in a finite range of
cohomological degrees (degrees $-i$ and $1-i$, to be exact). Then $R
\circ r$ is a filtered quasiisomorphism by
\eqref{mixed.filt.gr}. Moreover, the completion of the complex
$cp_\idot(E_\idot)$ can be computed as in Example~\ref{sum.exa}, and
it clearly coincides with $\bCP_\idot(E_\idot)$.
\endproof

\begin{corr}\label{cpf.Q.corr}
For any complex $E_\idot$ in $\Fun(\Lambda,R)$, the complex
$\bCP_\idot(E_\idot) \otimes \Q$ is acyclic.
\end{corr}

\proof{} Equip $E_\idot$ with the stupid filtration as in
Lemma~\ref{stu.cp.le}, so that the tensor product
$\bCP_\idot(E_\idot) \otimes \Q \cong \bCP_\idot(E_\idot \otimes
\Q)$ is complete with respect to the standard filtration. Then it
suffices to check that its associated graded quotient is acyclic. By
\eqref{gr.stu.cp}, this associated graded quotient is a sum of Tate
cohomology complexes, and as in Lemma~\ref{cpf.Q.le}, these become
acyclic after tensoring with $\Q$.
\endproof

\subsection{Representing objects.}\label{repr.subs}

Assume given a small category $I$ and a ring $R$, and denote by
$I^o$ the opposite category. Then for any objects $E \in \Fun(I,R)$,
$M \in \Fun(I^o,\Z)$, the {\em tensor product} $E \otimes_I M$ is the
cokernel of the natural map
$$
\begin{CD}
\displaystyle\bigoplus_{f:i \to i'}E(i) \otimes M(i')
@>{E(f) \otimes \id - \id \otimes M(i')}>>
\displaystyle\bigoplus_{i \in I}E(i) \otimes M(i),
\end{CD}
$$
where the sum in the right-hand side is over all objects $i \in I$,
and the sum in the left-hand side is over all morphisms $f:i \to i'$
in $I$. This tensor product and its derived functor $-
\lotimes_I -$ are very convenient for representing functors from
$\Fun(I,R)$ to $R$-modules or complexes of $R$-modules. For example,
for any $E \in \Fun(I,R)$, we have a natural identification
\begin{equation}\label{ho.ten}
E \lotimes_I \Z \cong H_\idot(I,E),
\end{equation}
where $\Z \in \Fun(I^o,\Z)$ is the constant functor. On the other
hand, for any object $i \in I$, we have a natural Yoneda-type
identification
\begin{equation}\label{yo.ten}
E(i) \cong E \otimes_I \Z_i \cong E \lotimes_I \Z_i,
\end{equation}
where $\Z_i \in \Fun(I^o,\Z)$ is the representable functor given by
$$
\Z_i(i') = \Z[I(i',i)],
$$
with $I(i',i)$ standing for the set of morphisms $f:i' \to i$ in
$I$. Thus to compute homology $H_\idot(I,E)$ by an explicit complex
functorial in $E$, it suffices to find a resolution of the constant
functor $\Z \in \Fun(I^o,\Z)$ by sums of objects $\Z_i$, $i \in I$.

Applying this formalism to cyclic homology, we can represent the
functorial complexes $cc_\idot(-)$, $CH_\idot(-)$, $CC_\idot(-)$ of
Subsection~\ref{cp.subs}. Namely, fix an integer $l \geq 1$, for any
$i \geq 0$, let
$$
Q_i = \Z_{[i+1]} \in \Fun(\Lambda_l^o,\Z),
$$
and let $q_i = (Q_i)_{\sigma^\dg}$ be cokernel of the natural map
$\id - \sigma^\dg:Q_i \to Q_i$, where $\sigma$ is the generator of
the cyclic group $\Z/(i+1)l\Z = \Aut([i+1])$, and $\sigma^\dg$ is as
in Subsection~\ref{cycl.grp.subs}. Then \eqref{d.del} defines a
differential $d:Q_{\idot+1} \to Q_\idot$ turning $Q_\idot$ into a
complex, and it descends to a differential $d:q_{\idot+1} \to
q_\idot$. By \eqref{yo.ten}, we then have
$$
cc_\idot(E) \cong (E \otimes_{\Lambda_l} Q_\idot)_\sigma \cong E
\otimes_{\Lambda_l} q_\idot
$$
for any object $E \in \Fun(\Lambda_l,R)$. Analogously,
\eqref{d.del.bis} defines a differential $d':Q_{\idot + 1} \to
Q_\idot$. Denote the resulting complex by $Q'_\idot$, and denote by
$P_\idot$ the cone of the natural map of complexes
\begin{equation}\label{q.qp}
Q'_\idot \to Q_\idot
\end{equation}
given by $(\id - \sigma^\dg)$ termwise. Then $P_\idot$ is a mixed
complex in $\Fun(\Lambda_l^o,\Z)$, and we have
$$
CH_\idot(E) \cong E \otimes_{\Lambda_l} P_\idot
$$
for any $E \in \Fun(\Lambda_l,R)$.

\begin{lemma}\label{P.mix.le}
The mixed complex $P_\idot$ is a mixed resolution of the constant
functor $\Z \in \Fun(\Lambda_l^o,\Z)$ in the sense of
Definition~\ref{mix.reso}.
\end{lemma}

\proof{} The statement is essentially \cite[Appendix, A3]{FT}, but
we reproduce the proof for the convenience of the reader. By
\eqref{yo.ten}, it suffices to prove that for any $[n] \in
\Lambda_l$, the mixed complex $CH_\idot(\Z_{[n]})$ is a mixed
resolution of the abelian group $\Z$. Fix $[n]$, and denote by
$E_\idot$, $E'_\idot$ the cokernel resp. kernel of the natural map
$$
\begin{CD}
CH'_\idot(j_l^*\Z_{[n]}) @>{\id - \sigma^\dg}>> CH_\idot(j_l^*\Z_{[n]})
\end{CD}
$$
whose cone is the complex $CH_\idot(\Z_{[n]})$. Then for any $[m]
\in \Lambda_l$, the action of the cyclic group $\Z/ml\Z = \Aut([m])$
on the set of morphisms $\Lambda_l([n],[m])$ is free. Therefore for
any $i \geq 0$, the differential $B$ in the mixed complex
$CH_\idot(\Z_{[n]})$ is an isomorphism between the cokernel $E_i$
and the kernel $E'_i$ of the endomorphism $\id - \sigma^\dg$ of the
$\Z[\Z/l(i+1)\Z]$-module $\Z_{[n]}([i+1]) =
\Z(\Lambda_l([n],[i+1]))$. We conclude that $B$ factors through an
isomorhism $E_\idot \to E'_\idot$, so that it suffices to prove
that $E_\idot$ is a resolution of $\Z$. But $E_\idot$ is precisely
the standard chain complex of the elementary simplex $\Delta^o_{[n]}
\in \Delta^o\Sets$.
\endproof

As a corollary of Lemma~\ref{P.mix.le}, we see that the expansion
$\Exp(P_\idot)$ of the complex $P_\idot$ is a resolution of the
constant functor $\Z \in \Fun(\Lambda_l,\Z)$, so that by
\eqref{ho.ten}, for any $E \in \Fun(\Lambda_l,R)$, the complex
\begin{equation}\label{cc.p}
CC_\idot(E) = \Exp(CH_\idot(E)) \cong \Exp(E \otimes_{\Lambda_l}
P_\idot) \cong E \otimes_{\Lambda_l} \Exp(P_\idot)
\end{equation}
computes the homology $H_\idot(\Lambda_l,R)$ of the category
$\Lambda_l$ with coefficients in $E$, and more generally, we have
$HC_\idot(E_\idot) \cong H_\idot(\Lambda_l,E_\idot)$ for any complex
$E_\idot$ in $\Fun(\Lambda_l,R)$ (this is of course well-known, see
e.g. \cite[Appendix, Corollary A3.2]{FT}). Lemma~\ref{cp.const.le}
then implies another well-known fact, namely, that the cohomology
$H^\hdot(\Lambda_l,\Z)$ is the free algebra $\Z[u]$ in one generator
$u$ of degree $2$, and the cohomology class of the mixed resolution
$P_\idot$ is exactly $u$.

For some applications of the tensor product formalism, one does not
even need to know the exact shape of the representing complexes. For
example, we will need the following easy observation.

\begin{lemma}\label{c.cc.le}
For any complex $E_\idot$ in $\Fun(\Lambda_l,R)$, the map
\begin{equation}\label{alpha.cmpl}
\alpha:CC_\idot(\K_\idot(E_\idot)) \to cc_\idot(\K_\idot(E_\idot))
\cong CH_\idot(E_\idot)
\end{equation}
induced by the augmentation map \eqref{aug.cc} is a
quasiisomorphism.
\end{lemma}

\proof{} The functor $CH_\idot(-)$ is represented by the complex
$P_\idot$ of Lemma~\ref{P.mix.le}, and the functor
$CC_\idot(\K_\idot(E_\idot))$ is also obviously represented by a
complex $M_\idot$ of projective objects in
$\Fun(\Lambda_l^o,\Z)$. Both $P_\idot$ and $M_\idot$ are bounded
from above (by $0$, but this is not important). The map
\eqref{alpha.cmpl} is induced by a map of complexes
$$
a:M_\idot \to P_\idot.
$$
Since for any $[n] \in \Lambda_l$, both $\K_0([n])$ and $\K_1([n])$
are free $\Z[\Z/nl\Z]$-modules, $H_i(\Z/nl\Z,\K_0(E)) =
H_i(\Z/nl,\K_1(E)) = 0$ for any $i \geq 1$ and any object $E$ in
$\Fun(\Lambda_l,R)$. Therefore the map \eqref{alpha.cmpl} is a
quasiisomorphism when $E_\idot$ is concentrated in a single
homological degree. By \eqref{yo.ten}, this means that $a$ is a
quasiisomorphism after evaluating at any object $[n] \in
\Lambda^o_l$, thus a quasiisomorphism of complexes in
$\Fun(\Lambda_l^o,\Z)$. But since both $M_\idot$ and $P_\idot$ are
complexes of projective objects bounded from above, $a$ must then be
a chain-homotopy equivalence. Therefore the map \eqref{alpha.cmpl}
is also a chain-homotopy equivalence, thus a quasiisomorphism for
any $E_\idot$.
\endproof

Observe now that in fact, both the differentials \eqref{d.del} and
\eqref{d.del.bis} only involve the standard face maps -- that is,
injective maps in $\Delta$. These correspond to surjective maps in
$\Delta^o$ and in $\Lambda_l$. Therefore if we denote by $\bLambda_l
\subset \Lambda_l$ the subcategory of surjective maps, and let
$e:\bLambda_l \to \Lambda_l$ be the natural embedding, then for any
$E \in \Fun(\Lambda_l,R)$, the mixed complex $CH_\idot(E)$ only
depends on the restriction $e^*E \in \Fun(\bLambda_l,R)$. In terms
of representing objects, denote
$$
\bQ_i = \Z_{[i+1]} \in \Fun(\bLambda_l^o,\Z)
$$
for any $i \geq 0$. Then we have $Q_i \cong e_!\bQ_i$, the differentials
$d,d':Q_{\idot+1} \to Q_\idot$ are induced by differentials
$d,d':\bQ_{\idot+1} \to \bQ_\idot$, and we have
$$
Q_\idot \cong e_!\bQ_\idot, \qquad Q'_\idot \cong e_!\bQ_\idot
$$
for some canonical complexes $\bQ_\idot$, $\bQ'_\idot$ in
$\Fun(\bLambda_l^o,\Z)$. Moreover, since $\sigma \in \Aut([m])$ is
surjective for any $[m] \in \Lambda_l$, the map \eqref{q.qp} is
induced by a map $\id-\sigma^\dg:\bQ'_\idot \to \bQ_\idot$, and we have
$$
P_\idot \cong e_!\bP_\idot,
$$
where $\bP_\idot$ is the cone of $\id-\sigma$.

\begin{lemma}\label{P.mix.bis.le}
The complex $\bP_\idot$ is a mixed resolution of the constant
functor $\Z \in \Fun(\bLambda_l^o,\Z)$ in the sense of
Definition~\ref{mix.reso}.
\end{lemma}

\proof{} Literally the same as Lemma~\ref{P.mix.le}, except that the
complex $E_\idot$ is not the standard chain complex of the
elementary simplex $\Delta^o_{[n]} \in \Delta^o\Sets$ but rather,
its normalized chain complex.
\endproof

\section{Projections and subdivisions.}\label{edge.sec}

\subsection{Edgewise subdivision.}

Fix an integer $l \geq 1$, and consider the functor $i_l$ of
\eqref{i.pi}. Recall the following result.

\begin{lemma}\label{edge.le}
For any ring $R$ and any $E \in \D(\Lambda,R)$, the natural map
$$
H_\idot(\Lambda_l,i_l^*E) \to H_\idot(\Lambda,E)
$$
induced by the functor $i_l$ is an isomorphism.\endproof
\end{lemma}

This result is known as ``edgewise subdivision'' and goes back to
Segal and Quillen (a short proof with exactly the same notation as
in this paper can be found in \cite[Lemma 1.14]{ka0}).

The homology groups $H_\idot(\Lambda,-)$, $H_\idot(\Lambda_l,-)$ can
be computed by the cyclic complexes $CC_\idot(-)$, and it is natural
to ask whether one can lift the edgewise subdivision isomorphism to
a map of complexes. To do this, one can use the representing
complexes $P_\idot$ of Subsection~\ref{repr.subs}. By
Lemma~\ref{P.mix.le}, both $P_\idot$ and $i_l^*P_\idot$ are mixed
resolutions of the constant functor $\Z \in
\Fun(\Lambda_l^o,\Z)$. Moreover, Lemma~\ref{edge.le} shows that the
natural map $i_l^*:H^2(\Lambda^o,\Z) \to H^2(\Lambda_l^o,\Z)$ is an
isomorphism, so that both resolutions have the same cohomology class
$u=i_l^*u$. Since $P_i$ is a projective object for any $i$, we can
apply Lemma~\ref{mix.map.le} and obtain a quasiisomorphism
$$
\nu_l:\Exp(P_\idot) \to i_l^*\Exp(P_\idot).
$$
Applying \eqref{cc.p}, for any complex $E_\idot$ in
$\Fun(\Lambda,R)$, we obtain a functorial map
\begin{equation}\label{a.l.cc}
\nu_l:CC_\idot(i_l^*E_\idot) \to CC_\idot(E_\idot)
\end{equation}
realizing the edgewise subdivision isomorphism on the chain level.

Moreover, we can do this construction using the complexes
$\bP_\idot$ instead of the complexes $P_\idot$ and
Lemma~\ref{P.mix.bis.le} instead of Lemma~\ref{P.mix.le}. Indeed,
the functor $i_l$ induces a functor $\bi_l:\bLambda_l \to \bLambda$
such that $e \circ \bi_l \cong i_l \circ e$, and since the
adjunction map
$$
\bP_\idot \to e^*e_!\bP_\idot \cong e^*P_\idot
$$
is a quasiisomorphism, $\bP_\idot$ and $\bi_l^*\bP_\idot$ are also
mixed resolutions of the constant functor $\Z$ with the same
cohomology class. Therefore Lemma~\ref{mix.map.le} provides a
natural map
\begin{equation}\label{a.l}
\nu_l:\Exp(\bP_\idot) \to i_l^*\Exp(\bP_\idot).
\end{equation}
This map induces a map
$$
\nu_l:CC_\idot(\bi_l^*\overline{E}_\idot) \to
CC_\idot(\overline{E}_\idot)
$$
for any complex $\overline{E}_\idot$ in $\Fun(\bLambda,R)$. In
general, this map is {\em not} a quasiisomorphism; however, if
$\overline{E}_\idot = e^*E_\idot$ for some complex $E_\idot$ in
$\Fun(\Lambda,R)$, it is a quasiisomorphism since it coincides with
the map \eqref{a.l.cc}.

Using $\bP_\idot$ instead of $P_\idot$ has the following
advantage. By \eqref{per.exp}, a map $\nu$ provided by
Lemma~\ref{mix.map.le} uniquely extends to a $u$-equivariant map
$\nu:\Per(P_\idot) \to \Per(M_\idot)$. In particular, the map
\eqref{a.l.cc} extends to a map
\begin{equation}\label{a.l.cp}
\nu_l:CP_\idot(i_l^*E_\idot) \to CP_\idot(E_\idot)
\end{equation}
for any complex $E_\idot$ in $\Fun(\Lambda,R)$.

\begin{lemma}\label{edge.filt.le}
Assume given a complex $E_\idot$ in $\Fun(\Lambda,R)$ equipped with
a termwise-split filtration, and assume that the map \eqref{a.l.cc}
is induced by a map \eqref{a.l}. Then the corresponding map $\nu_l$
of \eqref{a.l.cp} sends $cp_\idot(i_l^*E_\idot) \subset
CP_\idot(i_l^*E_\idot)$ into $cp_\idot(E_\idot) \subset
CP_\idot(E_\idot)$, and induces a filtered map
\begin{equation}\label{nu.l.filt}
\nu_l:cp_\idot(i_l^*E_\idot^{[l]}) \to cp_\idot(E_\idot)^{[l]},
\end{equation}
where the complexes $cp_\idot(-)$ are equipped with standard
filtrations of Subsection~\ref{per.co.subs}, and $[l]$ stand for the
filtrations rescaled by $l$, as in \eqref{scale}.
\end{lemma}

\proof{} It obviously suffices to prove both claims when $E_\idot$
is a single object $E$ concentrated in a single filtered degree, say
degree $n$. Let $m$ be the integer such that $lm \leq n < l(m+1)$,
so that $E^{[l]}$ is concentrated in the filtered degree $m$. Then
the standard filtration on $CP_\idot(i_l^*E_\idot)$ is the shift by
$m$ of the filtration induced by the filtrations on $\bQ_\idot$,
$\bQ'_\idot$ assigning filtered degree $-(i+1)$ to
$\bQ_i$. Analogously, the standard filtration on $CP_\idot$ is
obtained by assigning filtered degree $-(i+1)$ to $\bQ_i$ and
shifting by $n$. Then to prove the second claim, it suffices to show
that $\Hom(\bQ_i,i_l^*\bQ_j)=0$ for $j+1-n > l(i+1-m)$. Since $n
\geq lm$, it suffices to let $n=m=0$. Then the claim is clear -- by
Yoneda Lemma, we have
$$
\Hom(\bQ_i,i_l^*\bQ_j) \cong i_l^*\bQ_j([i+1]) \cong \bQ_j([l(i+1)])
\cong \Z[\bLambda([l(i+1)],[j+1])],
$$
and the set $\bLambda([l(i+1)],[j+1])$ is empty unless $j+1 \leq
l(i+1)$. The first claim immediately follows from the second, since
$cp_\idot(E) \cong CP^f(E)$ is then the filtered completion of
$CP_\idot(E)$, and analogously for $cp_\idot(i_l^*E)$.
\endproof

Note that since all the rescalings of a given filtration are
commensurable, Lemma~\ref{stu.cp.le} shows that if we equip a complex
$E_\idot$ with the stupid filtration, the completion of the target
of the map \eqref{nu.l.filt} coincides with
$\bCP_\idot(E_\idot)$. However, the induced filtration on
$i_l^*E_\idot^{[l]}$ is not the stupid filtration because of the
rescaling. To adjust for this, we introduce the following.

\begin{defn}\label{cpv.def}
For any complex $E_\idot$ in $\Fun(\Lambda_l,R)$ and any integer
$n$, the {\em complex $\bCP^{[n]}(E_\idot)$} is the completion of
the complex $cp_\idot(E_\idot)$ with respect to the standard
filtration of Definition~\ref{stand.filt.def} corresponding to the
$n$-th rescaling $F^\hdot_{[n]}$ of the stupid filtration on
$E_\idot$. The homology of the complex $\bCP^{[n]}(E_\idot)$ is
denoted $\bHP^{[n]}(E_\idot)$.
\end{defn}

Then if $n=1$, $\bCP^{[1]}_\idot(E_\idot)$ coincides with
$\bCP_\idot(E_\idot)$ by Lemma~\ref{stu.cp.le}, but for $n > 1$ this
need not be the case. Nevertheless, the map \eqref{nu.l.filt}
induces a map
\begin{equation}\label{bcp.nu}
\nu_l:\bCP^{[l]}(i_l^* E_\idot) \to \bCP(E_\idot)
\end{equation}
for any complex $E_\idot$ in $\Fun(\Lambda,R)$.

\subsection{Filtered refinement.}\label{edge.filt.subs}

Assume now that $l=p$ is an odd prime that annihilates our base ring
$R$. Then we have the following filtered counterpart of
Lemma~\ref{edge.le}.

\begin{prop}\label{edge.prop}
Assume that $l=p$ is an odd prime such that $pR=0$. Then for any map
$\nu_p$ of \eqref{a.l} and any filtered complex $E_\idot$ in
$\Fun(\Lambda,R)$, the corresponding filtered map \eqref{nu.l.filt}
provided by Lemma~\ref{edge.filt.le} is a filtered quasiisomorphism.
\end{prop}

In order to prove this, we need a preliminary lemma. Let $k=\Z/p\Z$
be the prime field corresponding to $p$, fix an integer $n$, and
consider the group algebra $k[\Z/np\Z]$. We have a natural
augmentation map
$$
\aug:k[\Z/np\Z] \to k
$$
sending the generator $\sigma \in \Z/np\Z$ to $(-1)^{np}$, so that
$\aug(\sigma^\dg)=1$. Moreover, assume given two elements $a,b \in
k[\Z/np\Z]$ such that
\begin{equation}\label{a-b}
\begin{aligned}
&(1-\sigma^\dg)a=b(1-\sigma^\dg), \\
&(1+\sigma^\dg+\dots+(\sigma^\dg)^{np-1})b
= b(1+\sigma^\dg+\dots+(\sigma^\dg)^{np-1}).
\end{aligned}
\end{equation}
Then for any $k[\Z/pn\Z]$-module $E$, we have a functorial
$u$-equivariant map of twisted Tate homology complexes
\begin{equation}\label{f.ab}
f(a,b):\vC^\dg_\idot(\Z/np\Z,E) \to \vC^\dg_\idot(\Z/np\Z,E)
\end{equation}
given by
$$
f(a,b) = \begin{cases} a \text{ on } \C^\dg_{2i}(\Z/np,E),\\
b \text{ on } \C^\dg_{2i+1}(\Z/np,E)\end{cases}
$$
for any integer $i$, and any $u$-equivariant functorial map
\eqref{f.ab} arises in this way.

\begin{lemma}\label{a-b.le}
In the assumptions above, the following are equivalent:
\begin{enumerate}
\item $\aug(a) \neq 0$,
\item $\aug(b) \neq 0$,
\item $f(a,b)$ is a quasiisomorphism for any $E$.
\end{enumerate}
\end{lemma}

\proof{} Note that \eqref{a-b} implies that $(1-\sigma^\dg)(a-b)=0$,
and since $\Z/pn\Z$ has no homology with coefficients in its regular
representation, this implies that
$(a-b)=(1+\sigma^\dg+\dots+(\sigma^\dg)^{np-1})h$ for some $h \in
k[\Z/np\Z]$. Since
$\aug(1+\sigma^\dg+\dots+(\sigma^\dg)^{pn-1})=np=0$, we have
$\aug(a) = \aug(b)$, so that \thetag{i} and \thetag{ii} are
equivalent. Since any group acts trivially on its homology,
$\Z/pn\Z$ acts on the homology of the twisted Tate complex
$\vC^\dg(E)$ via the sign representation, so that the map $f(a,b)$
acts on this homology by $\aug(a)$ in even degrees and $\aug(b)$ in
odd degrees. Thus \thetag{i} and \thetag{ii} together imply
\thetag{iii}. Conversely, to deduce \thetag{i} and \thetag{ii} from
\thetag{iii}, take $E$ to be the trivial representation $k$.
\endproof

\proof[Proof of Proposition~\ref{edge.prop}.] By
\eqref{gr.stand.cp}, it suffices to prove the claim when $E_\idot$
is a single object $E$ placed in a single filtered degree $i$. Since
$R$ is $p$-torsion, the twisted Tate homology complexes
$\vC^\dg_\idot(\Z/m\Z,M)$ are acyclic for any $m$ not divisible by
$p$ and any $R[\Z/m\Z]$-module $M$, so that by \eqref{gr.stand.cp},
$\gr^m_Fcp_\idot(E)$ is acyclic unless $m-i$ divides $p$. Therefore
for any integer $m$, we have a quasiisomorphism
$$
\gr^m_Fcp_\idot(E_\idot)^{[p]} \cong \gr^{pn+i}_Fcp_\idot(E_\idot),
$$
where $n$ is the integer such that $pn \leq m-i < p(n+1)$. For any
$n \geq 1$, the map $\nu_p$ induces a functorial $u$-equivariant map
$$
\gr^{i-np}_F\nu_p:\vC^\dg_\idot(\Z/np\Z,i_p^*E([n])) \to
\vC^\dg_\idot(\Z/np\Z,E([np]))[-n(p-1)],
$$
and we have to prove that all these maps are
quasiisomorphisms. Since $i_p^*E([n]) \cong E([np])$, and the Tate
complex is $2$-periodic, thus $n(p-1)$-periodic, we have
\begin{equation}\label{gr.ab.n}
\gr^{i-np}_F\nu_p = f(a_n,b_n)
\end{equation}
for some $a_n,b_n \in \Z[\Z/np\Z]$ satisfying \eqref{a-b}. By
Lemma~\ref{a-b.le}, we then have to prove that $\aug(a_n)$ and
$\aug(b_n)$ are invertible for any $n \geq 1$.

We note that by construction, $a_n$ and $b_n$ are universal
constants that only depend on $p$, and do not depend on $R$ and
$E$. Thus we can take $R = k$, and take as $E$ the constant functor
$k \in \Fun(\Lambda,k)$ placed in filtered degree $0$. Then $i_p^*E$
is also constant, so that all the terms in the complexes
$CH_\idot(E)$, $CH'_\idot(E)$, $CH_\idot(i_p^*E)$,
$CH'_\idot(i_p^*E)$ are identified with $k$, and one immediately
checks that the differentials $d$, $d'$ are given by
$$
d = \begin{cases}
\id &\text{ on }CH_{2i-1}(i_p^*E),\\
0 &\text{ on }CH_{2i}(i_p^*E),
\end{cases}
\qquad
d' = \begin{cases}
0 &\text{ on }CH'_{2i-1}(i_p^*E),\\
\id &\text{ on }CH'_{2i}(i_p^*E),
\end{cases}
$$
as in Lemma~\ref{cp.const.le}. Therefore for any $i \geq 1$, we have
$$
\begin{aligned}
\aug(b_{2i-1}) &= \aug(b_{2i-1}) \circ d'  = d \circ \aug(a_{2i}),\\
\aug(a_{2i}) &= \aug(a_{2i}) \circ d  = d' \circ \aug(b_{2i+1}),
\end{aligned}
$$
and by Lemma~\ref{a-b.le} and induction on $i$, it suffices to prove
that $\aug(a_1) \in k$ is non-zero. By Lemma~\ref{cp.const.le},
this amounts to checking that the map
$$
\nu_p:cp_\idot(i_p^*E) \cong \bCP_\idot(i_p^*E) \to cp_\idot(E)
\cong \bCP_\idot(E)
$$
is a quasiisomorphism, and moreover, by the same
Lemma~\ref{cp.const.le}, we have quasiisomorphisms
$\bCP_\idot(i_p^*E) \cong CP_\idot(i_p^*E)$ and $\bCP_\idot(E) \cong
CP_\idot(E)$. Then we are done by Lemma~\ref{edge.le}.
\endproof

In particular, we can equip any complex $E_\idot$ in
$\Fun(\Lambda,R)$ with the stupid filtration. Then the map $\nu_p$
gives by completion the map \eqref{bcp.nu}, and
Proposition~\ref{edge.prop} shows that it is a quasiisomorphim
(provided $l=p$ is an odd prime that annihilates $R$). We also have
the following statement for the restricted complexes $CP^f(-)$,
$\bCP^f(-)$ of Definition~\ref{restr.def}.

\begin{corr}\label{edge.corr}
Under the assumptions of Proposition~\ref{edge.prop}, the map
$\nu_p$ extends to quasiisomorphisms
\begin{equation}\label{bcpf.nu}
CP^f_\idot(i_p^*E_\idot) \cong CP^f_\idot(E_\idot), \qquad
\bCP^f_\idot(i_p^*E_\idot) \cong \bCP^f_\idot(E_\idot).
\end{equation}
\end{corr}

\proof{} Place the complex $E_\idot$ in filtered degree $0$, and
equip the complexes $cp_\idot(E_\idot)$, $cp_\idot(i_p^*E_\idot)$
with the corresponding standard filtrations $F^\hdot$. Moreover,
assume given an additional filtration $W^\hdot$ on $E_\idot$, and
equip $i_p^*E_\idot$ with the induced filtration. Then since $i_p^*$
is an exact functor, we have $i_p^*\gr^\hdot_WE_\idot \cong
\gr^\hdot_Wi_p^*E_\idot$, and for any $n \geq 1$, the
quasiisomorphism
$$
F^{-n}\nu_p:F^{-n}cp_\idot(i_p^*E_\idot) \cong F^{-n}cp_\idot(E_\idot)
$$
of Proposition~\ref{edge.prop} is a filtered quasiisomorphism with
respect to filtrations induced by $W^\hdot$. Therefore if we denote
$$
cp^W_\idot(E_\idot) =
\lim_{\overset{n}{\to}}\widehat{F^{-n}cp_\idot(E_\idot)}, \qquad
cp^W_\idot(i_p^*E_\idot) =
\lim_{\overset{n}{\to}}\widehat{F^{-n}cp_\idot(i_p^*E_\idot)},
$$
where the completions in the right-hand side are taken with respect to
the filtration induced by $W^\hdot$, then $\nu_p$ gives a
quasiisomorphism
$$
\nu_p:cp^W_\idot(i_p^*E_\idot) \cong cp^W_\idot(E_\idot).
$$
It remains to notice that if $W^\hdot$ is the canonical resp. stupid
filtration on $E_\idot$, then by Lemma~\ref{tate.compl.le} and
\eqref{gr.stand.cp}, \eqref{gr.stand.CP}, the complex
$cp^W_\idot(E_\idot)$ coincides with $CP^f(E_\idot)$
resp. $\bCP^f(E_\idot)$.
\endproof

\subsection{Projections.}\label{proj.subs}

Consider now the second functor of \eqref{i.pi}, namely, the
projection $\pi_l:\Lambda_l \to \Lambda$. By definition, for any $E
\in \D(\Lambda_l,R)$, we have
$$
H_\idot(\Lambda_l,E) \cong H_\idot(\Lambda,L^\hdot\pi_{p!}E).
$$
However, unlike $i_l^*$, the functor $\pi_l^*$ does not induce an
isomorphism on homology and cohomology. In fact, the natural map
$$
\pi_l^*:H^\hdot(\Lambda_l,R) \to H^\hdot(\Lambda,R)
$$
sends the generator $u$ of the algebra $H^\hdot(\Lambda_l,R) \cong
R[u]$ to $lu$, where $u$ is the generator of $R[u] \cong
H^\hdot(\Lambda,R)$.

In particular, assume that $l=p$ is a prime such that $pR=0$. Then
$\pi_p^*u$ vanishes. Therefore by adjunction, for any $E \in
\D(\Lambda,R)$ such that $E \cong L^\hdot\pi_{p!}E'$ for some $E'
\in \D(\Lambda_p,R)$, the connecting differential in the long exact
sequence
\begin{equation}\label{gamma.seq}
\begin{CD}
H_\idot(\Lambda,E)[1] @>>> H_\idot(\Delta^o,j^*E) @>{\gamma}>>
H_\idot(\Lambda,R) @>>>
\end{CD}
\end{equation}
vanishes, so that the natural map $\gamma:H_\idot(\Delta^o,j^*E) \to
H_\idot(\Lambda,E)$ admits a splitting
\begin{equation}\label{delta.eq}
H_\idot(\Lambda,E) \to H_\idot(\Delta^o,j^*E).
\end{equation}
The goal of this subsection is to refine this observation to a
statement on the level of complexes, similar to
Proposition~\ref{edge.prop}.

For any complex $E_\idot$ in $\Fun(\Lambda_p,R)$,
$\pi_{p!}K_\idot(E_\idot)$ is a mixed complex in $\Fun(\Lambda,R)$.
Denote by
\begin{equation}\label{c.p.eq}
\pi_{p\flat}E_\idot = \per(\pi_{p!}\K_\idot(E_\idot))
\end{equation}
its polynomial periodic expansion. By virtue of the identification
\eqref{cc.pi.0}, \eqref{cp.bis.def} gives a canonical isomorphism
\begin{equation}\label{cc.pi}
cp_\idot(E_\idot) \cong cc_\idot(\pi_{p\flat}E_\idot).
\end{equation}
Moreover, denote
$$
cph_\idot(E_\idot) = CH_\idot(\pi_{p\flat}E_\idot) =
cc_\idot(\K_\idot(\pi_{p\flat}E_\idot)).
$$
Then by the projection formula, we have
\begin{equation}\label{K.pi}
\K_\idot(\pi_{p\flat}E_\idot) \cong \pi_{p\flat}E_\idot \otimes
\K_\idot \cong \pi_{p\flat}(E_\idot \otimes \pi_p^*\K_\idot),
\end{equation}
so that we have a natural isomorphism
\begin{equation}\label{cph.cp}
cph_\idot(E_\idot) \cong cp_\idot(E_\idot \otimes \pi_p^*\K_\idot),
\end{equation}
and the natural map \eqref{gamma} gives a natural
functorial map
\begin{equation}\label{eps.eq}
\gamma:cph_\idot(E_\idot)  \to cp_\idot(E_\idot),
\end{equation}
a chain-level lifting of the map $\gamma$ of \eqref{gamma.seq}.  Under
the identification \eqref{cph.cp}, the map $\gamma$ is induced by the
natural map of complexes
\begin{equation}\label{gamma.loc}
\wt{\gamma}:\pi_p^*\K_\idot \to k
\end{equation}
given by the pullback $\pi^*_p(\kappa_0)$ of the map $\kappa_0$ of
\eqref{K.la.eq} in homological degree $0$.

The functor $\pi_{p\flat}$ of \eqref{c.p.eq} is exact, so that any
filtration $F^\hdot$ on a complex $E_\idot$ in $\Fun(\Lambda_p,R)$
induces a filtration
\begin{equation}\label{c.p.filt}
F^i\pi_{p\flat}E_\idot = \pi_{p\flat}F^iE_\idot \subset
\pi_{p\flat}E_\idot, \qquad i \in \Z
\end{equation}
on $\pi_{p\flat}E_\idot$. If $F^\hdot$ was termwise-split, then this
induced filtration is also termwise-split, and if we equip
$cp_\idot(E_\idot)$ with the standard filtration of
Definition~\ref{stand.filt.def}, the isomorphism \eqref{cc.pi} is a
filtered isomorphism. Analogously, if we equip the product $E_\idot
\otimes \pi^*_p\K_\idot$ with the filtration induced by $F^\hdot$,
and consider the corresponding standard filtration on
$cp_\idot(E_\idot \otimes \pi^*_p\K_\idot)$, the identification
\eqref{cph.cp} becomes a filtered isomorphism, and the map
\eqref{eps.eq} is a filtered map.

\begin{lemma}\label{proj.le}
Assume that $pR=0$. Then for any complex $E_\idot$ in the category
$\Fun(\Lambda_p,R)$, we have a natural complex
$\wt{cp}_\idot(E_\idot)$ and a map
$$
\delta:\wt{cp}_\idot(E_\idot) \to cph_\idot(E_\idot)
$$
such that $\gamma \circ \delta:\wt{cp}_\idot(E_\idot) \to
cp_\idot(E_\idot)$ is a quasiisomorphism. Both
$\wt{cp}_\idot(E_\idot)$ and $\delta$ are functorial in
$E_\idot$. Moreover, if $E_\idot$ carries a termwise-split
filtration, then $\wt{cp}_\idot(E_\idot)$ carries a functorial
standard filtration such that $\delta$ is a filtered map, and
$\gamma \circ \delta$ is a filtered quasiisomorphism.
\end{lemma}

\proof{} Since $p$ annihilates $R$, we may replace the objects
$\K_0$, $\K_1$ in \eqref{K.eq} and \eqref{K.la.eq} with their
reductions modulo $p$, and take the tensor product in \eqref{K.grp}
over $k=\Z/p\Z$. Then \eqref{K.la.eq} provides an exact sequence
$$
\begin{CD}
0 @>>> k @>{\kappa_1}>> \K_1 @>>> \K_0 @>{\kappa_0}>>
k @>>> 0
\end{CD}
$$
in $\Fun(\Lambda,k)$ that represents by Yoneda the generator $u$ of
the cohomology algebra $H^\hdot(\Lambda,k) \cong k[u]$. Since
$\pi_p^*u=0$, the sequence
$$
\begin{CD}
0 @>>> k @>{\kappa_1}>> \pi_p^*\K_1 @>>> \pi_p^*\K_0 @>{\kappa_0}>>
k @>>> 0
\end{CD}
$$
represents the trivial class. Therefore by the standard criterion,
there exists an object $\K_{01} \in \Fun(\Lambda_p,k)$ with a
three-step filtration $w^i\K_{01}$, $i=0,1,2$, such that $w^2\K_{01}
\cong k$, $w^1\K_{01} \cong \pi_p^*\K_1$, $\K_{01}/w^2\K_{01} \cong
\pi_p^*\K_0$, $\K_{01}/w^1\K_{01} \cong k$. Equivalently, denote by
$\wt{\K}_\idot$ the complex with terms $\wt{\K}_1 = \pi_p^*\K_1$,
$\wt{\K}_0 = \K_{01}$, and the differential given by the embedding
$\pi^*_p\K_1 \cong w^1\K_{01} \hookrightarrow \K_{01}$. Then the
projection $\K_{01} \to \K_0$ induces a map of complexes
\begin{equation}\label{delta.loc}
\wt{\delta}:\wt{\K}_\idot \to \pi_p^*\K_\idot,
\end{equation}
and its composition $\wt{\gamma} \circ \wt{\delta}:\wt{\K}_\idot \to
k$ with the map \eqref{gamma.loc} is a surjective map of complexes
with contractible kernel. To prove the lemma, it remains to denote
\begin{equation}\label{wt.E}
\wt{E}_\idot = E_\idot \otimes_k \wt{\K}_\idot
\end{equation}
and let $\wt{cp}_\idot(E_\idot) = cp_\idot(\wt{E}_\idot)$, $\delta =
cp_\idot(\id \otimes \wt{\delta})$. Indeed, $\wt{E}_\idot \to
E_\idot = \id \otimes (\wt{\gamma} \circ \wt{\delta})$ is also a
surjective map of complexes with contractible kernel, hence a strong
quasiisomorphism in the sense of Definition~\ref{strong.qis}, and
therefore the map $\delta \circ \gamma = cp_\idot(\id \otimes
(\wt{\delta} \circ \wt{\gamma}))$ is a quasiisomorphism. The same is
true in presence of filtrations.
\endproof

We note that since the map $\delta$ provided by Lemma~\ref{proj.le}
is compatible with filtrations, it automatically has a completed
counterpart. Namely, let $\bCPH^{[p]}_\idot(E_\idot)$ be the
completion of the complex $cph_\idot(E_\idot)$ with respect to the
standard filtration corresponding to the $p$-th rescaling of the
stupid filtration on $E_\idot$. Then the map $\gamma$ of
\eqref{eps.eq} gives by completion a natural map
$$
\gamma:\bCPH^{[p]}_\idot(E_\idot) \to \bCP^{[p]}_\idot(E_\idot),
$$
where $\bCP^{[p]}_\idot(E_\idot)$ is as in Definition~\ref{cpv.def}, and
map $\delta$ of Lemma~\ref{proj.le} extends to a map of complete
complexes such that $\delta \circ \gamma$ is a
quasiisomorphism. Analogously, denote $CPH_\idot(E_\idot) =
CP_\idot(E_\idot \otimes \pi_p^*\K_\idot)$. Then we have natural
maps
$$
\begin{CD}
CP_\idot(\wt{E}_\idot) @>{\delta}>> CPH_\idot(E_\idot) @>{\gamma}>>
CP_\idot(E_\idot)
\end{CD}
$$
whose completion is a quasiisomorphism. Moreover,
Lemma~\ref{proj.le} has an obvious counterpart for the restricted
periodic complexes $CP^f_\idot(E_\idot)$,
$\bCP^f_\idot(E_\idot)$. Namely, denote by
\begin{equation}\label{vpi}
\vpi_{p\flat}E_\idot = \Per(\pi_{p!}\K_\idot(E_\idot)), \qquad
\vbpi_{p\flat}E_\idot = \cPer(\pi_{p!}\K_\idot(E_\idot))
\end{equation}
the periodic and co-periodic expansions of the mixed complex
$\pi_{p!}(\K_\idot(E_\idot))$, and let
$$
\begin{aligned}
CPH^f_\idot(E_\idot) &= CH_\idot(\vpi_{p\flat}E_\idot) \cong
CP^f_\idot(E_\idot \otimes \pi_p^*\K_\idot),\\
\bCPH^f_\idot(E_\idot) &= CH_\idot(\vbpi_{p\flat}E_\idot) \cong
\bCP^f_\idot(E_\idot \otimes \pi_p^*\K_\idot).
\end{aligned}
$$
Then we have natural maps
$$
\begin{aligned}
\gamma:&CPH^f_\idot(E_\idot) \to CP^f_\idot(E_\idot) \cong
cc_\idot(\vpi_{p\flat}E_\idot),\\
\gamma:&\bCPH^f_\idot(E_\idot) \to \bCP^f_\idot(E_\idot) \cong
cc_\idot(\vbpi_{p\flat}E_\idot),
\end{aligned}
$$
and the same argument as in Lemma~\ref{proj.le} provides maps
$\delta$ such that $\gamma \circ \delta$ is a
quasiisomorphism. 

\section{Computational tools.}\label{conj.sec}

\subsection{Conjugate spectral sequence.}\label{conj.subs}

As in Subsection~\ref{edge.filt.subs} and
Subsection~\ref{proj.subs}, assume that the base ring $R$ is
annihilated by an odd prime $p$, and let $k=\Z/p\Z$. Then for any
complex $E_\idot$ in $\Fun(\Lambda,R)$, we have the filtered
quasiisomorphism $\nu_p$ of Proposition~\ref{edge.prop} and its
completed version \eqref{bcp.nu}. To study the complex
$\bCP^{[p]}(i_p^*E_\idot)$, we will use the truncation functors
$\tau^\hdot$, $\beta^\hdot$ of Subsection~\ref{t.subs}.

\begin{defn}\label{conj.def}
For any termwise-split filtered complex $E_\idot$ in
$\Fun(\Lambda_p,R)$, the filtrations $V^\hdot$ and $\bV^\hdot$ on
$cp_\idot(E_\idot)$ are given by
\begin{equation}\label{V.bV}
V^n = \tau^{2n-1}cp_\idot(E_\idot), \quad \bV^n =
\beta^{2n-2}cp_\idot(E_\idot), \qquad n \in \Z,
\end{equation}
where we equip $cp_\idot(E_\idot)$ with the standard filtration
induced by the filtration on $E_\idot$. The filtration $V^\hdot$ is
called the {\em conjugate filtration}.
\end{defn}

Note that both filtrations are $2$-periodic -- for every integer
$n$, the periodicity endomorphism $u:cp_\idot(E_\idot) \to
cp_\idot(E_\idot)[2]$ induces isomorphisms
\begin{equation}\label{V.per}
u:V^{n+1}cp_\idot(E_\idot) \cong V^ncp_\idot(E_\idot)[2], \
u:\bV^{n+1}cp_\idot(E_\idot) \cong \bV^ncp_\idot(E_\idot)[2],
\end{equation}
where $[2]$ stands for the cohomological shift. By
Lemma~\ref{tau.le}, we have $V^{n+1}cp_\idot(E_\idot) \subset
\bV^ncp_\idot(E_\idot) \subset V^ncp_\idot(E_\idot)$, and the first
embedding is a quasiisomorphism. The conjugate filtration $V^\hdot$
is by definition the rescaling by $2$ and shift by $1$ of the
filtration $\tau^\hdot$. In particular, $V^\hdot$ and $\tau^\hdot$
are commensurable, so that by Lemma~\ref{tau.le}, the conjugate
filtration is commensurable with the standard filtration. Thus if we
equip a complex $E_\idot$ in $\Fun(\Lambda_p,R)$ with the $p$-th
rescaling of the stupid filtration, the completion of
$cp_\idot(E_\idot)$ with respect to the conjugate filtration
coincides with the complex $\bCP^{[p]}(\wt{E}_\idot)$ of
Definition~\ref{cpv.def}. By virtue of the quasiisomorphism
\eqref{bcp.nu}, we then have a a spectral sequence
\begin{equation}\label{conj.sp}
H_\idot(\gr^\hdot_V(cp_\idot(i_p^*E_\idot^{[p]}))) \Rightarrow
\bHP^{[p]}_\idot(i_p^*E_\idot^{[p]}) = \bHP_\idot(E_\idot)
\end{equation}
for any complex $E_\idot$ in $\Fun(\Lambda,R)$. As it turns out,
this spectral sequence is quite useful, because under some
assumptions on $E_\idot$, one can find a rather effective
description of its initial term.

\medskip

To do this, we need to recall some material from \cite{ka1}. The
cohomology $H^\hdot(\Z/p\Z,k)$ with coefficients in $k=\Z/p\Z$ is
the graded-commutative algebra given by
\begin{equation}\label{z.pz.coh}
H^\hdot(\Z/p\Z,k) = k[u]\langle\eps\rangle,
\end{equation}
where $u$ is a generator of degree $2$ and $\eps$ is a generator of
degree $1$. The generator $\eps$ gives an extension $\wt{k}$ of the
trivial $k[\Z/p\Z]$-module $k$ by itself, so that for any
$R[\Z/p\Z]$-module $E$, we have a functorial short exact sequence
\begin{equation}\label{eps.seq}
\begin{CD}
0 @>>> E @>>> \wt{E} @>>> E @>>> 0,
\end{CD}
\end{equation}
where we let $\wt{E} = E \otimes_k \wt{k}$. Taking the Tate homology
complex $\vC_\idot(-)$ of \eqref{tate.eq}, we obtain a short exact
sequence of complexes
$$
\begin{CD}
0 @>>> \vC_\idot(\Z/p\Z,E) @>>> \vC_\idot(\Z/p\Z,\wt{E}) @>>>
\vC_\idot(\Z/p\Z,E) @>>> 0.
\end{CD}
$$
This defines a distinguished triangle in $\DF(R)$, so that we have a
connecting differential $\vC_\idot(\Z/p\Z,E) \to
\vC_\idot(\Z/p\Z,E)[1]$ and the corresponding maps
\begin{equation}\label{eps.i.0}
\eps_i:\vH_i(\Z/p\Z,E) \to \vH_{i-1}(\Z/p\Z,E)
\end{equation}
for any integer $i$. In effect, the cohomology algebra
$H^\hdot(\Z/p\Z,k)$ acts on $\vH_\idot(\Z/p\Z,E)$, and the maps
$\eps_i$ give the action of the generator $\eps$. The generator $u$
acts by the $2$-shift in the periodic complex $\vC_\idot(\Z/p\Z,E)$,
so that $\eps_i=\eps_{i+2}$ for any $i$. Thus effectively, we only
have two maps $\eps_{odd}$ and $\eps_{even}$, depending on the
parity of $i$.  An $R[\Z/p\Z]$-module $E$ is called {\em tight} if
the map $\eps_{odd}$ is an isomorphism. Since $\eps^2=0$, we have
$\eps_{odd} \circ \eps_{even}=0$; thus for a tight module, we
automatically have $\eps_{even}=0$.

\begin{defn}\label{tight.def}
A complex $E_\idot$ of $R[\Z/p\Z]$-modules is {\em tight} if $E_i$
is a tight $R[\Z/p\Z]$-module for any $i$, and $\I(E_i)=0$ unless
$i$ divides $p$.
\end{defn}

Now assume given a filtered complex $E_\idot$ of
$R[\Z/p\Z]$-modules, and consider the Tate complex
$\wt{C}_\idot(\Z/p\Z,E_\idot)$ of \eqref{tate.bis}. Equip it with
the filtration $F^\hdot$ induced by the filtration on $E_\idot$, and
for any integer $i$, denote
$$
\wt{H}_i(\Z/p\Z,E_\idot) \cong H_i(\wt{C}_\idot(\Z/p\Z,E_\idot)),
$$
where $H_i(-)$ is as in \eqref{h.i}. Then again, for any
integer $i$, \eqref{eps.seq} induces a natural map
\begin{equation}\label{eps.i}
\eps_i:\wt{H}_i(\Z/p\Z,E_\idot) \to \wt{H}_{i-1}(\Z/p\Z,E_\idot)[1],
\end{equation}
a generalization of \eqref{eps.i.0}. As before, $\eps_i$ only
depends on the parity of $i$, so that effectively, we only have two
maps $\eps_{odd}$ and $\eps_{even}$.

\begin{lemma}\label{tight.le}
Assume given a complex $E_\idot$ of $R[\Z/p\Z]$-modules tight in the
sense of Definition~\ref{tight.def}, and equip it with the $p$-th
rescaling $F^\hdot_{[p]}$ of the stupid filtration $F^\hdot$. Then
$\eps_{even}=0$, and $\eps_{odd}$ is an isomorphism.
\end{lemma}

\proof{} Immediately follows from Lemma~\ref{tau.le} (note that
since $p$ is odd by assumption, $(p-1)i$ is even for any $i$).
\endproof

By virtue of Lemma~\ref{tight.le}, for any tight complex $E_\idot$
of $R[\Z/p\Z]$-modules, the complexes $\wt{H}_i(\Z/p\Z,E_\idot)[-i]$
for all integers $i$ are canonically identified. We denote this
complex by $\I(E_\idot)$. By Lemma~\ref{tau.le}, the natural
filtration on $\I(E_\idot)$ coincides with the stupid filtration,
and we have
\begin{equation}\label{i.e.n}
\I(E_\idot)_n = \I(E_{pn})
\end{equation}
for any integer $n$.

Assume now given a complex $E_\idot$ in $\Fun(\Lambda_p,R)$. Then
for any object $[n] \in \Lambda_p$, $E_\idot([n])$ is a naturally a
complex of $R[\Z/p\Z]$-modules via the embedding $\Z/p\Z \subset
\Z/pn\Z = \Aut([n])$.

\begin{defn}\label{tight.la}
A complex $E_\idot$ in $\Fun(\Lambda_p,R)$ is {\em tight} if
$E_\idot([n])$ is a tight complex of $R[\Z/p\Z]$-modules for any
$[n] \in \Lambda_p$. A complex $E_\idot$ in $\Fun(\Lambda,R)$ is
{\em $p$-adapted} if $i_p^*E_\idot$ is tight.
\end{defn}

\begin{defn}
For any tight complex $E_\idot$ in $\Fun(\Lambda_p,R)$, the complex
$\I(E_\idot)$ in $\Fun(\Lambda,R)$ is given by
$$
\I(E_\idot)([n]) = \I(E_\idot([n]))
$$
for every $[n] \in \Lambda$.
\end{defn}

Note that for every $n \geq 1$, we have a base change isomorphism
\begin{equation}\label{bc.n}
\pi_{p\flat}E_\idot([n]) \cong \wC_\idot(\Z/p\Z,\Z/pn\Z,E_\idot([n])),
\end{equation}
where the right-hand side is the extended version \eqref{c.cpr} of
the Tate homology complex $\wC_\idot(\Z/p\Z,-)$. Since for any $n$,
the complex $\wC_\idot(\Z/p\Z,\Z/pn\Z,-)$ is canonically
chain-homotopy equivalent to the usual Tate homology complex
$\wC_\idot(\Z/p\Z,-)$, the isomorphism \eqref{bc.n} provides a
natural identification
$$
\I(E_\idot) \cong H_i(\pi_{p\flat}E_\idot^{[p]})[-i]
$$
for any integer $i$, where $E_\idot^{[p]}$ is $E_\idot$ equipped
with the $p$-th rescaling of the stupid filtration, and $H_i(-)$ is
the truncation functor \eqref{h.i} in the category of termwise-split
filtered complexes in $\Fun(\Lambda,R)$. By abuse of notation, for
any complex $E_\idot$ in $\Fun(\Lambda,R)$, we will denote
$$
\I(E_\idot) = \I(i_p^*E_\idot).
$$
With these definitions, the main result concerning the spectral
sequence \eqref{conj.sp} is the following.

\begin{prop}\label{conj.prop}
Assume that $pR=0$ for an odd prime $p$. Then for any tight complex
$E_\idot$ in $\Fun(\Lambda_p,R)$, we have a natural quasiisomorphism
\begin{equation}\label{gr.v.pr}
\gr^0_V(cp_\idot(E_\idot^{[p]}))) \cong CH_\idot(\I(E_\idot))
\end{equation}
functorial in $E_\idot$, so that for any $p$-adapted complex
$E_\idot$ in $\Fun(\Lambda,R)$, \eqref{conj.sp} induces a functorial
spectral sequence
\begin{equation}\label{conj.sp.1}
HH_\idot(\I(E_\idot))((u^{-1})) \Rightarrow \bHP_\idot(E_\idot).
\end{equation}
\end{prop}

Note that by virtue of the periodicity isomorphism \eqref{V.per},
the identification \eqref{gr.v.pr} describes the whole associated
graded quotient $\gr^\hdot_Vcp_\idot(E_\idot)$. In
\eqref{conj.sp.1}, we use the same shorthand notation as in
Definition~\ref{per}. We will call the spectral sequence
\eqref{conj.sp.1} the {\em conjugate spectral sequence} for the
tight complex $E_\idot$.

\subsection{Localization.}\label{loc.subs}

To prove Proposition~\ref{conj.prop}, we first need to localize the
conjugate filtration onto the category $\Lambda$. For any complex
$E_\idot$ in $\Fun(\Lambda,R)$ equipped with a termwise-split
filtration $F^\hdot$, we denote by $V^\hdot E_\idot$ and $\bV^\hdot
E_\idot$ the filtrations given by
\begin{equation}\label{V.bV.loc}
V^n E_\idot = \tau^{2n}E_\idot, \qquad \bV^n E_\idot =
\beta^{2n-1}E_\idot, \qquad n \in \Z,
\end{equation}
where $\tau^\hdot$ and $\beta^\hdot$ are as in
\eqref{tau.beta.F}. By Lemma~\ref{tau.le}, for any $E_\idot$, we
have $V^{n+1}E_\idot \subset \bV^nE_\idot \subset V^nE_\idot$, and
the first embedding is an isomorphism.

In particular, assume given a complex $E_\idot$ in the category
$\Fun(\Lambda_p,R)$, equip it with the $p$-th rescaling of the
stupid filtration, and to simplify notation, denote by
$E_\idot^\flat = \pi_{p\flat}E^{[p]}_\idot$ the complex
\eqref{c.p.eq} with the filtration \eqref{c.p.filt}. Then we have
filtrations $V^\hdot$, $\bV^\hdot$ on $E^\flat_\idot$, and both
filtrations are periodic in the same sense as
\eqref{V.per}. Analogously, let $\wt{E}_\idot^\flat =
\pi_{p\flat}\wt{E}_\idot$, where $\wt{E}_\idot$ is the complex
\eqref{wt.E} with the filtration induced from $E_\idot$. Then
$\wt{E}^\flat_\idot$ also carries periodic filtrations $V^\hdot$,
$\bV^\hdot$. By virtue of the identification \eqref{K.pi}, the map
\eqref{delta.loc} induces a map
\begin{equation}\label{delta.flat}
\delta^\flat:\wt{E}^\flat_\idot \to E^\flat_\idot \otimes \K_\idot.
\end{equation}

\begin{lemma}\label{I.le}
Assume that a complex $E_\idot$ in $\Fun(\Lambda_p,R)$ is tight in
the sense of Definition~\ref{tight.la}. Then the following is true.
\begin{enumerate}
\item For any integer $n$, the map $\delta^\flat$ of \eqref{delta.flat}
  sends $V^n\wt{E}^\flat_\idot \subset \wt{E}^\flat_\idot$ into
  $\bV^nE^\flat_\idot \otimes \K_\idot \subset E^\flat_\idot \otimes
  \K_\idot$.
\item Moreover, denote $\I_\idot = \gr^{-1}_\beta E^\flat_\idot$, and
  let $\rho:\bV^0E^\flat_\idot = \beta^{-1}E^\flat_\idot \to \I_\idot$
  be the natural projection. Then the composition map
$$
\begin{CD}
\gr^0_V\wt{E}^\flat_\idot @>{\delta^\flat}>> \gr^0_{\bV}E^\flat_\idot
\otimes \K_\idot @>{\rho \otimes \id}>> \I_\idot \otimes \K_\idot
\end{CD}
$$
is a quasiisomorphism.
\end{enumerate}
\end{lemma}

\proof{} Since by definition, the filtration on $E_\idot$ is
termwise-split, both claims commute with passing to the associated
graded quotients. Therefore we may assume right away that $E_\idot$
is concentrated in a single filtered degree, say $0$. Moreover, by
periodicity, it suffices to prove \thetag{i} for $n=0$. We have
$V^0=\tau^0$, $W^0=\beta^{-1}$, and
$$
\tau^0(E^\flat_\idot \otimes \K_\idot)/((\beta^{-1}E^\flat_\idot
\otimes \K_\idot)\cap \tau^0(E^\flat_\idot \otimes \K_\idot)) \cong
H_{-1}(E^\flat_\idot) \otimes H_1(\K_\idot).
$$
Since $H_1(\K_\idot) \cong \Z$, we have $H_{-1}(E^\flat_\idot)
\otimes H_1(\K_\idot) \cong H_{-1}(E^\flat_\idot)$, and the map
$$
\begin{CD}
\tau^0\wt{E}^\flat_\idot @>{\delta^\flat}>> \tau^0(E^\flat_\idot
\otimes \K_\idot) @>>> H_{-1}(E^\flat_\idot) \otimes H_1(\K_\idot)
\cong H_{-1}(E^\flat_\idot)
\end{CD}
$$
factors through a map
\begin{equation}\label{eps.loc}
H_0(\wt{E}^\flat_\idot) =
\tau^0\wt{E}^\flat_\idot/\beta^0\wt{E}^\flat_\idot \to
H_{-1}(E^\flat_\idot).
\end{equation}
We have to show that this map is equal to $0$. This claim is local
with respect to $\Lambda$, so that it suffices to prove it after
evaluation at an arbitrary object $[n] \in \Lambda$. We have
$$
\wt{E}^\flat_\idot([n]) \cong \wt{C}_\idot(\Z/p\Z,M_\idot \otimes_k
K_\idot), \qquad E^\flat_\idot([n]) \cong
\wt{C}_\idot(\Z/p\Z,M_\idot),
$$
where we denote $M_\idot = E_\idot([n])$, $K_\idot =
\wt{\K}_\idot([n])$, and \eqref{eps.loc} evaluates to a map
\begin{equation}\label{eps.loc.bis}
\vH_0(\Z/p\Z,M \otimes K_\idot) \to \vH_{-1}(\Z/p\Z,M),
\end{equation}
where $M = H_0(M_\idot)$ stand for the only non-trivial homology
group of the complex $M_\idot$.

However, since the $\Z/p\Z$-action on $K_1=\pi^*_p\K_1([n])$ is
trivial, the embedding $\kappa_1:k \to K_1$ splits as a map of
$k[\Z/p\Z]$-modules. Choosing such a splitting gives a
quasiisomorphism between $K_\idot$ and the complex
$\overline{K}_\idot$ with terms $\overline{K}_1 = k$,
$\overline{K}_0 = \wt{k}$, and the differential given by the
embedding $k \to \wt{k}$ of \eqref{eps.seq}. Then in defining the
map \eqref{eps.loc.bis}, we may replace $K_\idot$ with
$\overline{K}_\idot$, and the map becomes precisely the map
$\eps_{even}$ for the $k[\Z/p\Z]$-module $M$. Since $M$ is tight by
assumption, $\eps_{even}=0$. This proves \thetag{i}.

The argument for \thetag{ii} is similar --- the non-trivial part is
to check that the natural map
$$
H_1(\wt{E}^\flat_\idot) \to H_0(\I_\idot) \otimes H_1(\K_\idot)
\cong H_0(E^\flat_\idot)
$$
induced by $\delta^\flat$ is an isomorphism, this is a local fact,
and after evaluation at $[n] \in \Lambda$ and choosing a
quasiisomorphism $K_\idot \cong \overline{K}_\idot$, the map becomes
the map $\eps_{odd}$ for the tight $k[\Z/p\Z]$-module $M$. We leave
the details to the reader.
\endproof

\begin{remark}
Lemma~\ref{I.le}~\thetag{ii} is a strengthening of \cite[Lemma
  3.6]{ka1}. The proof is also essentially the same; the only
difference is that instead of the complex $\wt{\K}_\idot$,
\cite{ka1} uses an arbitrary resolution of the constant functor $k$
by objects in $\Fun(\Lambda_p,R)$ acyclic for the functor
$\pi_{p!}$.
\end{remark}

Now note that by Lemma~\ref{tau.func.le}, for any integer $n$, the
inclusion $V^n\wt{E}^\flat_\idot \subset \wt{E}^\flat_\idot$ induces
a surjective map
$$
\xi:cc_\idot(V^n\wt{E}^\flat_\idot) \to
V^ncc_\idot(\wt{E}^\flat_\idot) \cong V^ncp_\idot(\wt{E}_\idot),
$$
where $V^n$ in the right-hand side is the conjugate filtration of
Definition~\ref{conj.def}, and the shift by $1$ between \eqref{V.bV}
and \eqref{V.bV.loc} compensates for the shift by one in
Definition~\ref{stand.filt.def}. Composing $\xi$ with the map
$\alpha$ of \eqref{aug.cc}, we obtain a natural map
\begin{equation}\label{alp.V}
CC_\idot(V^n\wt{E}^\flat_\idot) \to V^ncp_\idot(\wt{E}_\idot).
\end{equation}
Since $CC_\idot(-)$ is an exact functor, it induces a map
\begin{equation}\label{CC.cp}
CC_\idot(V^{[n,m]}\wt{E}^\flat_\idot) \to V^{[n,m]}(\wt{E}_\idot)
\end{equation}
for any integer $m \geq n$.

\begin{lemma}\label{V.CC.le}
For any tight filtered complex $E_\idot$ in $\Fun(\Lambda_p,R)$, and
any integer $n \leq m$, the map \eqref{CC.cp} is a quasiisomorphism.
\end{lemma}

\proof{} By induction, it suffices to consider the case $m=n+1$, and
then by periodicity, it suffices to prove that the map
$$
CC_\idot(\gr^0_V\wt{E}^\flat_\idot)
\to \gr^0_Vcp_\idot(\wt{E}_\idot)
$$
induced by \eqref{alp.V} is a quasiisomorphism. Lemma~\ref{proj.le}
provides maps
$$
\begin{CD}
cp_\idot(\wt{E}_\idot) = \wt{cp}_\idot(E_\idot) @>{\delta}>>
cph_\idot(E^\hdot) @>{\gamma}>> cp_\idot(E_\idot)
\end{CD}
$$
whose composition is a filtered quasiisomorphism. By definition, we
have $cph_\idot(E_\idot) \cong CH_\idot(E^\flat_\idot)$, and since
$CH_\idot(-)$ is an exact functor, the natural map
$\iota:CH_\idot(\bV^0E^\flat_\idot) \to CH_\idot(E^\flat_\idot) =
cph_\idot(E_\idot)$ is injective. By Lemma~\ref{I.le}~\thetag{i}, we
then have a commutative square of complexes
$$
\begin{CD}
cc_\idot(V^0\wt{E}^\flat_\idot) @>{cc_\idot(\delta^\flat})>>
CH_\idot(\bV^0E^\flat_\idot)\\
@V{\xi}VV @VV{\iota}V\\
V^0cc_\idot(\wt{E}^\flat_\idot) @>{\delta}>> cph_\idot(E_\idot).
\end{CD}
$$
Since $\xi$ is surjective and $\iota$ is injective, there exists a
map $\bdelta:V^0cc_\idot(\wt{E}^\flat_\idot) \to
CH_\idot(\bV^0E^\flat_\idot)$ such that $cc_\idot(\delta^\flat) =
\bdelta \circ \xi$ and $\delta = \iota \circ \bdelta$. Now note that
the map $\wt{\gamma}$ of \eqref{gamma.loc} induces a map
$\gamma^\flat:E^\flat_\idot \otimes \K_\idot \to E^\flat$ such that
$\gamma^\flat \circ \delta^\flat$ is a filtered quasiisomorphism, so
that altogether, we have a commutative diagram
$$
\begin{CD}
CC_\idot(\gr^0_V\wt{E}^\flat_\idot) @>{CC_\idot(\delta^\flat)}>>
  CC_\idot(\gr^0_\bV E^\flat_\idot \otimes \K_\idot)
  @>{CC_\idot(\gamma^\flat})>> CC_\idot(\gr^0_\bV E^\flat_\idot)\\
@VVV @V{\alpha}VV @VVV\\
\gr^0_Vcp_\idot(\wt{E}_\idot) @>{\bdelta}>>
CH_\idot(\gr^0_\bV E^\flat_\idot) @>{\gamma \circ \iota}>>
\gr^0_\bV cp_\idot(E_\idot).
\end{CD}
$$
The composition $\gamma \circ \iota \circ \bdelta = \gamma \circ
\delta$ is a quasiisomorphism. Since $CC_\idot(-)$ is an exact
functor, $CC_\idot(\gamma^\flat) \circ CC_\idot(\delta^\flat) =
CC_\idot(\gamma^\flat \circ \delta^\flat)$ is also a quasiisomorphism,
and the map $\alpha$ is a quasiisomorphism by
Lemma~\ref{c.cc.le}. Therefore the remaining vertical arrows are
quasiisomorphisms, and we are done.
\endproof

\proof[Proof of Proposition~\ref{conj.prop}.]
Assume given a tight complex $E_\idot$ in $\Fun(\Lambda_p,R)$, equip
it with the $p$-th rescaling of the stupid filtration, and keep the
notation introduced earlier in this Subsection. Then
Lemma~\ref{I.le} and Lemma~\ref{c.cc.le} provide a quasiisomorphism
$$
CC_\idot(\gr^0_V\wt{E}^\flat_\idot) \cong CC_\idot(\I_\idot \otimes
\K_\idot) \cong CH_\idot(\I_\idot),
$$
and Lemma~\ref{V.CC.le} further identifies this with
$\gr^0_Vcp_\idot(\wt{E}_\idot)$. It remains to notice that
since the embedding $\tau^0E^\flat_\idot
\subset \beta^{-1}E^\flat_\idot$ is a quasiisomorphism, the
embedding $\I(E_\idot) \subset \I_\idot$ is a quasiisomorphism, so
that $CH_\idot(\I_\idot)$ is quasiisomorphic to
$CH_\idot(\I(E_\idot))$, and moreover, since $E_\idot$ is
quasiisomorphic to $\wt{E}_\idot$ as a filtered complex,
$\gr^0_Vcp_\idot(\wt{E}_\idot)$ is quasiisomorphic to
$\gr^0_Vcp_\idot(E_\idot)$.
\endproof

\subsection{Comparison maps.}\label{comp.subs}

We now turn to the other functorial complexes in \eqref{5.dia} and
the comparison maps $r$, $R$, $l$, $L$ between them. Keep the
assumption $pR=0$, $p$ an odd prime.  Recall that for any complex
$E_\idot$ in $\Fun(\Lambda_p,R)$, we have the complex
$\pi_{p\flat}E_\idot$ of \eqref{c.p.eq} and its completed versions
\eqref{vpi}.

\begin{defn}
A complex $E_\idot$ in $\Fun(\Lambda_p,R)$ is {\em locally bounded
  from below} resp.\ {\em locally strongly bounded from above} if
for any $[n] \in \Lambda_p$, the complex $E_\idot([n])$ of
$R[\Z/pn\Z]$-modules is bounded from below resp.\ strongly bounded
from above in the sense of Definition~\ref{strong.bnd}.
\end{defn}

\begin{lemma}\label{loc.le}
Assume given a complex $E_\idot$ in $\Fun(\Lambda_p,R)$. Then the
natural map $\alpha:CC_\idot(\pi_{p\flat}E_\idot) \to
cp_\idot(E_\idot)$ induced by the map \eqref{aug.cc} is a
quasiisomorphism, and so are the natural maps
\begin{equation}\label{alpha.f}
\alpha:CC_\idot(\vpi_{p\flat}E_\idot) \to CP^f_\idot(E_\idot), \quad
\alpha:CC_\idot(\vbpi_{p\flat}E_\idot) \to \bCP^f_\idot(E_\idot)
\end{equation}
\end{lemma}

\proof{} As in the proof of Lemma~\ref{V.CC.le}, Lemma~\ref{proj.le}
shows that to prove the first claim, it suffices to prove that the
natural map
$$
CC_\idot(\pi_{p\flat}E_\idot \otimes \K_\idot) \to
cph_\idot(E_\idot) \cong CH_\idot(\pi_{p\flat}E_\idot)
$$
is a quasiisomorphism. This immediately follows from
Lemma~\ref{c.cc.le}. For the maps \eqref{alpha.f}, use the
counterpart of Lemma~\ref{proj.le} for the restricted complexes, and
again apply Lemma~\ref{c.cc.le}.
\endproof

\begin{corr}\label{bnd.la.corr}
Assume given a complex $E_\idot$ in $\Fun(\Lambda_p,R)$. If
$E_\idot$ is locally bounded from below, then the natural map $l$ of
\eqref{5.dia} is a quasiisomorphism.  If $E_\idot$ is locally
strongly bounded from above, then the map $r$ of \eqref{5.dia} is a
quasiisomorphism.
\end{corr}

\proof{} By Lemma~\ref{loc.le}, the maps $l$ and $r$ are obtained by
applying the exact functor $CC_\idot(-)$ to their local counterparts
$$
l_\flat:\pi_{p\flat}E_\idot \to \vpi_{p\flat}E_\idot, \qquad
r_\flat:\pi_{p\flat}E_\idot \to \vbpi_{p\flat}E_\idot.
$$
Therefore both claims immediately follow from
Lemma~\ref{tate.bnd.le}.
\endproof 

The situation for $L$ and $R$, the other two comparison maps in
\eqref{5.dia}, is more difficult. In fact, we can only prove any
useful statements in the special case of cyclic complexes coming
from DG algebras; this we do later in Section~\ref{alg.sec}. For
now, we prepare the ground by axiomatizing the situation and proving
some easy auxiliary results.

First of all, what we are really interested in are complexes
$E_\idot$ in the category $\Fun(\Lambda,R)$, but we study them by
applying the quasiisomorphism of Proposition~\ref{edge.prop} and its
completed versions \eqref{bcp.nu}, \eqref{bcpf.nu}. Thus for a
complex $E_\idot$ in $\Fun(\Lambda_p,R)$, we need to
consider the completion $\bCP^{[p]}(E_\idot)$ of the complex
$cp_\idot(E_\idot)$ introduced in Definition~\ref{cpv.def}.

\begin{lemma}\label{L.v.le}
For any complex $E_\idot$ in $\Fun(\Lambda_p,R)$, there exists a
functorial map $L^{[p]}:\bCP^f(E_\idot) \to
\bCP^{[p]}_\idot(E_\idot)$ such that $L^{[p]} \circ
l:cp_\idot(E_\idot) \to \bCP^{[p]}_\idot(E_\idot)$ is the completion
map.
\end{lemma}

\proof{} Recall that by definition, we have
$$
cp_\idot(E_\idot) = cc_\idot(\per(\K_\idot(E_\idot)) = \bigoplus_{i
  \geq 0}cc_i(\per(\K_\idot(E_{\idot-i})))),
$$
where the right-hand side is the decomposition \eqref{phi.c}. For
any $i \geq 0$, the restriction of the standard filtration on
$cp_\idot(E_\idot^{[p]})$ to the summand $cc_i(-)$ in this decomposition
is a shift of the filtration induced by
the $p$-th rescaling of the stupid filtration on
$E_\idot$. Therefore by Lemma~\ref{tate.compl.le}, the completion of
$cc_i(-)$ with respect to this restricted filtration is precisely
$cc_i(\cPer(\K_\idot(E_{\idot-i})))$. Summing up over all $i$, we
obtain a map
$$
L^{[p]}:\bigoplus_{i \geq 0}cc_i(\cPer(\K_\idot(E_{\idot-i}))) \to
\bCP^{[p]}(E_\idot),
$$
and again by \eqref{phi.c}, the left-hand is exactly
$\bCP^f_\idot(E_\idot)$.
\endproof

Next, we note that the embedding $j_p:\Delta^o \to \Lambda_p$
extends to an embedding $\wj_p:\Delta^o \times \ppt_p \to
\Lambda_p$, where $\ppt_p$ is the groupoid with one object with
automorphism group $\Z/p\Z$. Therefore the pullback functor $j_p^*$
can be refined to a functor
\begin{equation}\label{wj.p}
\wj_p^*:\Fun(\Lambda_p,R) \to \Fun(\Delta^o \times \ppt_p,R) \cong
\Fun(\Delta^o,R[\Z/p\Z]).
\end{equation}
Moreover, if we denote by $\pi_p:\Delta^o \times \ppt_p \to
\Delta^o$ the projection onto the first component, then we have a
base change isomorphism $\pi_{p!}\circ\wj_p^* \cong j^* \circ
\pi_{p!}$. For any complex $E_\idot$ in the category
$\Fun(\Delta^o,R[\Z/p])$, we then denote by $\K_\idot(E_\idot)$ the
product $E_\idot \otimes \wj_p^*K_\idot$, and we define
$\pi_{p\flat}E_\idot$, $\vpi_{p\flat}E_\idot$,
$\vbpi_{p\flat}E_\idot$ by \eqref{c.p.eq} and \eqref{vpi}. We let
$cph_\idot(E_\idot) = CH_\idot(\pi_{p\flat}E_\idot) =
\per(CH_\idot(\pi_{p!}\K_\idot(E_\idot)))$ and
\begin{equation}\label{CPH.f}
CPH^f_\idot(E_\idot) = CH_\idot(\vpi_{p\flat}E_\idot), \quad
\bCPH^f_\idot(E_\idot) = CH_\idot(\vbpi_{p\flat}E_\idot).
\end{equation}
We have $\pi_{p\flat}\circ\wj_p^* \cong j^* \circ \pi_{p\flat}$,
$\vpi_{p\flat}\circ\wj_p^* \cong j^* \circ \vpi_{p\flat}$,
$\vbpi_{p\flat}\circ\wj_p^* \cong j^* \circ \vbpi_{p\flat}$, so that
for any complex $E_\idot$ in $\Fun(\Lambda_p,R)$, we have a natural
map
\begin{equation}\label{cph.j}
cph_\idot(\wj_p^*E_\idot) \to cph_\idot(E_\idot)
\end{equation}
and natural maps
\begin{equation}\label{CPH.f.j}
J^f:CPH_\idot^f(\wj_p^*E_\idot) \to CPH_\idot(E_\idot), \ 
\bJ^f:\bCPH_\idot^f(\wj_p^*E_\idot) \to \bCPH_\idot(E_\idot).
\end{equation}
Moreover, we let
$CPH_\idot(E_\idot)=\Per(CH_\idot(\pi_{p!}\K_\idot(E_\idot)))$, and
we denote by $CPH^{[p]}_\idot(E_\idot)$ the completion of the
filtered complex $cph^f(E_\idot^{[p]})$, where $E^{[p]}_\idot$ is
$E_\idot$ equipped with the $p$-th rescaling of the stupid
filtration, and $cph^f(-)$ is the filtered extension \eqref{phi.f}
of the functor $cph_\idot(-)$. Then for any complex $E_\idot$ in
$\Fun(\Lambda_p,R)$, the map \eqref{cph.j} induces natural maps
\begin{equation}\label{CPH.j}
J:CPH_\idot(\wj_p^*E_\idot) \to CPH_\idot(E_\idot), \ 
\bJ:\bCPH_\idot^{[p]}(\wj_p^*E_\idot) \to \bCPH^{[p]}_\idot(E_\idot).
\end{equation}
The same argument as in Lemma~\ref{L.v.le} shows that for any
complex $\wt{E}_\idot$ in $\Fun(\Delta^o,R[\Z/p\Z])$, we have
natural maps
\begin{equation}\label{R.L.j}
R:CPH^f_\idot(\wt{E}_\idot) \to CPH_\idot(\wt{E}_\idot), \quad
L^{[p]}:\bCPH^f_\idot(\wt{E}_\idot) \to \bCPH^{[p]}_\idot(\wt{E}_\idot),
\end{equation}
such that for any complex $E_\idot$ in the category
$\Fun(\Lambda_p,R)$, we have $J \circ R = R \circ J^f$ and $\bJ
\circ L^{[p]} = L^{[p]} \circ \bJ^f$.

\begin{defn}\label{fin.j.def}
A complex $E_\idot$ in the category $\Fun(\Delta^o,R[\Z/p\Z])$ is
{\em convergent} if the maps \eqref{R.L.j} are quasiisomorphisms.
\end{defn}

\begin{lemma}\label{L.R.le}
Assume given a complex $E_\idot$ in the category $\Fun(\Lambda_p,R)$
such that $\wj_p^*E_\idot$ is convergent in the sense of
Definition~\ref{fin.j.def}. Then the maps
$$
R:CP^f_\idot(E_\idot) \to CP_\idot(E_\idot), \qquad
L^{[p]}:\bCP^f(E_\idot) \to \bCP^{[p]}(E_\idot)
$$
are quasiisomorphisms.
\end{lemma}

\proof{} By Lemma~\ref{proj.le} and its completed versions, it
suffices to prove that the maps $J$, $\bJ$, $J^f$, $\bJ^f$ are
quasiisomorphisms: then the maps $R$, $L^{[p]}$ are retracts of the maps
\eqref{R.L.j}, and we are done. For the maps $J$, $J^f$ and $\bJ^f$
this is clear, since they are all versions of the quasiisomorphism
\eqref{ch.ch}. The case of the map $\bJ$ is more delicate, since the
definition of $CPH_\idot^{[p]}(E_\idot)$ and
$CPH_\idot^{[p]}(\wj_p^*E_\idot)$ involves completions, and the map
\eqref{ch.ch} is not a filtered quasiisomorphism with respect to the
standard filrations. However, we can also filter $cph_\idot(E_\idot)
= CH_\idot(\pi_{p\flat}E_\idot \otimes \K_\idot)$ by setting
$$
V^ncph_\idot(E_\idot) = cc_\idot(V^nE^\flat_\idot \otimes
\K_\idot),
$$
where the filtration $V^n$ on $E^\flat_\idot = \pi_{p\flat}E_\idot$
is as in Subsection~\ref{loc.subs}. Then $V^\hdot$ is conmensurable
with the standard filtration in every cohomogical degree, hence
gives the same completion $\bCP^{[p]}_\idot(E_\idot)$, and the induced
filtration on $cph_\idot(\wj_p^*E_\idot)$ has completion
$\bCPH^{[p]}_\idot(\wj_p^*E_\idot)$. The map $\bJ$ is now another
instance of the quasiisomorphism \eqref{ch.ch}.
\endproof

\subsection{Convergent complexes.}\label{conve.subs}

By virtue of Lemma~\ref{L.R.le}, studying the comparison maps $R$,
$L^{[p]}$ reduces to studying convergent complexes in the category
$\Fun(\Delta^o,R[\Z/p\Z])$. We will need two results in this
direction.

For any bicomplex $E_{\idot,\idot}$ of $R[\Z/p\Z]$-modules, consider
the Tate complex $\wC_\idot(\Z/p\Z,E_{\idot,\idot})$ of
\eqref{tate.bis}, and equip it with filtrations $F^\hdot$,
$F_1^\hdot$, $F_2^\hdot$ by
\begin{equation}\label{F.12}
\begin{aligned}
F^n\wC_\idot(\Z/p\Z,E_{\idot,\idot}) &= \bigoplus_{i \geq
  n}\wC_\idot(\Z/p\Z,E_{i,\idot}) \subset
\wC_\idot(\Z/p\Z,E_{\idot,\idot})\\
F_1^n\wC_\idot(\Z/p\Z,E_{\idot,\idot}) &= \bigoplus_{i \geq
  n}\wC_i(\Z/p\Z,E_{\idot,\idot}) \subset
\wC_\idot(\Z/p\Z,E_{\idot,\idot})\\
F_2^n\wC_\idot(\Z/p\Z,E_{\idot,\idot}) &= \bigoplus_{i + pj \geq
  n}\wC_\idot(\Z/p\Z,E_{i,j}) \subset
\wC_\idot(\Z/p\Z,E_{\idot,\idot}).
\end{aligned}
\end{equation}
Denote the completions of $\wC_\idot(\Z/p\Z,E_{\idot,\idot})$ with
respect to filtrations $F^\hdot_1$, $F^\hdot_2$ by
$C^1_\idot(E_{\idot,\idot})$ and
$C^2_\idot(E_{\idot,\idot})$. Then $F^\hdot$ induces a
filtration on each of these two complexes, and we have natural maps
\begin{equation}\label{f.loc}
\lim_{\overset{n}{\to}}F^{-n}C^1_\idot(E_{\idot,\idot}) \to
C^1_\idot(E_{\idot,\idot}), \qquad
\lim_{\overset{n}{\to}}F^{-n}C^2_\idot(E_{\idot,\idot}) \to
C^2_\idot(E_{\idot,\idot}). 
\end{equation}
If $E_{\idot,\idot} = CH_\idot(E_\idot)$ for a complex $E_\idot$ in
the category $\Fun(\Delta^o,R[\Z/p\Z])$, then we have
$\wt{C}_\idot(E_{\idot,\idot}) \cong cph_\idot(E_\idot)$. We then
have
$$
C^1_\idot(E_{\idot,\idot}) \cong CPH_\idot(E_\idot), \qquad
C^2_\idot(E_{\idot,\idot}) \cong \bCPH^{[p]}_\idot(E_\idot),
$$
and the maps \eqref{f.loc} are exactly the maps
\eqref{R.L.j}.

\begin{lemma}\label{hmo.le}
Assume given a complex $E_\idot$ in $\Fun(\Delta^o,R[\Z/p\Z])$, and
consider the complex $CH_\idot(E_\idot)$ as a complex in
$C_\idot(R[\Z/p\Z])$, with terms $E_\idot([n])$, $[n] \in
\Delta$. If $CH_\idot(E_\idot)$ is strongly bounded from below in
the sense of Definition~\ref{mix.qis}, then $E_\idot$ is convergent
in the sense of Definition~\ref{fin.j.def}.
\end{lemma}

\proof{} It suffices to prove that more generally, for any bicomplex
$E_{\idot,\idot}$ of $R[\Z/p\Z]$-modules that is strongly bounded
from below with respect to the first index, the maps \eqref{f.loc}
are quasiisomorphisms. This is obvious: if $E_{\idot,\idot}$ is
contractible with respect to the first index, then all the complexes
in \eqref{f.loc} are contractible, and if $E_{i,\idot} = 0$ for $i
\gg 0$, the limits in \eqref{f.loc} stabilize at a finite step.
\endproof

Now note that for any bicomplex $E_{\idot,\idot}$ of
$R[\Z/p\Z]$-bimodules, the filtration $F^\hdot_1$ of \eqref{F.12}
induces a spectral sequence, and if $E_{\idot,\idot} =
CH_\idot(E_\idot)$ for some complex $E_\idot$ in
$\Fun(\Delta^o,R[\Z/p\Z])$, this spectral sequences reads as
\begin{equation}\label{sp.d.1}
CH_\idot(E_\idot)((u))\langle\eps\rangle \Rightarrow CPH_\idot(E_\idot),
\end{equation}
where $\eps$ is a formal generator of cohomological degree $1$, as
in \eqref{z.pz.coh}. This is an analog of the standard spectral
sequence of a cyclic object. To obtain an analog of the conjugate
spectral sequence, we need a version of Definition~\ref{tight.la}.

\begin{defn}\label{tight.de}
A complex $E_\idot$ in $\Fun(\Delta^o,R[\Z/p\Z])$ is {\em tight} if
$E_\idot([n])$ is a tight complex of $R[\Z/p\Z]$-modules for any
$[n] \in \Delta$. For any tight complex $E_\idot$ in
$\Fun(\Delta^o,R[\Z/p\Z])$, the complex $\I(E_\idot)$ is given by
$$
\I(E_\idot)([n]) = \I(E_\idot([n]))
$$
for any $[n] \in \Delta$.
\end{defn}

Then for any tight complex $E_\idot$ in $\Fun(\Delta^o,R[\Z/p\Z])$,
we can consider the filtration on $cph_\idot(E_\idot)$ induced by
the filtration $\tau^\hdot$ on $\pi_{p\flat}E_\idot^{[p]}$. By
Lemma~\ref{tau.le}, in every cohomological degree, this filtration
is commensurable to the filtration $F^\hdot_2$ of \eqref{F.12}, so
that we obtain a spectral sequence
\begin{equation}\label{sp.d.2}
CH_\idot(\I(E_\idot))((u^{-1}))\langle\eps\rangle \Rightarrow
\bCPH^{[p]}_\idot(E_\idot),
\end{equation}
with the same meaning of $\eps$ as in \eqref{sp.d.1}.

\begin{lemma}\label{conve.le}
\begin{enumerate}
\item A finite extension of convergent complexes is convergent.
\item Assume given three tight complexes $E_\idot$, $E'_\idot$,
  $E''_\idot$ in $\Fun(\Delta^o,R[\Z/p\Z])$, and two maps
  $a:E'_\idot \to E_\idot$, $b:E_\idot \to E''_\idot$ such that both
  $b \circ a$ and $\I(b \circ a)$ are quasiisomorphisms. Then if
  $E_\idot$ is convergent, $E'_\idot$ and $E''_\idot$ are convergent
  as well.
\end{enumerate}
\end{lemma}

\proof{} \thetag{i} is clear. For \thetag{ii}, note that by
Lemma~\ref{tate.compl.le}, the maps
$$
\vpi_{p\flat}E'_\idot \to \vpi_{p\flat}E''_\idot, \qquad
\vbpi_{p\flat}E'_\idot \to \vbpi_{p\flat}E''_\idot
$$
induced by $b \circ a$ are quasiisomorphisms (for the second map,
note that the stupid filtrations on $E'_\idot$, $E''_\idot$ are
commensurable to their $p$-th rescaling, and by Lemma~\ref{tau.le},
the induced filtrations on $\pi_{p\flat}E'_\idot$,
$\pi_{p\flat}E''_\idot$ are in turn commensurable to the filtration
$\tau^\hdot$ in every cohomological degree). Therefore by
\eqref{CPH.f}, the corresponding maps
$$
CPH^f_\idot(E'_\idot) \to CPH^f_\idot(E''_\idot), \qquad
\bCPH^f_\idot(E'_\idot) \to \bCPH^f_\idot(E''_\idot)
$$
are also quasiisomorphisms. Moreover, we have convergent spectral
sequences \eqref{sp.d.1}, \eqref{sp.d.2}, so that the maps
$$
CPH_\idot(E'_\idot) \to CPH_\idot(E''_\idot), \qquad
\bCPH^{[p]}_\idot(E'_\idot) \to \bCPH^{[p]}_\idot(E''_\idot)
$$
induced by $b \circ a$ are also quasiisomorphisms. Thus the maps
\eqref{R.L.j} for $E'_\idot$ and $E'_\idot$ are retracts of the
corresponding maps for $E_\idot$, and we are done.
\endproof

\subsection{Characteristic $2$.}

In Subsection~\ref{edge.filt.subs} and throughout
Section~\ref{conj.sec}, we have assumed that our base ring $R$ is
annihilated by an odd prime $p$. Let us now describe what happens if
$p=2$.

Note right away that since we are in characteristic $2$, signs do
not matter. In particular, there is no difference between the
twisted complexes \eqref{tate.tw} and their untwisted versions.

The first problem occurs in the proof of
Proposition~\ref{edge.prop}, and specifically, in
\eqref{gr.ab.n}. If $p=2$, then $(p-1)n = n$ is no longer
necessarily divisible by $2$, and it can happen that
$\gr^{i-np}_F\nu_p$ is a map from the Tate complex to its odd
homological shift. However, this only happens if $n$ is odd. For any
$n \geq 1$ and any $R[\Z/2n\Z]$-module $E$, we have a functorial map
of complexes
\begin{equation}\label{eps.2}
\eps_\idot:C_\idot(\Z/2n\Z,E) \to C_{\idot+1}(\Z/2n\Z,E)
\end{equation}
given by
$$
\eps_i = \begin{cases}
\id + \sigma^2 + \dots + \sigma^{2(n-1)}, &\quad i = 2j,\\
\id, &\quad i=2j+1,
\end{cases}
$$
where $\sigma \in \Z/2n\Z$ is the generator. If $n$ is odd, the
map $\eps_\idot$ is a quasiisomorphism. Therefore for odd $n$, we
redefine the constants $a_n$, $b_n$ by setting
$$
\gr^{i-np}_F\nu_p \circ \eps_\idot = f(a_n,b_n)
$$
instead of \eqref{gr.ab.n}, and again, it suffices to prove that all
the maps $f(a_n,b_n)$ are quasiisomorphisms. Since in the case $E=k$
we obviously have $\eps_\idot = \id$ for all odd $n$, the rest of
the proof of Proposition~\ref{edge.prop} goes through without any
changes. The other results in Subsection~\ref{edge.filt.subs} also
go through with exactly the same proofs.

The real problem occurs in Section~\ref{conj.sec} and affects
Subsection~\ref{conj.subs} and Subsection~\ref{loc.subs} (everything
in Subsection~\ref{comp.subs} and Subsection~\ref{conve.subs} works
for any prime). If $p=2$, then the cohomology $H^\hdot(\Z/p\Z,k)$ is
still given by \eqref{z.pz.coh}, but the multiplication is
different: instead of $\eps^2=0$, we have $\eps^2=u$, the
periodicity generator. In fact, the action of the generator $\eps$
on the Tate complex is explicitly given by the quasiisomorphism
\eqref{eps.2}. Thus in particular, $\eps_{odd}$ is always an
isomorphism, so that every $R[\Z/2\Z]$-module $E$ is automatically
tight. But on the downside, $\eps_{even}$ is never equal to $0$ --
conversely, it is also an invertible map. Thus Lemma~\ref{I.le},
Lemma~\ref{V.CC.le} and Proposition~\ref{conj.prop} completely break
down.

To alleviate the situation, let us prove a weaker version of
Lemma~\ref{V.CC.le} that works for all primes.

Assume that the base ring $R$ is annhilated by a prime $p$, and
assume given a complex $E_\idot$ in $\Fun(\Lambda_p,R)$. Assume
further that $E_\idot$ is tight in the sense of
Definition~\ref{tight.def} (if $p=2$, this just means that for any
$[n] \in \Lambda_p$ and any odd integer $i$, $E_i([n])$ is a free
$R[\Z/p\Z]$-module). As in Subsection~\ref{loc.subs}, denote
$E^\flat_\idot = \pi_{p\flat}E_\idot^{[p]}$, $\wt{E}_\idot = E_\idot
\otimes_k \wt{\K}_\idot$, $\wt{E}^\flat_\idot =
\pi_{p\flat}\wt{E}_\idot$, and filter $CC_\idot(E^\flat_\idot)$ and
$CC_\idot(\wt{E}^\flat_\idot)$ by setting
\begin{equation}\label{v.cc.p}
V^\hdot CC_\idot(E^\flat_\idot) =
CC_\idot(V^\hdot E^\flat_\idot), \qquad
V^\hdot CC_\idot(\wt{E}^\flat_\idot) =
CC_\idot(V^\hdot \wt{E}^\flat_\idot),
\end{equation}
where the filtration $V^\hdot$ in the right-hand side is the
filtration \eqref{V.bV.loc}. Denote by $\wCC_\idot(E^\flat_\idot)$
the completion of the complex $CC_\idot(E^\flat_\idot)$ with respect
to the filtration \eqref{v.cc.p}, and analogously for
$\wt{E}_\idot$. Then as in Subsection~\ref{loc.subs}, the map
\eqref{aug.cc} induces a filtered map
$$
\langle CC_\idot(E^\flat_\idot),V^\hdot \rangle \to \langle
cp_\idot(E_\idot),V^\hdot \rangle,
$$
where $V^\hdot$ in the right-hand side is the conjugate filtration
of Definition~\ref{conj.def}. Passing to the completions, we obtain
a functorial map
\begin{equation}\label{wcc.cp}
\wCC_\idot(E^\flat_\idot) \to \bCP^{[p]}_\idot(E_\idot).
\end{equation}

\begin{lemma}\label{V.CC.le.bis}
The map \eqref{wcc.cp} is a quasiisomorphism.
\end{lemma}

\proof{} If we equip the complex $E^\flat_\idot \otimes \K_\idot$
with the shift $V^\hdot_{-1}E^\flat_\idot \otimes \K_\idot$ of the
filtration $V^\hdot$, then the map $\delta^\flat$ of
\eqref{delta.flat} is tautologically a filtered map, and the
composition $\gamma^\flat \circ \delta^\flat$ factors as
$$
\begin{CD}
\langle \wt{E}^\flat_\idot,V^\hdot \rangle @>{\iota}>> \langle
E^\flat_\idot,V^\hdot \rangle @>{\id}>> \langle
E^\flat_\idot,V^\hdot_{-1}\rangle,
\end{CD}
$$
where $\iota$ is a filtered quasiisomorphism, and $\id$ is the
identity map considered as a filtered map from $V^\hdot$ to
$V^\hdot_{-1}$. Since all shifts of a given filtration are
commensurable, the map
$$
\gamma \circ \delta:\wCC_\idot(\wt{E}^\flat_\idot) \to
\wCC_\idot(E^\flat_\idot)
$$
induced by the composition $\gamma^\flat \circ \delta^\flat$ is a
quasiisomorphism. Then as in the proof of Lemma~\ref{V.CC.le}, the
map
$$
\wCC_\idot(E^\flat_\idot \otimes \K_\idot) \to
\bCPH^{[p]}_\idot(E_\idot)
$$
is a quasiisomorphism by Lemma~\ref{c.cc.le}, and the map
\eqref{wcc.cp} is a retract of this map.
\endproof

As a corollary of Lemma~\ref{V.CC.le.bis}, we obtain a functorial
spectral sequence
\begin{equation}\label{conj.2}
HC_\idot(\gr^\hdot_VE^\flat_\idot) \Rightarrow
\bHP^{[p]}_\idot(E_\idot)
\end{equation}
induced by the filtration \eqref{v.cc.p}. Unfortunately, if $p=2$,
the identification of its first term presents a problem. Namely, let
$a_{odd},a_{even} \in \Ext^2(\I(E_\idot),\I(E_\idot))$ be the
extension classes given by successive associated graded quotients of
the filtration $\tau^\hdot$ on $2$-the periodic complex
$E^\flat_\idot$, so that $\gr^0_VE^\flat_\idot$ is an extension of
$\I(E_\idot)$ by $\I(E_\idot)$ given by the class $a_{even}$. If $p$
is odd, then Lemma~\ref{I.le} shows that $a_{even}=u \cdot id$,
where $u \in H^2(\Lambda,k)$ is the periodicity class. This implies
that $\gr^0_VE^\flat_\idot$ is isomorphic to $\I(E_\idot) \otimes
\K_\idot$ in the derived category $\D(\Lambda,R)$, so that the
left-hand side of \eqref{conj.2} is identified with
$HH_\idot(\I(E_\idot))$. If $p=2$, then it is easy to modify the
arguments of Lemma~\ref{I.le} to show that
$$
a_{even} + a_{odd} = u \cdot \id.
$$
Ideally, we would have $a_{odd}=0$, so that $a_{even} = u \cdot
\id$, as in the case of an odd prime. However, it certainly cannot
be true for any tight complex $E_\idot$. Indeed, we can just take a
single object $E \in \Fun(\Lambda_p,R)$, so that the tightness
condition becomes trivial, and consider a short exact sequence
$$
\begin{CD}
0 @>>> E' @>>> P @>>> E @>>> 0
\end{CD}
$$
with some projective object $P \in \Fun(\Lambda_p,R)$. Then
$\I(P)=0$ and $\I(E)=\I(E')[1]$, with the roles of the classes
$a_{odd}$ and $a_{even}$ for $E'$ played by $a_{even}$
resp.\ $a_{odd}$. Thus we cannot have $a_{odd}=0$ for $E$ and for
$E'$ at the same time.

We expect that for complexes $E_\idot$ that come from DG algebras
and DG categories, we do have $a_{odd}=0$, so that the spectral
sequence \eqref{conj.2} has the form prescribed by
Proposition~\ref{conj.prop}. However, we have not been able to prove
it. Therefore we restrict ourselves to the following
observation. Since the filtration $\tau^\hdot$ on $E^\flat_\idot$ is
commensurable to its rescaling $V^\hdot$, we can consider a spectral
sequence induced by $\tau^\hdot$ instead of $V^\hdot$. Thus for any
prime $p$ and tight complex $E_\idot$ in $\Fun(\Lambda_p,R)$, the
isomorphism \eqref{wcc.cp} trivially yields a spectral sequence
\begin{equation}\label{conj.2.bis}
HC_\idot(\I(E_\idot))((u^{-1}))\langle \eps \rangle \Rightarrow
\bHP^{[p]}_\idot(E_\idot)
\end{equation}
where $u$, $\eps$ are formal variables of cohomological degrees $2$
and $1$, as in \eqref{sp.d.1}.

\section{DG categories.}\label{alg.sec}

\subsection{The setup.}

We recall some basic facts about DG algebras and DG categories (the
standard reference is \cite{kel}).

By a {\em DG algebra} $A_\idot$ over a commutative ring $R$ we will
understand an associative unital differential-graded algebra over
$R$ considered up to a quasiisomorphism. In particular, we will
tacitly assume that all DG algebras are $h$-projective as complexes
of modules over $R$, and all individiual terms $A_i$, $i \in \Z$ are
flat $R$-modules (this can be achieved, for example, by choosing a
cofibrant representative with respect to the standard model
structure).

For any DG algebra $A_\idot$, we denote by $\D(A_\idot)$ the derived
category of left DG modules over $A_\idot$. Any DG algebra map
$f:A_\idot \to A'_\idot$ induces by pullback a natural triangulated
functor $f_*:\D(A'_\idot) \to \D(A_\idot)$, and it has a
left-adjoint $f^*:\D(A_\idot) \to \D(A'_\idot)$ sending $M_\idot$ to
$A'_\idot \otimes_{A_\idot} M_\idot$. The map $f$ is a {\em derived
  Morita equivalence} if $f^*$ and $f_*$ are an adjoint pair of
equivalences of categories. Every quasiisomorphism is a derived
Morita equivalence.

A {\em DG category} $A_\idot$ over $R$ is ``an algebra with several
objects'' --- we have a collection $S$ of objects, and a collection
of $\Hom$-complexes $A_\idot(s,s')$ for any two objects $s$, $s'$,
equipped with the associative composition maps and identity elements
$\id_s \in A_0(s,s)$ for any object $s$. If $S$ has exactly one
element, then $A_\idot$ is simply a DG algebra.  A DG category is
small if $S$ is a set. As in the DG algebra case, we consider small
DG categories up to a quasiequivalence, and we will tacitly assume
that the $\Hom$-complexes $A_\idot(-,-)$ in a small DG category
$A_\idot$ are $h$-projective complexes of flat $R$-modules. A {\em
  module} $M_\idot$ over a DG category $A_\idot$ is a contravariant
DG functor from $A_\idot$ to the DG category of complexes of
$R$-modules; explicitly, it is given by a collection of complexes
$M_\idot(s)$, $s \in S$, and structure maps $A_\idot(s',s) \otimes_R
M_\idot(s) \to M_\idot(s')$.  DG modules over a small DG category
$A_\idot$ form a triangulated derived category $\D(A_\idot)$, and
for any DG functor $f:A_\idot \to A'_\idot$ between small DG
categories, we have an adjoint pair of natural triangulated functors
$f^*:\D(A_\idot) \to \D(A'_\idot)$, $f_*:\D(A'_\idot) \to
\D(A_\idot)$. The functor $f$ is a {\em derived Morita equivalence}
if $f^*$, $f_*$ are adjoint equivalences of categories. More
generally, a sequence
\begin{equation}\label{loc.seq}
\begin{CD}
A'_\idot @>{f}>> A_\idot @>{g}>> A''_\idot
\end{CD}
\end{equation}
of DG categories and DG functors is a {\em localization sequence} if
$g^* \circ f^*=0$, $f^*$ and $g_*$ are fully faithful, and $\langle
f^*(\D(A'_\idot)),g_*(\D(A''_\idot))\rangle$ is a semiorthogonal
decomposition of the triangulated category $\D(A_\idot)$. For
example, one can take a small DG category
$A_\idot$ and a set $S$ of objects in $A_\idot$, let $A'_\idot$ be
the full subcategory $A_\idot(S) \subset A_\idot$ spanned by $S$,
and let $A''_\idot = A_\idot/A'_\idot$ be the {\em Drinfeld
  quotient} \cite{dr} obtained by formally adding to $A_\idot$
morphisms $h_s:s \to s$, $s \in S$ of homological degree $1$ such
that $d(h_s) = \id_s$. Then the sequence
\begin{equation}\label{dr.loc.seq}
\begin{CD}
A'_\idot @>>> A_\idot @>>> A_\idot/A'_\idot
\end{CD}
\end{equation}
is a localization sequence. This example is universal: every
localization sequence \eqref{loc.seq} is derived-Morita equivalent to
a sequence \eqref{dr.loc.seq}.

A DG module $M_\idot$ over a DG category $A_\idot$ is {\em perfect}
if it is compact as an object in $\D(A_\idot)$ (that is,
$\Hom(M_\idot,-)$ commutes with arbitraty sums). For any object $s
\in S$ of the category $A_\idot$, we have the representable module
$A_\idot^s$ given by
$$
A_\idot^s(s') = A_\idot(s',s), \qquad s,s' \in S.
$$
Representable modules are perfect, and any perfect module is a
retract of a finite extension of shifts of representable DG modules
(``retract'' here is understood in the derived category sense, that
is, $M_\idot$ is a retract of $M'_\idot$ if we have a module
$M''_\idot$ and maps $a:M''_\idot \to M'_\idot$, $b:M'_\idot \to
M_\idot$ such that $b \circ a$ is a quasiisomorphism). A small DG
category $A_\idot$ is {\em compactly generated} if there exists a
perfect DG module $M_\idot$ that weakly generates $\D(A_\idot)$ --
that is, for any $N_\idot \in \D(A_\idot)$, $\Hom(M_\idot,N_\idot) =
0$ implies $N_\idot=0$.

For any set $S$ of objects in a DG category $A_\idot$, we denote by
$A^S_\idot$ the complex given by
\begin{equation}\label{A.S}
A^S_\idot = \bigoplus_{s,s' \in S}A_\idot(s,s').
\end{equation}
If $S$ is finite, then this is a unital associative DG algebra over
$R$, and we have a natural equivalence of derived categories
$$
\D(A^S_\idot) \cong \D(A_\idot(S)),
$$
where $A_\idot(S) \subset A_\idot$ is the full DG subcategory
spanned by $S$. The category $A_\idot$ is compactly generated if and
only if the embedding $A_\idot(S) \to A_\idot$ is a derived Morita
eqivalence for a finite set $S$; in this case, $\D(A_\idot)$ is
equivalent to the derived category $\D(A^S_\idot)$ of the DG algebra
$A^S_\idot$.

For any two small DG categories $A_\idot$, $A'_\idot$ with sets of
objects $S$, $S'$, the {\em tensor product} $A_\idot \otimes_R
A'_\idot$ is the DG category with set of objects $S \times S'$ and
$\Hom$-complexes given by
$$
(A_\idot \otimes_R A'_\idot)(s_1 \times s_1',s_2 \times s_2') =
A_\idot(s_1,s_2) \otimes _R A'_\idot(s_1',s_2')
$$
for any $s_1,s_2 \in S$, $s_1',s_2' \in S'$. A {\em DG bimodule}
$M_\idot$ over $A_\idot$ is a DG-module over the product $A^o_\idot
\otimes_R A_\idot$ of $A_\idot$ with its opposite DG category
$A^o_\idot$. An example of a DG bimodule is the diagonal bimodule
$A_\idot$. A DG bimodule $M_\idot$ is perfect if it is perfect as a
DG module over $A^o_\idot \otimes_R A_\idot$, and a small DG
category $A_\idot$ is {\em smooth} if the diagonal bimodule
$A_\idot$ is perfect. Smoothness is a derived-Morita invariant
property: for any derived Morita equivalence $f:A_\idot \to
A'_\idot$ of small DG categories, $A_\idot$ is smooth if and only if
$A'_\idot$ is smooth. It is convenient to introduce another
Morita-invariant property of DG categories that does not seem to
have a standard name in the literature (although it did appear, for
example, in \cite{orlo}, where it was emphasized as an important
feature of DG categories of geometric origin).

\begin{defn}\label{dgcat.bnd}
A DG category $A_\idot$ over a ring $R$ is {\em bounded from above}
resp.\ {\em bounded from below} if for any objects $s,s' \in
A_\idot$ and any $R$-module $M$, the complex $A_\idot(s,s')
\otimes_R M$ is bounded from above resp.\ below in the sense of
Definition~\ref{strong.bnd}. A DG category is {\em bounded} if it is
bounded both from above and from below.
\end{defn}

Again, for any derived Morita equivalence $f:A_\idot \to A'_\idot$
of small DG categories, $A_\idot$ is bounded from above resp.\ below
if and only if $A'_\idot$ is bounded from above resp.\ below.

\medskip

For any DG algebra $A_\idot$, one defines a complex $A_\idot^\hs$ in
$\Fun(\Lambda,R)$ by setting
$$
A^\hs_\idot([n]) = A_\idot^{\otimes_R n}, \qquad [n] \in \Lambda,
$$
where terms in the product are numbered by vertices $v \in
V([n])$. For any map $f:[n'] \to [n]$, the corresponding map
$A_\idot^\hs(f)$ is given by
\begin{equation}\label{A.hs.f}
A^\hs_\idot(f) = \bigotimes_{v \in V([n])}m_{f^{-1}(v)},
\end{equation}
where for every finite totally ordered set $S$, we let
$m_S:A^{\otimes_R S}_\idot \to A_\idot$ be the multiplication
map. If $S$ is empty, we let $A^{\otimes_R S}_\idot = R$, and $m_S:R
\to A_\idot$ is the embedding map of the unity element $1 \in A_0$.

To extend it to DG categories, for any small DG category $A_\idot$
with the set of objects $S$, one sets
$$
A_\idot^\hs([n]) = \bigoplus_{s_1,\dots,s_n \in S}A_\idot(s_1,s_2)
\otimes_R \dots \otimes_R A_\idot(s_{n-1},s_n) \otimes_R
A_\idot(s_n,s_1)
$$
for any object $[n] \in \Lambda$. Then if $S$ is finite,
$A_\idot^\hs([n])$ is canonically a direct summand of
$(A^S_\idot)^\hs([n])$, where $A^S_\idot$ is the DG algebra
\eqref{A.S}, and the structure maps \eqref{A.hs.f} induce maps
between $A_\idot^\hs([n])$ turning $A^\hs_\idot$ into a complex in
$\Fun(\Lambda,R)$. In the general case, we let
$$
A^\hs_\idot = \lim_{\overset{S}{\to}}A_\idot(S)^\hs,
$$
where the limit is over all finite sets $S$ of objects in $A_\idot$.

The {\em Hochschild homology} $HH_\idot(A_\idot)$, the {\em cyclic
  homology} $HC_\idot(A_\idot)$, and the {\em periodic cyclic
  homology} $HC_\idot(A_\idot)$ of a small DG category $A_\idot$ is
defined by means of the complex $A^\hs_\idot$: we denote
$$
CH_\idot(A_\idot) = CH_\idot(A^\hs_\idot), \quad CC_\idot(A_\idot)
= CC_\idot(A^\hs_\idot), \quad CP_\idot(A_\idot) =
CP_\idot(A^\hs_\idot),
$$
and we let $HH_\idot(A_\idot)$, $HC_\idot(A_\idot)$,
$HP_\idot(A_\idot)$ be the homology groups of these complexes.  Any
DG functor $f:A_\idot \to A'_\idot$ between small DG categories
induces maps $HH_\idot(A_\idot) \to HH_\idot(A'_\idot)$,
$HC_\idot(A_\idot) \to HC_\idot(A'_\idot)$, $HP_\idot(A_\idot) \to
HP_\idot(A'_\idot)$. If $f$ is a derived Morita equivalence, then
all three maps are isomorphisms. More generally, a localization
sequence \eqref{loc.seq} induces a distinguished triangle
$$
\begin{CD}
CH_\idot(A'_\idot) @>>> CH_\idot(A_\idot) @>>> CH_\idot(A''_\idot)
@>>>
\end{CD}
$$
of Hochschild homology complexes and the induced long exact sequence
of Hochschild homology groups, and we have analogous distinguished
triangles and long exact sequences for cyclic and for periodic
cyclic homology. One shortens this by saying that $HH_\idot(-)$,
$HC_\idot(-)$, and $HP_\idot(-)$ are {\em additive invariants} of
small DG categories.

\subsection{Statements.}

Assume given a small DG category $A_\idot$ over a commutative ring
$R$, and consider the corresponding complex $A^\hs_\idot$ in
$\Fun(\Lambda,R)$.

\begin{defn}
The {\em co-periodic cyclic homology} $\bHP_\idot(A_\idot)$ is given
by
$$
\bHP_\idot(A_\idot) = \bHP_\idot(A_\idot^\hs),
$$
where the right-hand side is as in Definition~\ref{cp.def}, and the
{\em co-periodic cyclic complex} $\bCP_\idot(A_\idot)$ is the
complex $\bCP_\idot(A^\hs_\idot)$. The {\em polynomial periodic
  cyclic homology} $hp_\idot(A_\idot)$ is the homology of the
complex $cp_\idot(A_\idot^\hs)$ of Definition~\ref{cp.bis.def}, and
and the {\em restricted periodic and co-periodic cyclic homology}
$HP^f_\idot(A_\idot)$, $\bHP^f_\idot(A_\idot)$ are the homology of
the complexes $CP^f(A_\idot^\hs)$, $\bCP^f(A_\idot^\hs)$ of
Definition~\ref{restr.def}.
\end{defn}

\begin{lemma}\label{tors.alg.le}
For any small DG category $A_\idot$ over a commutative ring $R$, we
have
$$
hp_\idot(A_\idot) \otimes \Q = HP^f_\idot(A_\idot) \otimes \Q =
\bHP^f_\idot(A_\idot) \otimes \Q = \bHP_\idot(A_\idot) \otimes \Q =
0.
$$
\end{lemma}

\proof{} Lemma~\ref{cpf.Q.le} and Corollary~\ref{cpf.Q.corr}.
\endproof

Now we formulate our main technical result about co-periodic cyclic
homology. Assume that the base ring $R$ is annihilated by a prime
$p$. For any $R$-module $M$, denote by $M^{(1)}$ its Frobenius twist
--- that is, the module $V \otimes_R R^{(1)}$, where $R^{(1)}$ is
$R$ considered as an algebra over itself via the absolute Frobenius
map $R \to R$, $r \mapsto r^p$. For any DG category $A_\idot$ over
$R$, denote by $HH_\idot^{(1)}(A_\idot)$
resp.\ $HC^{(1)}_\idot(A_\idot)$ the homology groups of the
complexes $CH_\idot(A_\idot)^{(1)}$
resp.\ $CC_\idot(A_\idot)^{(1)}$.

\begin{prop}\label{conj.alg.prop}
Assume given a small DG category $A_\idot$ over a ring $R$
annihilated by a prime $p$. Then we have a convergent spectral
sequence
\begin{equation}\label{sp.alg.0}
HC^{(1)}_\idot(A_\idot)((u^{-1}))\langle \eps \rangle \Rightarrow
\bHP_\idot(A_\idot),
\end{equation}
where $u$, $\eps$ are formal variables of cohomological degrees $2$
and $1$, as in \eqref{sp.d.1}. Moreover, if $p \neq 2$, we have a
convergent spectral sequence
\begin{equation}\label{sp.alg}
HH^{(1)}_\idot(A_\idot)((u^{-1})) \Rightarrow \bHP_\idot(A_\idot).
\end{equation}
Both spectral sequences are functorial in $A_\idot$.
\end{prop}

The spectral sequence \eqref{sp.alg} is called the {\em conjugate
  spectral sequence} for the DG category $A_\idot$. We will prove
Proposition~\ref{conj.alg.prop} in
Subsection~\ref{alg.proof.subs}. For now, we use it to prove the
following result.

\begin{theorem}\label{morita.thm}
Assume that the commutative ring $R$ is Noetherian. Then for any
derived Morita equivalence $f:A_\idot \to A'_\idot$ between small DG
categories over $R$, the induced map
$$
\bHP_\idot(A_\idot) \to \bHP_\idot(A'_\idot)
$$
is an isomorphism.
\end{theorem}

\proof{} Since $R$ is Noetherian, the derived category of
$R$-modules is generated by residue fields of localizations of
$R$. Therefore is suffices to prove the claim after taking product
with such a residue field. In other words, we may assume right away
that $R=k$ is a field. Then if $k$ contains $\Q$, we are done by
Lemma~\ref{tors.alg.le}, and if not, we can apply the specral
sequence \eqref{sp.alg.0} and the corresponding property of the
cyclic homology functor $HC_\idot(-)$.
\endproof

Moreover, we can also prove a stronger statement --- not only is
co-periodic cyclic homology derived-Morita invariant, but it also
gives an additive invariant of small DG categories over a fixed
Noetherian ring.

\begin{theorem}\label{loc.thm}
Assume that the commutative ring $R$ is Noetherian. Then any
localization sequence \eqref{loc.seq} of small DG categories over
$R$ induces a long exact sequence
$$
\begin{CD}
\bHP_\idot(A'_\idot) @>>> \bHP_\idot(A_\idot) @>>> \bHP_\idot(A''_\idot)
@>>>
\end{CD}
$$
of co-periodic cyclic homology.
\end{theorem}

\proof{} By Theorem~\ref{morita.thm}, we may assume that the
localization sequence in question is of the form
\eqref{dr.loc.seq}. Then by definition, the composition $g \circ
f:A'_\idot \to A''_\idot = A_\idot/A'_\idot$ factors through natural
projection $q:A'_\idot \to A'_\idot/A'_\idot$ to the Drinfeld
quotient $A'_\idot/A'_\idot$. Therefore if we denote by
$\bCP'_\idot(A''_\idot)$ the cone of the natural map of complexes
$$
\begin{CD}
\bCP_\idot(A'_\idot) @>{f - q}>> \bCP_\idot(A_\idot)
\oplus \bCP_\idot(A'_\idot/A'_\idot),
\end{CD}
$$
then $g$ induces a natural map of complexes
\begin{equation}\label{del.eq}
\bCP'_\idot(A''_\idot) \to \bCP_\idot(A''_\idot),
\end{equation}
and since $\bCP_\idot(A'_\idot/A'_\idot)$ is acyclic by
Theorem~\ref{morita.thm}, it suffices to prove that \eqref{del.eq}
is a quasiisomorphism. As in the proof of Theorem~\ref{morita.thm},
it suffices to prove this when $R=k$ is a field, and the statement
then immediately follows from Lemma~\ref{tors.alg.le} if $k$
contains $\Q$ and from \eqref{sp.alg.0} otherwise.
\endproof

Finally, we also have a comparison result about different versions
of periodic cyclic homology (the proof is in
Subsection~\ref{alg.proof.subs}). Note that by definition, taking
$E_\idot = A_\idot^\hs$ in \eqref{5.dia} induces a corresponding
commutative diagram for any small DG category $A_\idot$.

\begin{theorem}\label{comp.thm}
Assume given a DG algebra $A_\idot$ over a Noetherian commutative
ring $R$.
\begin{enumerate}
\item If $A_\idot$ is bounded from above, then the map
  $r:hp_\idot(A_\idot) \to \bHP^f_\idot(A_\idot)$ of \eqref{5.dia}
  is an isomorphism.
\item If $A_\idot$ is bounded from below, then the map
  $l:hp_\idot(A_\idot) \to HP^f_\idot(A_\idot)$ of \eqref{5.dia} is
  an isomorphism.
\item Assume that $A_i=0$ unless $i \geq 0$. Then the map
  $R:\bHP^f_\idot(A_\idot) \to \bHP_\idot(A_\idot)$ of \eqref{5.dia}
  is an isomorphism.
\item Assume that $A_\idot$ is smooth. Then the map
  $R:\bHP^f_\idot(A_\idot) \to \bHP_\idot(A_\idot)$ of \eqref{5.dia}
  is an isomorphism, and the map $L:HP^f_\idot(A_\idot) \to
  HP_\idot(A_\idot)$ fits into a long exact sequence
$$
\begin{CD}
HP^f_\idot(A_\idot) @>{L}>> HP_\idot(A_\idot) @>>> HP_\idot(A_\idot)
\otimes \Q @>>>.
\end{CD}
$$
\end{enumerate}
\end{theorem}

Here DG algebra is treated as a DG category with one object, and
``bounded from above/below'' is understood in the sense of
Definition~\ref{dgcat.bnd}. The condition in \thetag{iii} can also
be replaced with its cohomological version (we can always replace
$A_\idot$ with its truncation $\tau^0 A_\idot$). Either way, the
condition is pretty strong. However, it does hold in some
interesting cases, for example when $A_\idot$ is a unital
associative algebra $A$ placed in homological degree $0$. In this
case, Theorem~\ref{comp.thm}~\thetag{i},\thetag{ii},\thetag{iii}
actually gives
$$
HP^f_\idot(A) \cong hp_\idot(A) \cong \bHP^f_\idot(A) \cong
\bHP_\idot(A).
$$
However, these groups are different from $HP_\idot(A)$ unless $A$
has finite homological dimension.

\begin{corr}
Assume given a small compactly generated DG category $A_\idot$ over
a Noetherian commutative ring $R$.
\begin{enumerate}
\item If $A_\idot$ is bounded from above and smooth, then the
  natural map $R \circ r:hp_\idot(A_\idot) \to \bHP_\idot(A_\idot)$
  is an isomorphism.
\item If $A_\idot$ is bounded from below and smooth, then the
  natural map $L \circ l$ fits into a long exact sequence
$$
\begin{CD}
hp_\idot(A_\idot) @>{L \circ l}>> HP_\idot(A_\idot) @>>> HP_\idot(A_\idot)
\otimes \Q @>>>.
\end{CD}
$$
\end{enumerate}
In particular, if $A_\idot$ is bounded and smooth, and $R \otimes \Q
= 0$, we have a natural isomorphism $HP_\idot(A_\idot) \cong
\bHP_\idot(A_\idot)$.
\end{corr}

\proof{} If the set $S$ of objects in $A_\idot$ is finite, all
claims for the category $A_\idot(S)$ immediately follow from the
Theorem~\ref{comp.thm} applied to the DG algebra $A^S_\idot$. In the
general case, finite subsets $S_0 \subset S$ such that $A_\idot(S_0)
\subset A_\idot$ is a derived Morita equivalence are cofinal in the
partially ordered set of all finite subsets in $S$. For every such
subset $S_0$, we know the claims for $A_\idot(S_0)$; to finish the
proof, notice that $hp_\idot(-)$ obviously commutes with filtered
direct limits, while $HP_\idot(-)$ and $\bHP_\idot(-)$ are
derived-Morita invariant.
\endproof

\subsection{Tensor powers.}\label{ten.subs}

Before we prove Theorem~\ref{comp.thm} and
Proposition~\ref{conj.alg.prop}, we need to make a digression about
tensor power functors.

Fix a prime $p$ and a commutative ring $R$ such that $pR=0$. For any
flat $R$-module $M$, consider the $p$-fold tensor power
$M^{\otimes_R p}$, and let the cyclic group $\Z/p\Z$ act on it by
permutations, with the generator $\sigma$ corresponding to the
order-$p$ permutation $\sigma:M^{\otimes_R p} \to M^{\otimes_R
  p}$. We then have a natural trace map
$$
\begin{CD}
M^{\otimes_R p} @>{\id + \sigma + \dots + \sigma^{p-1}}>>
M^{\otimes_R p}.
\end{CD}
$$
Denote by $C(M)$ the cokernel of this map, and consider the map
\begin{equation}\label{psi.eq}
\psi:M^{(1)} \to C(M)
\end{equation}
given by $\psi(\lambda \otimes m) = \lambda m^{\otimes p}$, $\lambda
\in R$, $m \in M$.

\begin{lemma}\label{ten.ti.le}
The $R[\Z/p\Z]$-module $M^{\otimes_R p}$ is tight. The map $\psi$ of
\eqref{psi.eq} is an additive $R$-linear map, and it induces an
isomorphism $M^{(1)} \cong \I(M)$.
\end{lemma}

\proof{} All the claims are obviously compatible with filtered
colimits and retracts, so that we may assume that $M = R[S]$ is the
free $R$-module generated by a set $S$. The rest of the proof
proceeds exactly as the special case when $R$ is a field considered
in \cite[Lemma 4.1]{ka1}. Namely, decompose
\begin{equation}\label{spl.m}
M^{\otimes_R p} = R[S^p] = M \oplus M',
\end{equation}
where $M' = R[S']$ is the free module spanned by complement $S' =
S^p \setminus S$ to the diagonal $S \subset S^p$. Then $M'$ is a
free $R[\Z/p\Z]$-module, and the $\Z/p\Z$-action on $M$ is trivial,
so that both are tight. Moreover, since $\sigma(m^{\otimes p}) =
m^{\otimes p}$ for any $m \in M$, the map $\psi$ takes values in
$\I(M^{\otimes_R p}) \subset C(M)$. Moreover, its composition with
the projection $\I(M^{\otimes_R p}) \to \I(M)$ onto the first
summand in \eqref{spl.m} is obviously an isomorphism. To finish the
proof, it remains to notice that since $M'$ is free over
$R[\Z/p\Z]$, we have $\I(M')=0$.
\endproof

Assume now given a complex $M_\idot$ of flat $R$-modules, and
consider the $p$-fold tensor power $M_\idot^{\otimes_R p}$. Again,
equip it with the permutation action of the group $\Z/p\Z$.

\begin{prop}\label{ten.ti.prop}
The complex $M^{\otimes_R p}_\idot$ is a tight complex of
$R[\Z/p\Z]$-mo\-du\-les in the sense of Definition~\ref{tight.def},
and we have a functorial isomorphism
\begin{equation}\label{ten.ti.eq}
\I(M^{\otimes_R p}) \cong M_\idot^{(1)}.
\end{equation}
\end{prop}

\proof{} As in Lemma~\ref{ten.ti.le}, we may assume right away that
$M_i$ is a free $R$-module for any integer $i$. For any integer $l$,
the degree-$l$ term of the complex $M^{\otimes_R p}_\idot$
decomposes as
$$
\left(M_\idot^{\otimes_R p}\right)_l = \bigoplus_{i_1+\dots+i_p=l}
M_{i_1} \otimes_R \dots \otimes_R M_{i_p},
$$
and the cyclic group $\Z/p\Z$ permutes the indices
$i_1,\dots,i_p$. If $l$ is not divisible by $p$, then the
$\Z/p\Z$-action on the set of indices such that $i_1+\dots+i_p=l$ is
stabilizer-free, so that $\left(M^{\otimes_R p}_\idot\right)_l$ is a free
$R[\Z/p\Z]$-module. If $l=np$ for some integer $n$, then we have
\begin{equation}\label{M.pr}
\left(M^{\otimes_R p}_\idot\right)_l = M_n^{\otimes_R p} \oplus M',
\end{equation}
where $M'$ is a free $R[\Z/p\Z]$-module. Since $M^{\otimes_R p}_n$
is tight by Lemma~\ref{ten.ti.le}, $M_\idot^{\otimes_R p}$ is
therefore indeed a tight complex of $R[\Z/p\Z]$-modules. Moreover,
\eqref{M.pr} together with \eqref{i.e.n} give an identification
$$
\I\left(M^{\otimes_R p}_\idot\right)_n \cong \I(M_n^{\otimes_R p})
\oplus \I(M') \cong \I(M^{\otimes_R p}_n),
$$
and Lemma~\ref{ten.ti.le} then provides a functorial isomorphism
\begin{equation}\label{wt.psi}
\wt{\psi}:M_\idot^{(1)} \cong \I(M_\idot^{\otimes_R p})
\end{equation}
of graded $R$-modules.

To check whether the map \eqref{wt.psi} commutes with the
differential, consider first the special case when $M_1=M_0=R$,
$M_i=0$ otherwise, and $d:M_1 \to M_0$ is the identity map. Then the
complex $M^{\otimes_R p}_\idot$ is strongly acyclic as a complex of
$R[\Z/p\Z]$-modules in the sense of Definition~\ref{strong.qis}, so
that $\wC_\idot(\Z/p\Z,M^{\otimes_R p})$ is an acyclic
complex. Equip it with the filtration $F^\hdot$ induced by the
$p$-th rescaling of the stupid filtration on $M_\idot^{\otimes_R
  p}$, and consider the spectral sequence induced by the filtration
$\tau^\hdot$ of \eqref{tau.beta.F}. Then for dimension reasons, it
degenerates at the term $E_3$, and the term $E_2$ is the sum of
shifts of the complex $\I(M^{\otimes_R p}_\idot)$. Thus the
differential in this complex must be an isomorphism, and we have
\begin{equation}\label{a.psi}
\wt{\psi} \circ d = a d \circ \wt{\psi}
\end{equation}
for some invertible element $a \in R$. But since the map
\eqref{wt.psi} is functorial, the same equality must then hold for
any complex $M_\idot$ of flat $R$-modules. Now to finish the proof,
it suffices to define a graded isomorphism
\begin{equation}\label{psi.comp}
\psi:M^{(1)} \cong \I(M_\idot^{\otimes_R p})
\end{equation}
by setting $\psi = a^{-i}\wt{\psi}$ in graded degree $i$, and
observe that $\psi$ commutes with the differentials by
\eqref{a.psi}.
\endproof

\begin{remark}
The construction of the map \eqref{psi.comp} is obviously functorial
in $R$, so that the constant $a$ in \eqref{a.psi} is actually an
invertible element in the prime field $k=\Z/p\Z$. If $p=2$, then of
course $a=1$, and if $p$ is odd, $a$ can be more explicitly
described as follows: take $M_\idot$ with $M_1=M_0=k$, $M_i=0$
otherwise, $d=\id$, and note that its $p$-th tensor power represents
by Yoneda a completely canonical class in the group
$\Exp^{p-1}_{k[\Z/p\Z]}(k,k) = H^{p-1}(\Z/p\Z,k)$. This class must
be of the form $au^{(p-1)/2}$, where $u \in H^2(\Z/p\Z,k)$ is the
natural generator, and $a \in k$ is some coefficient. This is our
element $a$. Alternatively, the Steenrod $p$-th power of the
generator $\eps \in H^1(S,k)$ of the first homology group of the
circle $S^1$ is equal to $a\eps$. In view of the latter description,
the value of $a$ must be in the literature, but I could not find it.
\end{remark}

Finally, we will need the following result about tensor powers of
$h$-projective complexes.

\begin{lemma}\label{ten.abs.le}
Assume given a quasiisomorphism $f:N_\idot \to M_\idot$ between
$h$-projective complexes of flat $R$-modules. Then the $p$-th tensor
power
$$
a^{\otimes p}:N_\idot^{\otimes_R p} \to M_\idot^{\otimes_R p}
$$
is a strong quasiisomorphism of complexes of $R[\Z/p\Z]$-modules in
the sense of Definition~\ref{strong.qis}.
\end{lemma}

\proof{} Since a quasiisomorphism of $h$-projective complexes is a
chain-ho\-mo\-topy equivalence, it suffices to prove that for two
chain-homotopic maps $f_1,f_2:N_\idot \to M_\idot$, the tensor
powers $f_1^{\otimes p}$, $f_2^{\otimes p}$ give the same maps in
the absolute derived category $\D_{abs}(R[\Z/p\Z])$. It obviously
suffices to consider the universal situation: we let $N_\idot =
M_\idot \otimes I_\idot$, where $I_\idot$ is the complex with terms
$I_0=\Z \oplus \Z$, $I_{-1}=\Z$, $I_i=0$ otherwise, $d:I_0 \to I_1$ the
difference map, and we let $f_1$, $f_2$ be the maps induced by the
projections $I_0 \to \Z$ onto the two summands. Then we also have the
map $e:M_\idot \to N_\idot$ induced by the diagonal embedding $\Z
\to I_0$, and $f_1 \circ e = f_2 \circ e = \id$. Therefore
$f_1^{\otimes p} \circ e^{\otimes p} = f_2^{\otimes p} \circ
e^{\otimes p} = \id$, and it suffices to prove that $e^{\otimes p}$
is a strong quasiisomorphism. This is clear: the diagonal
embedding $\Z \to I_\idot^{\otimes p}$ is a quasiisomorphism, and
since the complex $I_\idot^{\otimes p}$ sits in a finite range of
degrees, the quasiisomorphism must be strong.
\endproof

\begin{corr}\label{ten.bnd.corr}
Assume that an $h$-projective complex $M_\idot$ of $R$-modules is
bounded from above in the sense of Definition~\ref{strong.bnd}. Then
its $p$-th tensor power $M_\idot^{\otimes_R p}$ is strongly bounded
from above as a complex of $R[\Z/p\Z]$-modules.
\end{corr}

\proof{} Choose $n$ such that $\tau^nM_\idot \to M_\idot$ is a
quasiisomorphism, note that one can choose an $h$-projective
replacement for $\tau^nM_\idot$ that is trivial in homological
degrees $< n$, and apply Lemma~\ref{ten.abs.le}.
\endproof

\subsection{Proofs.}\label{alg.proof.subs}

We can now prove Proposition~\ref{conj.alg.prop} and
Theorem~\ref{comp.thm}. In fact, Proposition~\ref{conj.alg.prop}
immediately follows from \eqref{conj.2.bis},
Proposition~\ref{conj.prop} and the following result.

\begin{lemma}\label{i.a.la.le}
For any small DG category $A_\idot$ over a commutative ring $R$
annihilated by a prime $p$, the complex $A_\idot^\hs$ in
$\Fun(\Lambda,R)$ is $p$-adapted in the sense of
Definition~\ref{tight.la}, and we have a natural identification
\begin{equation}\label{i.ten}
\I(i_p^*A^\hs_\idot) \cong (A^\hs_\idot)^{(1)}.
\end{equation}
\end{lemma}

\proof{} If $A_\idot$ is a DG algebra --- that is, the set $S$ of
its objects consists of one element --- then for every $[n] \in
\Lambda_p$, we have
$$
A_\idot^\hs(i_p([n])) \cong \left(A_\idot^{\otimes_R
  n}\right)^{\otimes_R p}.
$$
Proposition~\ref{ten.ti.prop} immediately shows that $A^\hs_\idot$
is $p$-adapted, and the isomorphisms \eqref{ten.ti.eq} provide the
identification \eqref{i.ten}. If the set $S$ is finite, then for any
$[n] \in \Lambda_p$, $A_\idot^\hs(i_p([n]))$ is a retract of
$A^{S\hs}_\idot(i_p([n]))$. Since a retract of a tight complex is
obviously tight, $A^\hs_\idot$ is again $p$-adapted, and the
identification \eqref{i.ten} for $A_\idot$ is induced by the
corresponding identification for the DG algebra $A_\idot$. In the
general case, note that a filtered colimit of tight objects is
obviously tight, and the functor $\I(-)$ commutes with filtered
colimits.
\endproof

The proof of Theorem~\ref{comp.thm} --- specifically, of
Theorem~\ref{comp.thm}~\thetag{iv} --- requires some
preparations. Recall that the embedding $j:\Delta^o \to \Lambda$ of
\eqref{j.l} identifies $\Delta^o$ with the category of objects $[n]
\in \Lambda$ equipped with a distinguished vertex $v \in V([n])$. For
any bimodule $M$ over an associative unital flat $R$-algebra $A$, we
can define the simplicial $R$-modules $(M/A)^\hs \in
\Fun(\Delta^o,R)$ by setting
\begin{equation}\label{m/a}
(M/A)^\hs([n]) = A^{\otimes_R n-1} \otimes_R M, \qquad
[n] \in \Delta^o
\end{equation}
where the terms $A$ in the product correspond to vertices $v' \in
V(j([n]))$ different from the distinguished vertex $v$, and $M$
corresponds to $v$. The structure maps $(M/A)^\hs(f)$, $f:[n'] \to
[n]$ are given by \eqref{A.hs.f}. The complex
$$
CH_\idot(A,M) = CH_\idot((M/A)^\hs)
$$
is then the standard Hochschild homology complex of the group
algebra $A$ with coefficients in $M$. In particular, if $M = A^o
\otimes_R A$ is a free $A$-bimodule, the complex $CH_\idot(A,M)$ is
chain-homotopy equivalent to $HH_0(A,M) = A$ placed in degree
$0$. Moreover, for any integer $l \geq 1$, the $l$-fold tensor
product $A^{\otimes_R l}$ is an associative algebra equipped with an
action of the group $\Z/l\Z$ generated by the longet permutation
$\sigma:A^{\otimes_R l} \to A^{\otimes_R l}$. The $l$-fold tensor
product $M^{\otimes_R l}$ is naturally a bimodule over $A^{\otimes_R
  l}$. Denote by $M^{\otimes_R l}_\sigma$ the $R$-module
$M^{\otimes_R l}$ considered as an $A^{\otimes_R l}$-bimodule in the
following way: the left multiplication is the standard one, and the
right multiplication is twisted by $\sigma$ --- in other words, we
have
$$
a \cdot m \cdot a' = am\sigma(a'), \qquad a,a' \in A^{\otimes_R l},m
\in M^{\otimes_R l}.
$$
Then the object
\begin{equation}\label{ch.l}
(M/A)^\hs_l = (M^{\otimes_R l}_\sigma/A^{\otimes_R l})^\hs
\end{equation}
carries a natural action of $\Z/l\Z$ generated by the same
permutation $\sigma$ acting both on $A^{\otimes_R l}$ and on
$M^{\otimes_R l}$, so that is lies naturally in the category
$\Fun(\Delta^o,R[\Z/l\Z])$. We denote by
$$
CH_\idot^l(A,M) = CH_\idot((M/A)^\hs_l)
$$
the corresponding complex of $R[\Z/l\Z]$-modules, and we denote its
homology by $HH_\idot^l(A,M)$. For example, if $M = A^o \otimes_R
A$, we have
$$
HH_0^l(A,M) = A^{\otimes_R l},
$$
and the complex $CH_\idot^l(A,M)$ is chain-homotopy-equivalent to
its degree-$0$ homology placed in degree $0$.

For a bimodule $M_\idot$ over a DG algebra $A_\idot$, we can apply
\eqref{m/a} termwise and obtain the complex $(M_\idot/A_\idot)^\hs$
in the category $\Fun(\Delta^o,R)$. By definition, if $A_\idot$ is
the diagonal bimodule, we have $(A_\idot/A_\idot)^\hs \cong
j^*A^\hs_\idot$. Moreover, for any integer $l \geq 1$, we can define
the complex $(M_\idot/A_\idot)^\hs_l$ in $\Fun(\Delta^o,R[\Z/l\Z])$
by \eqref{ch.l}, and we have
\begin{equation}\label{ch.l.j}
(A_\idot/A_\idot)^\hs_l \cong \wj_l^*i_l^*A^\hs_\idot,
\end{equation}
where $\wj^*_l$ is the pullback functor \eqref{wj.p}.

\begin{lemma}\label{p.conv.le}
Assume that the commutative ring $R$ is annihilated by a prime
$p$. Then for any DG algebra $A_\idot$ over $R$ and any perfect
$A_\idot$-bimodule $M_\idot$ termwise flat over $R$, the complex
$(M_\idot/A_\idot)^\hs_p$ in $\Fun(\Delta^o,R[\Z/p\Z])$ is
convergent in the sense of Definition~\ref{fin.j.def}.
\end{lemma}

\proof{} Say that an $A_\idot$-bimodule is {\em finite free} if it
is a finite sum of homological shifts of the bimodule $A^o_\idot
\otimes_R A_\idot$, and say that it is {\em finite semifree} if it
admits a filtration whose associatied graded quotient is finite
free.  Since $M_\idot$ is perfect, we have a diagram
$$
\begin{CD}
M''_\idot @>{a}>> M'_\idot @>{b}>> M_\idot
\end{CD}
$$
of $A_\idot$-bimodules such that $M'_\idot$ is finite semifree and
$b \circ a$ is a quasiisomorphism. We can also choose it in such a
way that all the bimodules are complexes of flat $R$-modules. Then
the composition $b^{\otimes p} \circ a^{\otimes p}$ is also a
quasiisomorphism. Moreover, by the same argument as in
Lemma~\ref{i.a.la.le}, Proposition~\ref{ten.ti.prop} shows that for
any $A_\idot$-bimodule $N_\idot$ termwise flat over $R$, the complex
$(N_\idot/A_\idot)^\hs_p$ is tight, and we have
$$
\I((N_\idot/A_\idot)^\hs_p) \cong
\left((N_\idot/A_\idot)^\hs\right)^{(1)}.
$$
Therefore $I(b^{\otimes p}) \circ \I(a^{\otimes p})$ is also a
quasiisomorphism, and by Lemma~\ref{conve.le}~\thetag{ii}, it
suffices to prove the claim for the bimodule $M'_\idot$. In other
words, we may assume right away that the bimodule $M_\idot$ is
finite semifree.

Let $F^\hdot$ be the filtration on $M_\idot$ such that
$\gr^\hdot_FM_\idot$ is finite free. Then $F^\hdot$ induces a
$\Z/p\Z$-invariant filtration $F^\hdot$ on the tensor power
$M^{\otimes_R p}$, hence on $(M_\idot/A_\idot)^\hs_p$, and by
Lemma~\ref{conve.le}~\thetag{i}, it suffices to prove that
$\gr^i_F(M_\idot/A_\idot)^\hs_p$ is convergent for any integer
$i$. Since a direct summand of a convergent complex is obviously
convergent, it in fact suffices to prove that
$$
\gr^\hdot_F(M_\idot/A_\idot)^\hs_p \cong
(\gr^\hdot_FM_\idot/A_\idot)^\hs_p
$$
is convergent --- that is, we have to prove the claim for
$\gr^\hdot_F M_\idot$ instead of $M_\idot$.  In other words, we may
assume right away that $M_\idot$ is finite free.

In this case, the bimodule $CH_\idot^p(A_\idot,M_\idot)$ considered
as a complex in $C_\idot(R[\Z/p\Z])$ is chain-homotopy equivalent to
its degree-$0$ homology placed in degree $0$, so that
$(M_\idot/A_\idot)^\hs_p$ is convergent by Lemma~\ref{hmo.le}.
\endproof

\proof[Proof of Theorem~\ref{comp.thm}.] \thetag{iii} is obvious:
under the assumptions, the complex $A_\idot^\hs$ in
$\Fun(\Lambda,R)$ is trivial in negative homological degrees, so
that the map $L$ is actually an isomorphism. As in the proof of
Theorem~\ref{morita.thm}, when proving the other claims, we may
assume right away that $R$ is a field, and if $R$ contains $\Q$,
everything immediately follows from Lemma~\ref{tors.alg.le}. Thus we
may assume that $R$ is annihilated by a prime $p$. Then \thetag{i}
immediately follows from Corollary~\ref{bnd.la.corr}, and
\thetag{ii} immediately follows from Corollary~\ref{bnd.la.corr} and
Corollary~\ref{ten.bnd.corr}. For \thetag{iv}, apply
Lemma~\ref{p.conv.le}, \eqref{ch.l.j}, and Lemma~\ref{L.R.le}.
\endproof

\subsection{Comparison with de Rham cohomology.}\label{derham.subs}

Now assume given a unital associative flat algebra $A$ over a
commutative ring $k$, and assume further that $A$ is commutative. In
this setting, the classic theorem of Hochschild, Kostant and
Rosenberg asserts the following.

\begin{theorem}[Hochschild-Kostant-Rosenberg]
In the assumptions above, assume further that $A$ is finitely
generated, and $X = \Spec A$ is smooth over $k$. Then
for any integer $i$, we have a natural isomorphism
\begin{equation}\label{hkr}
HH_i(A) \cong \Omega^i(A),
\end{equation}
where $\Omega^i(A) = H^0(X,\Omega^i_X)$ is the space of global
differential $i$-forms on $X$.\endproof
\end{theorem}

This results requires no assumptions whatsoever on $k$ --- in
particular, it can be a field of positive characteristic, or a ring
such as $\Z$. Moreover, the Hochschild-to-cyclic spectral sequence
induces a differential $B:HH_\idot(A) \to HH_{\idot+1}(A)$, and in
the Hochschild-Kostant-Rosenberg situation, it coincides with de
Rham differential (see e.g.\ \cite[Theorem 2.2]{ka1} for a short
proof).

If $k$ contains $\Q$, then one can do much more: the isomorphisms
\eqref{hkr} on individual Hochschild homology groups can be lifted
to the whole Hochschild homology complex $CH_\idot(A)$, and they are
compatible with the Connes-Tsygan differential $B$. As the result,
one obtains a natural quasiisomorphism
\begin{equation}\label{hp.dr}
HP_\idot(A) \cong H^\hdot_{DR}(A)((u)),
\end{equation}
where as always, $u$ is a formal generator of cohomological degree
$2$.

In positive characteristic, the situation is not as good. The
standard isomorphisms \eqref{hkr} typically do not lift to
$CH_\idot(A)$. Somewhat surprisingly, if the algebra $A$ is a
polynomial algebra in several generators, then there is a different
system of isomorphisms that does lift to the level of complexes in a
way compatible with $B$ (see \cite[Remark 3.2.3]{Lo} and
\cite[Theorem 3.2.5]{Lo} --- in fact, on the level of homology
groups, these isomorphisms are inverse to those of
\eqref{hkr}). Therefore if $A = k[x_1,\dots,x_n]$, we do have the
identification \eqref{hp.dr}. However, the isomorphism depends not
only on the algebra $A$ but on the actual choice of the generators
$x_1,\dots,x_n \in A$, and it is not known whether it holds if we
only assume that $\Spec A$ is smooth.

If the algebra $A$ is still commutative but $\Spec A$ is no longer
smooth, then one can always find a simplicial resolution of $A$ ---
that is, a functor $A_\idot$ from $\Delta^o$ to smooth commutative
$k$-algebras such that the $H_i(\Delta^o,A_\idot)=0$ for $i \geq 1$,
and $H_0(\Delta^o,A_\idot)$ is identified with $A$ by means of an
algebra map $A_0 \to A$. In this case, one can take the de Rham
complex termwise, and obtain a complex $\Omega^\hdot(A_\idot)$ in
the category $\Fun(\Delta^o,k)$. Taking the standard complex of the
simplicial complex, as in \eqref{d.del}, we obtain a bicomplex that
we also denote $\Omega^\hdot(A_\idot)$ by abuse of notation. 

The {\em derived de Rham cohomology groups} $H^\hdot_{DR}(A)$ of
Illusie are then the cohomology groups of the product-total complex
$\Tot(\Omega^\hdot(A_\idot))$. They do not depend on the choice of a
resolution $A_\idot$, there is a comparison theorem relating them to
cristalline cohomology of $X = \Spec A$, and they are related to the
derived exterior powers of the cotangent complex $\Omega_\idot(A)$
by a spectral sequence of a Hodge-to-de Rham type.

One can also consider the cohomology groups
$\overline{H}^\hdot_{DR}(A)$ of the sum-total complex
$\tot(\Omega^\hdot(A_\idot))$ and call them, for example, {\em
  restricted derived de Rham cohomology groups}. Then it has been
shown by Bhatt \cite{Bh}, following Beilinson, that if the base ring
$k$ is an algebra over $\Z/p^n$ for some $n$, then
$\overline{H}^\hdot_{DR}(A)$ also do not depend on the resolution
$A_\idot$ and satisfy even better compatibility results with respect
to cristalline cohomology. Moreover, if $k$ is an algebra over the
prime field $\F_p$, then $\overline{H}^\hdot_{DR}(A)$ are related to
the derived exterior powers of $\Omega_\idot(A_\idot)$ by a spectral
sequence of the conjugate type.

It would be natural to expect that our theory is parallel to that of
Beilinson and Bhatt, in that we have a version of the isomorphism
\eqref{hp.dr} relating co-periodic cyclic homology and restricted de
Rham cohomology, and our conjugate spectral sequence \eqref{sp.alg}
recovers Bhatt's conjugate spectral sequence.

Unfortunately, at present, such a picture seems beyond our
reach. Indeed, one can always choose a resolution $A_\idot$ so that
its terms $A_i$ are polynomial algebras over $k$, and if $A$ is
finitely generated, then one can further arrange so that all the
$A_i$ are finitely generated. Therefore for each individual $A_i$,
we do have the isomorphism \eqref{hp.dr}. However, these
isomorphisms depend on the choice of generators, and therefore do
not necessarily patch together with respect to the face and
degeneracy maps in the simplicial algebra $A_\idot$.

One example where Beilinson and Bhatt's approach is particularly
successful --- in fact, the motivating example for their whole
theory -- is the ring of integers $\overline{Z}_p \subset
\overline{\Q}_p$ in the algebraic closure $\overline{\Q}_p$. This is
not annihilated by a power of $p$, so one has to modify the
definition slightly --- one considers the quotient algebras $A_n =
\overline{Z}_p/p^n$, $n \geq 1$, and one defines the $p$-adically
completed cohomology by the (derived) inverse limit
$$
\wh{H}_{DR}^\hdot(\Z_p) =
\dlim_{\overset{n}{\gets}}\overline{H}^\hdot_{DR}(A_n).
$$
In this case, Bhatt proves \cite[Proposition 9.9]{Bh} that
$\wh{H}^i_{DR}(\overline{Z}_p)$ vanishes for $i \neq 0$, and we have
$$
\wh{H}^0_{DR}(\overline{Z}_p) = A_{cris},
$$
the Fontaine's cristalline period ring. The usual derived de Rham
cohomology $H_{DR}^\hdot(\overline{Z}_p)$, similarly $p$-adically
completed, recovers the bigger de Rham period ring $A_{dR}$.

We would venture to suggest that for dimension reasons, in this
particular case $A=\overline{Z}_p/p^n$, an isomorphism
\eqref{hp.dr} ought to exist, and we have a natural identification
$$
\bHP_\idot(A) \cong \overline{H}^\hdot_{DR}(A)((u^{-1})).
$$
Then by Bhatt's result, the derived inverse limit of the groups
$\bHP_\idot(A_n)$ should be computable in terms of Fontaine's
ring. However, given our existing technology, we do not see any way
to prove it.

\bigskip

\noindent
{\sc
Steklov Math Institute, Algebraic Geometry section\\
\mbox{}\hspace{30mm}and\\
Laboratory of Algebraic Geometry, NRU HSE\\
\mbox{}\hspace{30mm}and\\
Center for Geometry and Physics, IBS, Pohang, Rep. of Korea
}

\medskip

\noindent
{\em E-mail address\/}: {\tt kaledin@mi.ras.ru}

\end{document}